\documentclass[12pt]{amsart}
\usepackage{graphicx}
\usepackage[all]{xy}
\usepackage{amssymb, amsmath}
\usepackage{setspace}
\usepackage{eufrak}
\usepackage{amsfonts}
\usepackage{float}
\usepackage{amsmath}

\usepackage[super]{nth}

\usepackage{geometry}
\usepackage{colortbl}
\usepackage{color}

\def\qed{\hfill $\Box$}
\def\proof{\noindent {\sl Proof} :\;  }
\def\t{\noindent}

\newcommand{\Mcal}{\mathcal{M}}

\newcommand{\Zcal}{\mathcal{Z}}
\newcommand{\A}{\mathcal{A}}
\newcommand{\K}{\mathcal{K}}
\newcommand{\Ical}{\mathcal{I}}

\newcommand{\Rcal}{\mathcal{R}}

\newcommand{\C}{\mathbb{C}}
\newcommand{\Z}{\mathbb{Z}}
\newcommand{\N}{\mathbb{N}}
\newcommand{\Q}{\mathbb{Q}}
\newcommand{\R}{\mathbb{R}}
\newcommand{\bA}{\mathbb{A}}

\newcommand{\Ost}{\mathcal{O}}
\newcommand{\bP}{\mathbb{P}}
\newcommand{\Proj}{\mathbb{P}}
\newcommand{\bL}{\mathbb{L}}

\newcommand\codim{ \mbox{\rm codim}\, }

\def\Hom{{\rm Hom}}
\def\qed{\hfill $\Box$}
\def\proof{\noindent {\sl Proof} :\;  }
\def\t{\noindent}
\def\rd{\partial}
\def\Hilb{{\rm Hilb}}
\def\mDelta{{\varDelta}}
\def\mm{\mathfrak{m}}

\def\SM{{\rm SM}}

\def\Spec{{\rm  Spec}}
\def\Sm{{\rm \bf Sm}_k}
\def\Sch{{\rm \bf Sch}_k}
\def\CH{{\rm CH}}
\def\Aut{{\rm Aut}}

\def\bv{\mbox{\boldmath $v$}}

\def\b0{\mbox{\boldmath $0$}}

\def\bbv{\mbox{\tiny$\bv$}}

\newtheorem{thm}{\bf Theorem}[section]
\newtheorem{"thm"}[thm]{\bf `Theorem'}
\newtheorem{cor}[thm]{\bf Corollary} 
\newtheorem{lem}[thm]{\bf Lemma} 
\newtheorem{prop}[thm]{\bf Proposition} 
\newtheorem{definition}[thm]{\bf Definition} 
 
\newtheorem{rem}[thm]{\bf Remark} 
\newtheorem{exam}[thm]{\bf Example}

\begin{document}

\title[Universal polynomials for multi-singularities]
{Universal polynomials for \\multi-singularity loci of maps
} 
\author[T.~Ohmoto]{Toru Ohmoto}
\address[T.~Ohmoto]{Department of Pure and Applied Mathematics, Graduate School of Fundamental Science and Engineering, Waseda University,
Tokyo 169-8555, Japan}
\email{toruohmoto@waseda.jp}
%
%
%
%
\maketitle
\begin{abstract}
In the present paper, we prove the existence of universal polynomials which express {\em multi-singularity loci classes}    
of prescribed types for proper morphisms $f: X \to Y$ between smooth schemes 
over an algebraically closed field of characteristic zero -- we call them {\em Thom polynomials for multi-singularity types of maps}. 
It has been referred to as the {\em Thom-Kazarian principle} and unsolved for a long time. 
This result solidifies the foundation for a general enumerative theory of singularities of maps which is applicable to a broad range of problems in classical and modern algebraic geometry. 
For instance, it extends curvilinear multiple-point theory to be adapted to maps with arbitrary singularities; 
it guarantees universal expressions for counting divisors 
with prescribed isolated singular points in a given family of divisors; 
various counting problems of intersection and contact are formulated in a united way, 
and so on. 
In particular, it would contribute to a satisfactory answer to the rest of (an advanced form of) Hilbert's 15th problem and connect such classics to recent new interests in enumerations inspired by mathematical physics and other fields. 
First, we introduce the multi-singularity loci class for an arbitrary proper morphism via Intersection Theory on Hilbert schemes of points combined with the Thom-Mather Theory of singularities of maps. A main feature of our proof is a striking use of algebro-geometric {\em cohomology operations}.  
Somewhat surprisingly, when trying to grasp a full perspective of classical enumerative geometry,  we will inevitably encounter algebraic cobordism and motivic cohomology. 
\end{abstract}
%
%

\setcounter{tocdepth}{2}
\tableofcontents


\section {Introduction}

From classical to modern enumerative geometry, numerous problems are related to counting multiple-points of singular maps associated to given geometric situations. For instance, counting hyperplanes tangent at $r$ points to a hypersurface in the $r$-dimensional projective space $\bP^r$ is a typically classical problem, while a modern question posed in string theory,  {\em instanton counting}, is concerning enumeration of $r$-nodal curves in a sub-linear system on a projective surface. 
Although these two questions seem quite different at first glance, they are actually very close to each other; it is because they can be interpreted as counting $r$ multiple-points of certain singular maps  (\S \ref{appl}). 
A little surprisingly, however, an adequate general theory for this kind of enumerative problems has not yet been well established so far -- that is the main theme of the present paper. 

Roughly saying, a {\em multi-singularity} of a map $f: X \to Y$ between manifolds means the germ of $f$ at distinct  points $p_1, \cdots, p_r$ of $X$ satisfying $f(p_1)=\cdots =f(p_r) \in Y$, while to distinguish,  a singular map-germ at a single point ($r=1$) is referred to as a {\em mono-singularity} or {\em local singularity}. 
Put $$\kappa:=\dim Y-\dim X$$
and we call it the {\em codimension of $f$}. 
A {\em singularity type} is an equivalence class of map-germs up to local isomorphisms, and  {\em  multi-singularity loci} of prescribed type for $f: X \to Y$ are defined in a certain way as subvarieties of $X$ and $Y$.  Precise definitions of those terms are provided in the later sections. 
Of our particular interest is to universally express the cohomology classes represented by those multi-singularity loci in terms of intrinsic invariants of $X$, $Y$ and $f$,  and also to reveal the whole hidden combinatorial structure among those cohomological expressions. 
A solid foundation for solving such problems has been intensively pursued over the past decades within the framework of Intersection Theory developed by Fulton, Kleiman and MacPherson \cite{Fulton, Kleiman87} -- in \S \ref{multiple}, we will briefly summarize the {\em multiple-point theory} and its limitation as our motivational background. 

Meanwhile, at the beginning of this century, a remarkably novel approach was proposed from   topology, that suggests a significant breakthrough in this field. 
Based on his topological argument using complex cobordism theory $MU^*$, M. Kazarian \cite{Kaz03, Kaz06} has stated the existence of a unique universal polynomial which represents the (source) multi-singularity locus of prescribed type for any generic proper holomorphic map $f: X \to Y$ between complex manifolds -- it is called the {\em Thom polynomial for a multi-singularity type of maps}. We will use this terminology to mean a universal polynomial in  variables $c_i \; (i \ge 1)$ and $s_I$ with rational coefficients such that its shape depends only on the singularity type, where the {\em quotient Chern class} $c_i=c_i(f)$ means the $i$-th Chern class of the difference bundle $f^*TY-TX$, i.e., 
$$c(f):=1+c_1(f)+c_2(f)+\cdots =\frac{1+f^*c_1(TY)+f^*c_2(TY)+\cdots}{1+c_1(TX)+c_2(TX)+\cdots}$$ 
and the {\em Landweber-Novikov class} $s_I=s_I(f)$ with index $I=(i_1, i_2, \cdots, i_k)$ is given by  
$$s_I(f):=f_*(c^I(f))=f_*(c_1(f)^{i_1}c_2(f)^{i_2}\cdots c_k(f)^{i_k}).$$
The information of the map $f$ is only encoded in the substitution for variables $c_i=c_i(f)$ and $f^*s_I$ of the polynomial, that infers the principle that 
{\em local symmetry of singularities governs the global appearance of the singularity locus for a given map $f$}. 
Table 1 is a sample list of Thom polynomials for a few typical singularity types in case of codimension $\kappa=1$. For instance, let $f: X^2 \to Y^3$ be an appropriately generic map from a projective surface into a $3$-fold,  and especially  $Y=\Proj^3$ for the simplicity, then the number of crosscap points ($A_1$), the number of triple-points ($A_0^3$) and the degree of the double-point curve ($A_0^2$) are precisely obtained by computing the corresponding Thom polynomials in Table 1 evaluated by $c_i(f)$ and $f^*s_I(f)$ for $c_i$ and $s_I$. Similarly, for a map $f: X^3 \to \Proj^4$ of a $3$-fold into $4$-space, not only the degree of the curve loci formed by singularities of types $A_1$, $A_0^2$ and $A_0^3$ but also the number of singular points of each type of $A_1A_0$, $A_0A_1$ and $A_0^4$ is counted in entirely the same way (the meaning of symbols will be explained in the latter section). 
This picture is drastically extended by Kazarian -- many new Thom polynomials are systematically provided and applied to classical/modern enumerative problems, that leads to rediscover many known formulas and to make hundreds of generalizations possible. A key observation, not apparent at all, indicates a well-organized combinatorial structure among all multi-singularity Thom polynomials. It would be regarded as a taxonomy or systematics for `counting formulas' in enumerative geometry.

\begin{table}
$$
\begin{array}{c  |  l}
\hline 
\mbox{\footnotesize type} &  \mbox{\footnotesize Thom polynomials}\\
\hline 
A_0^2 &  s_0-c_1 \\
A_1 &  c_2 \\
A_0^3 &  \frac{1}{2}(s_{0}^2-s_1-2s_0c_1+2c_1^2+2c_2)\\
A_0A_1 &  s_{01}-2c_1c_2-2c_3\\
A_1A_0 &   s_0c_2-2c_1c_2-2c_3 \\
A_0^4 &   
\frac{1}{3!}\left(
\begin{array}{l}
s_0^3-3s_0s_1+2s_2+2s_{01}-3s_0^2c_1+3s_1c_1\\
+6s_0c_1^2+6s_0c_2-6c_1^3-18c_1c_2-12c_3
\end{array}
\right)\\
\hline
\end{array}
$$
\caption{\small Samples of Thom polynomials ($\kappa=1$)}
\label{tp_codim1}
\end{table}

Actually, the simplest case is treated in a well-established theory, 
originated by Ren\'e Thom \cite{Thom}, which shows the existence of a universal polynomial of Chern classes $c_i(f)$ that represents the locus of a given {\em mono-singularity type} \cite{Porteous, Damon, Ronga73, Rimanyi01, FR02, FR12, Kaz03, Kaz06, Ohmoto06, Ohmoto08, BercziSzenes}. 
A typical example is the famous {\em Thom-Porteous formula}; 
a certain Schur polynomial in Chern classes expresses the degeneracy locus of a vector bundle morphism or the $\Sigma^i$-locus of a map, see \S \ref{TP}.  
In this theory, most essential is the existence of the {\em classifying stack for mono-singularity types}, which is given by the quotient stack $[V/G]$ for a certain algebraic group action of $G$ on some affine space $V$ 
(for the Thom-Porteous case, we consider the natural action of $G=GL(l) \times GL(m)$ on the affine space $V$ of $l\times m$ matrices). 
Then the Thom polynomial of a mono-singularity type is given by the $G$-equivariant characteristic class representing the orbit closure of that type. That is fully established in the framework of Equivariant Intersection Theory based on the algebraic Borel construction \cite{Totaro, EdidinGraham}. 

In contrast, for counting multiple-points of singular maps ($r\ge 2)$,  
certain moduli spaces of $r$-tuples of points on $X$ (e.g., Hilbert schemes, Fulton-MacPherson configuration spaces etc.) must be involved, that makes the problem much deeper and harder than the mono-singularity case. 
Critical is, at least for now, the lack of the {\em classifying stack for multi-singularity types}. 
In fact, there has not yet been found any full proof of the existence of multi-singularity Thom polynomials, so we call it {\em Kazarian's conjecture} or the {\em Thom-Kazarian principle} throughout the present paper. 
Here is an impressive phrase from an article of S. Kleiman on Hilbert's 15th problem  \cite{Kleiman76}: 
\begin{quote}
``{\it ... In the offing, there is the exciting hope of the development in algebraic geometry of a general enumerative theory of singularities of mappings, a theory of Thom polynomials, which will, among other things, unify and justify the classical work dealing with prescribed conditions of intersection and contact imposed on linear spaces.}"
\end{quote} 
Nearly half a century later, the story is not complete yet  -- from our viewpoint, Kazarian's conjecture is exactly the missing piece of the expected theory ! 

In the present paper, we prove the conjecture with a sufficient generality 
in the context of algebraic geometry over an algebraically closed field $k$ of characteristic zero. 
We deal with {\em arbitrary} proper morphisms $f: X \to Y$ between smooth schemes, i.e., without any restriction of singularities of $f$ and dimensions. 
A key ingredient of our argument is Intersection Theory on relative Hilbert schemes of points associated to $f$ and its relation with the Thom-Mather theory of singularities of mappings under (formal) algebraic geometry setting. 
Hilbert schemes are highly singular in general, that has long been a major obstacle to establishing a general multiple-point theory (\S \ref{multiple}). 
In order to overcome it, we employ an entirely new approach -- it heavily relies on algebraic cobordism theory $\Omega^*$ constructed by  Levine-Morel \cite{LevineMorel} and further simplified by Levine-Pandharipande \cite{LevinePandhari}. 
The theory $\Omega^*$ requires {\em resolution of singularities}, and therefore (just for this), we put the hypothesis of characteristic zero. 
Most notable in our approach is that we find a `shortcut proof' of the conjecture by making use of {\it cohomology operations} in algebraic geometry, which has recently been established by A. Vishik \cite{Vishik} as an algebraic counterpart to cohomology operations in complex cobordism $MU^*$.  

We will state precisely our main results in \S \ref{thms} later (Theorems \ref{main_thm1} and \ref{main_thm2}). 
Instead, in the steps (i) to (vi) below, we outline a very clear flow of our arguments and  results. Each step requires many technical details and ideas, which will be fully discussed in the main body from \S 2. There are two parts: the first part (i--iii) is devoted to defining multi-singularity loci classes for {\em arbitrary} proper morphisms, and the second part (iv--vi) is devoted to our proof of the existence of universal polynomial expressions for those classes. A major prototype is {\em degeneracy loci class} in the Thom-Porteous formula and {\em double-point loci class} in the double-point formula (see \S \ref{motivation}). 

Let $\kappa \in \Z$ be fixed. Assume that $X$ and $Y$ are smooth schemes with 
$$m=\dim X\;\;\;\; \mbox{and} \;\;\;\; l=m+\kappa =\dim Y$$ 
(letters $l, m, \kappa$ are used with this meaning throughout the present paper). 
Let $\Ost_{X,p}$ denote the stalk at a closed point $p \in X$ of the structure sheaf of $X$ and $\mathfrak{m}_{X,p}$ the maximal ideal. A minimal basics of the Thom-Mather Theory for singularities of maps will be reviewed in \S \ref{TM}.

\begin{itemize}
\item[(i)]  
Given a finitely determined mono-germ $f: (\bA^m,0) \to (\bA^{l},0)$ with the determinacy degree $k_f<\infty$, 
put 
\begin{equation*}
I_{f}:=f^*(\mathfrak{m}_{\bA^{l},0})\Ost_{\bA^m,0}+\mathfrak{m}_{\bA^m,0}^{k_f+1} 
\end{equation*}
(see also \S \ref{gs}). Let $\eta$ denote the $\K$-equivalence class of $f$, $\K$-type of $f$ for short, i.e., the isomorphism class of the quotient algebra $\Ost_{\bA^m,0}/I_f$  (\S \ref{map-germ}). 
Abusing the notation, we often write $I_\eta$ to mean the ideal $I_f$. 
We denote the colength by 
$$n(\eta):=\dim_k \Ost_{\bA^m,0}/I_{\eta} < \infty$$ 
and also set 
$$\ell(\eta):=\mbox{\em $\K_e$-codimension of $\eta$} < \infty$$ 
which is equal to the minimal dimension $l'$ of the target space of a stable mono-germ $(\bA^{m'},0) \to (\bA^{l'},0)$ with $\kappa=l'-m'$ whose quotient algebra is isomorphic to $\eta$ (\S \ref{codim}). 
We may also use the letter $\eta$ to mean an {\em equisingularty type}, a certain moduli family of $\K$-types 
(then the definition of $\ell(\eta)$ is modified by taking account of the modality). 
\item[(ii)] 
Throughout this paper, 
we let a {\em multi-singularity type} 
$$\underline{\eta}=(\eta_1, \eta_2, \cdots, \eta_r)$$ 
mean an {\em ordered} $r$-tuple of $\K$-types $\eta_i$  
of finitely determined mono-germs $(\bA^m,0) \to (\bA^{l},0)$ $(1\le i \le r)$. 
Types of entries $\eta_i$ are allowed to be duplicate, e.g., in case of $\kappa > 0$, we denote by $A_0^r=(A_0, \cdots, A_0)$ ($r$ times)  the $r$-tuple of immersion-germs, which will be used for expressing the $r$ multiple-point locus. 
In case of $\kappa\le 0$, we always assume that every entry $\eta_i$ is {\em not submersive}, namely, $\underline{\eta}$ is a combination of types of {\em isolated complete intersection singularities} (ICIS).  
In any case, we put 
$$n=n(\underline{\eta}):=\sum_{i=1}^r n(\eta_i), \qquad \ell=\ell(\underline{\eta}):=\sum_{i=1}^r \ell(\eta_i) $$
where $n$ is the total multiplicity and 
$\ell$ is equal to the expected codimension of the multi-singularity locus of type $\underline{\eta}$ in the target space. 
For a multi-singularity type $\underline{\eta}$, 
we define the {\em geometric subset} of the (ordered) Hilbert scheme of $n$ points in $X$
$$\Xi(X; \underline{\eta}) \subset X^{[[n]]}$$
to be the Zariski closure of the subset parametrizing ideal sheaves supported at $r$ distinct points of $X$ with stalks being isomorphic to $I_{\eta_i}$'s ($1\le i\le r$) (Definition \ref{geom_subset}). 
\item[(iii)]  
Given an {\em arbitrary} morphism $f: X \to Y$ from an $m$-fold $X$ to an $l$-fold $Y$, 
we introduce the {\em Hilbert extension map}
$$f^{[[n]]}: X^{[[n]]} \to (X\times Y)^{[[n]]}$$
(cf. \cite{Iarrobino, Gaffney, DL}). 
This can be thought of as a ``compactification" of the multi-jet extension map introduced by J. Mather \cite{MatherV}. 
Although the Hilbert schemes may be highly singular and have many components of different dimension, the map $f^{[[n]]}$ itself is always a regular embedding and the normal bundle is isomorphic to the tautological bundle $(f^*TY)^{[[n]]}$ (Proposition \ref{emb}). Then, Fulton-MacPherson's {\em refined Gysin map} \cite[\S 6]{Fulton}  is well adapted; it produces the intersection product of $f^{[[n]]}$ with the closed subscheme $\Xi(X; \underline{\eta})\times Y \hookrightarrow (X\times Y)^{[[n]]}$. 
Projecting it to $X$ and $Y$ provided $f$ is proper, 
we define {\em multi-singularity loci classes} 
as rational equivalence classes of algebraic cycles 
$$m_{\underline{\eta}}(f) \in \CH^{\ell-\kappa}(X), \quad n_{\underline{\eta}}(f) \in \CH^\ell(Y),$$
which are supported on the multi-singularity loci $M_{\underline{\eta}}(f) \subset X$ and $N_{\underline{\eta}}(f) \subset Y$, respectively (Definitions \ref{multi-singularity_loci} and \ref{multi-singularity_class}). It holds that 
$$f_*m_{\underline{\eta}}(f)=n_{\underline{\eta}}(f).$$ 
In case that $f$ is appropriately generic,  the loci $M_{\underline{\eta}}(f)$ and $N_{\underline{\eta}}(f)$ have the right dimension ($=\dim Y -\ell $) and represent exactly the classes $m_{\underline{\eta}}(f)$ and $n_{\underline{\eta}}(f)$, respectively, up to certain multiplicities (Proposition \ref{multi-sing}). For instance, if $f$ is generic and the target locus $N_{\underline{\eta}}(f)$ is of dimension $0$, then we have the formula counting $\underline{\eta}$-singular points in the target space
$$\deg N_{\underline{\eta}}(f)=\frac{1}{\#\Aut(\underline{\eta})} \int_Y n_{\underline{\eta}}$$ 
where $\Aut(\underline{\eta}) \subset \mathfrak{S}_r$ is the subgroup of permutations preserving the collection $\underline{\eta}$.
\item[(iv)] 
After a brief review of algebraic cobordism theory in \S \ref{alg_cob}, 
we show that the class $n_{\underline{\eta}}(f)$ depends only on the algebraic cobordism class 
$$[f: X \to Y] \in \Omega^{\kappa}(Y)$$ 
(Proposition \ref{prop1}). 
The proof is based on careful consideration of refined intersection caused by a double-point degeneration 
(Proposition \ref{degeneration_free}). 
Furthermore, we prove that the resulting map 
$$n_{\underline{\eta}}: \Omega^{\kappa}(Y) \to \CH^\ell(Y)$$ 
is a (unstable) cohomology operation in an analogous sense as Quillen's in topology (Proposition \ref{prop2} and Corollary \ref{prop3}). 
\item[(v)] 
In \S \ref{target_tp}, we complete the proof of our main theorem (Theorem \ref{main_thm1}). 
Thanks to a recent result of Vishik \cite{Vishik}, 
it turns out that the operation $n_{\underline{\eta}}$ tensored by $\Q$ is uniquely written as a polynomial with $\Q$-coefficients of the Landweber-Novikov operations with indices $I=(i_1, i_2, \cdots)$ 
$$s_I: \Omega^{\kappa}(Y) \to \CH^\ell(Y),  \;\; [f] \mapsto s_I(f)=f_*(c^I(f)). $$
That gives the desired multi-singularity Thom polynomial expression in the target space $Y$. 
A $K$-theoretic refinement of the theorem is also obtained in entirely the same way (\S \ref{multi_cob}). 
\item[(vi)]  
In  \S \ref{source_tp}, for an appropriately nice map (later called an {\em admissible map}) $f:X \to Y$ of codimension $\kappa \ge 0$, 
we identify the class $m_{\underline{\eta}}(f)$ in terms of $c_i(f)$ and $f^*s_I(f)$,  
that leads to the existence of the source multi-singularity Thom polynomial (Theorem \ref{main_thm2}). 
The proof uses the result of (v) above but also heavily relies on another aspect coming from local symmetry of singularity types, that is based on the so-called {\em restriction method} due to R. Rim\'anyi in topology \cite{Rimanyi01} (Remark \ref{restriction}). 
\end{itemize}

To be precise, what we prove in the step (v) is the existence of {\em cobordism Thom polynomials in Landweber-Novikov classes} defined in $\Omega^\ell(Y)$ (Theorem \ref{main_thm}). 
Since algebraic cobordism $\Omega^*$ is the universal oriented cohomology theory \cite{LevineMorel}, 
the Chow theory version in $\CH^*$ and the $K$-theory version in $K^0[\beta, \beta^{-1}]$ immediately follow. Notice that in another oriented cohomology theory $A^*$, there is no longer {\em a priori} `fundamental classes' for singular subvarieties in general. In particular, finding an $\underline{\eta}$-singularity loci class $n^\Omega_{\underline{\eta}}(f) \in \Omega^\ell(Y)$ is more or less equivalent to that one picks up a `universal' desingularization of the $\underline{\eta}$-locus. It is not unique, so $n^\Omega_{\underline{\eta}}(f)$ depends on the choice. For instance, there has been studied the Thom-Porteous formula in algebraic cobordism $\Omega^*$ by using specific resolutions of orbit closures of $\Sigma^i$ \cite{Hudson12}. 
This remark is important for properly employing Vishik's result on cohomology operations \cite[Thm.3.18]{Vishik}. 
In fact, we will show that given a certain distinguished class of the geometric subset $\Xi(X; \underline{\eta})$ in $\Omega_*(X^{[[n]]})$, called an {\em $\Omega$-assignment} later (Definition \ref{good}), there exists a unique universal polynomial in $s_I$'s expressing the corresponding $\underline{\eta}$-singularity loci class $n^{\Omega}_{\underline{\eta}}(f)$ in $\Omega^\ell(Y)_\Q$. It should be emphasized that in our approach, the existence of multi-singularity Thom polynomials is implicitly indebted on resolution of singularities of $\Xi(X; \underline{\eta})$. 

Notice that the existence of the target Thom polynomial for $n_{\underline{\eta}}(f)$ does not directly imply that of the source Thom polynomial for $m_{\underline{\eta}}(f)$. 
In the step (v), our main formula of the target class $n_{\underline{\eta}}(f)$ is valid for {\em arbitrary} proper maps $f$, 
while in the step (vi) for the source class $m_{\underline{\eta}}(f)$, 
we consider only appropriately nice maps with $\kappa \ge 0$ because of some technical reason. 
It is highly expected that this restriction can be dropped, but it is surely related to some essential issues, and therefore this part of the conjecture still remains open (see \S \ref{bi}, \S \ref{critical} and \S \ref{classifying_stack}). 

We consider separately the two cases of positive and negative codimension $\kappa$ of maps (cf. \S \ref{multiple} and \S \ref{icis}). 
From the above construction, the source and target multi-singularity loci classes are essentially determined by the geometric subset $\Xi(X; \underline{\eta})$ in $X^{[[n]]}$. Namely, the roll of $Y$ is somehow auxiliary, and it infers that the form of Thom polynomial should depend only on the isomorphism class of the quotient algebra associated to $\underline{\eta}$ {\em with $\kappa$ as a parameter} (see \S \ref{series}). That is very related to tautological integrals over Hilbert schemes  \cite{BercziSzenes, Kaz17, Rennemo, Berczi, BercziSzenes21, Berczi23}. In fact, along this line, B\'erczi-Szenes \cite{BercziSzenes21} recently announced that they have proven {\em target multiple-point formulas} of $n_{A_0^r}(f)$ for generic maps  (see \S \ref{multiple}), i.e.,  in our side, it corresponds to Theorem \ref{main_thm1} for $\underline{\eta}=A_0^r$ for generic maps with $\kappa>0$. 


Our proof using cohomology operations is not constructive, and therefore it does not directly yield any explicit forms of universal polynomials. 
Nevertheless, once we assume the existence of universal polynomials, there is an effective method for computing their explicit forms based on equivariant cohomology using local symmetries of singularity types, that is called the {\it interpolation} or {\it restriction method} invented by R. Rim\'anyi \cite{Rimanyi01, FR02, FR12}. 
That reduces the problem to solving a system of linear equations of unknown coefficients in  the universal polynomial of a given type. 
Indeed, by this method, explicit computations have already been given for several mono and multi-singularity types  \cite{Rimanyi01, Rimanyi02, FR02, Kaz03, Kaz06, Ohmoto16} -- for instance, Marangell-Rim\'anyi \cite{MarangellRimanyi} calculated the {\em general quadruple-point formula of $m_{A_0^4}(f)$} which should be valid even for non-curvilinear maps (\S \ref{multiple}). 
Notice that our existence theorems ensure all of those computational results for multi-singularities of maps.  
We would not repeat such computations here. 
Rather, at the last step (vi), we essentially employ the restriction method to identify the source multi-singularity Thom polynomials. 

We end this introduction by briefly mentioning a more general framework. 
Besides singularity theory of maps, the whole perspective envisioned by Kazarian in \cite{Kaz03} involves the Thom-Pontryagin construction, Thom spectra, classifying spaces,  infinite loop spaces,  cohomology operations, and even more, Gromov's h-principle, all of which are deep concepts in differential topology and homotopy theory. 
However, the complex cobordism theory $MU^*$, which is described in $C^\infty$ category, does not fit the classification theory of singularities in complex analytic category. 
Therefore, as suggested in \cite{Kaz03} for finding an entire proof of the conjecture, we should unite and integrate those topological ideas suitably in algebraic geometry side (possibly, including positive characteristic case as well); especially, most desirable is to find an algebro-geometric {\em classifying stack for multi-singularity types of maps}. 
That is very challenging -- in such an attempt, one would encounter Voevodsky-Morel's motivic cohomology and Artin stacks for classifying multi-singularities of maps (see \S \ref{classifying_stack} and \S \ref{motivic}). 
As mentioned before, however, our strategy in the present paper avoids this deeper part of the problem and instead takes a shortcut by directly using Vishik's theory of cohomology operations for $\Omega^*$ over a field of characteristic $0$. 
So, in this sense, the full scope of the Thom-Kazarian principle is still left for future work, see \S \ref{perspective} for more comments. 

Throughout this paper, 
$k$ will denote the fixed base field of characteristic $0$, 
and $\Sch$  and $\Sm$ will denote respectively 
the category of separated schemes of finite type over $k$ and 
the full subcategory of smooth quasi-projective varieties over $k$. 
Also we denote by $\CH_*(X)$ the Chow group for $X \in \Sch$ and by 
$$\CH^*(X):=\CH_{m - *}(X)$$ 
the Chow ring for $X \in \Sm$ with $m=\dim X$. 
For $X \in \Sch$, $\CH^*(X)$ stands for the operational Chow ring \cite{Fulton}. 
In case of rational coefficeints, we write $\CH^*(X)_\Q$ etc. 

The rest of the present paper is organized as follows. 
In the first part, \S 2, \S 3 and \S 4, we prepare necessary materials, and especially we introduce multi-singularity loci classes in $\CH^*$ and state main results precisely. The second part, \S 5 and \S 6, is devoted to the proofs of theorems using algebraic cobordism $\Omega^*$ and cohomology operations. In the final section \S 7 we discuss future perspectives.

\section{Preliminaries}

\subsection{Motivational background} \label{motivation}
Before entering the main body, we briefly explain our motivation, historical remarks and the current state of the Kazarian conjecture. 

\subsubsection{\bf Degeneracy loci formula}\label{TP} 
Let $E$ and $F$ be vector bundles over $B$ in $\Sch$ of rank $m$ and $l$, respectively. 
For a bundle map $h: E \to F$, 
the {\em $k$-th degeneracy loci class} $\mathfrak{D}_k(h) \in \CH_*(B)$ is defined so that it is represented by the locus 
$$D_k(h)=\{\, p \in X \, | \, \dim \ker [h_p: E_p \to F_p] \ge k \, \}$$ 
if the locus has the right dimension \cite[\S 14]{Fulton}. The {\em Thom-Porteous formula} says that 
the class is expressed by a special Schur polynomial: 
\begin{equation*}
\mathfrak{D}_k(h) =  
\Delta_{(\kappa+k)^k}(c(F-E)) \frown [B]
\end{equation*}
where $\kappa:= l-m$ and 
$$\Delta_{(\kappa+k)^k}(c)=\det[c_{\kappa+k+j-i}]_{1\le i, j \le k}.$$ 
See \cite{Porteous, Ronga72, Damon, Fulton, FultonPragacz} for the detail. 
In case of maps $f: X \to Y$ in $\Sm$, 
applying the formula to the differential $df: TX \to f^*TY$, 
the Schur polynomial is just 
the Thom polynomial for the 1st order Thom-Boardman singularity type of $f$.

It is very natural to deal with this formula in the framework of {\em Equivariant Intersection Theory} \cite{Totaro, EdidinGraham}. 
For a reductive linear algebraic group $G$, the universal principal bundle $EG \to BG$ is associated by the {\em algebraic Borel construction} \cite{Totaro}:  any $G$-principal bundle over $B$ can be induced from the universal bundle via the so-called {\em classifying map}\footnote{Precisely speaking, the algebraic classifying map $\rho$ is defined on the total space of some affine bundle $B' \to B$ by a tautological construction \cite[Lem.1.6, Thm.1.5]{Totaro}. Since it does not affect any computations of characteristic classes, we simplify our notations so as not to touch on this technical issue.}  to $BG$. 
Now, let us consider 
$$G=GL(m) \times GL(l), \quad V=\bA^{lm},$$ 
the space of $l\times m$ matrices, 
where $G$ acts on $V$ by $(P, Q).H:=QHP^{-1}$.  Then $V$ is decomposed into finitely many $G$-orbits;  let $\Sigma^k \subset V$ denote the orbit of linear maps with kernel dimension $k$. Consider associated vector bundles with fiber $V$ and group $G$ and the classifying map, namely, the vector bundle $\Hom(E, F)$ is induced from the universal vector bundle $EG \times_G V$: 

$$
\xymatrix{
\Hom(E, F) \ar[r]^{\bar{\rho}} \ar[d]^{\pi}& EG \times_G V \ar[d]^{\pi}  \\
B   \ar@/^10pt/[u]^{h}  \ar[r]_{\rho} & BG 
}
$$
There are also associated cone bundles with fiber $\overline{\Sigma^k}$ (the bar means the closure): 
$$\overline{\Sigma^k}(E, F) \subset \Hom(E, F), \quad 
EG \times_G \overline{\Sigma^k} \subset  EG \times_G V$$ 
where the former is induced from the latter universal one. 
A bundle map $h: E\to F$ corresponds to a section of $\Hom(E, F)$, and indeed, the degeneracy loci class $\mathfrak{D}_k(h)$ is defined to be the intersection product of $h$ and $\overline{\Sigma^k}(E, F)$ in the sense of Fulton-MacPherson 
(see \cite[Rem.14.3]{Fulton} and also \S \ref{Intersection} later). 
By the above diagram, $\mathfrak{D}_k(h)$ is obtained from the {\em equivariant fundamental class of the orbit closure}\footnote{The Chow ring of the stack $[V/G]$ is defined to be the equivariant Chow ring $\CH^*_G(V)$ \cite{EdidinGraham}. } 
$$[EG \times_G \overline{\Sigma^k}] \in \CH^*_G(V) \stackrel{\simeq}{\longleftarrow} \CH^*_G(pt) = \CH^*(BG).$$
Since $\CH^*(BG)$ is a polynomial ring generated by universal Chern classes of the source and the target bundles, the above fundamental class is uniquely written by a universal polynomial,  that is the {\em Thom polynomial of type $\Sigma^k$}, and a slightly technical computation shows up that it is $\Delta_{(\kappa+k)^k}(c)$.  
The reason why the Thom polynomial is actually written in the quotient Chern classes $c_i=c_i(F-E)$ lies in a certain stabilization of the group actions with fixing $\kappa=l-m$ -- essential is that the embedding $\bA^{lm} \subset \bA^{(l+s)(m+s)}$ of the spaces of linear maps via the trivial suspension is always transverse to any orbit $\Sigma^k$.

The above picture is completely applicable to another suitable classification theories, e.g.,   for  quiver representations and for mono-singularities of maps.  
The former case has been rapidly developping in modern Schubert calculus related to mathematical physics. 
The latter is, more precisely, for the classification of map-germs $(\bA^m,0) \to (\bA^l, 0)$ under the group $\mathcal{G} (=\A, \K)$ of local non-linear coordinate changes, that is the core of the Thom-Mather theory which we will review in \S \ref{map-germ}. In this case, we take again the linear group (the $1$-jet part), $G=\mathcal{G}^1=GL(m)\times GL(l)$, which acts on the jet space $V=J^k(m, l)$ (=the space of Taylor coefficients of order $\le k$), and obviously $G$ preserves any $\mathcal{G}^k$-orbit $\eta_V \subset V$.  Then, in entirely the same way as the Thom-Porteous case, we define the {\em mono-singularity loci class of type $\eta_V$ for maps $f$} and show the existence of the {\em Thom polynomial associated to $\eta_V$} \cite{Thom, HaefligerKosinski, Porteous, Damon, Ronga73, Rimanyi01, FR02, Kaz03, Kaz06, Ohmoto06, BercziSzenes, FR12, Ohmoto16, Kaz17}.

\

\subsubsection{\bf Multiple-point formulas}\label{multiple} 
Another important prototype is the so-called {\it multiple-point formulas}, 
which have been explored in depth around 80's  
\cite{Colley, DL, Kleiman81, Kleiman90, LeBarz, Piene, Ran85, Ronga73, Ronga84, Vainsencher}; readers are referred to excellent surveys \cite{Kleiman77, Kleiman87}. 

Let $f: X \to Y$ be a proper morphism in $\Sm$ 
with positive codimension $\kappa=l-m>0$ ($l=\dim Y$, $m=\dim X$).  

The double-point locus of $f$ is the Zariski closure of the set of points $p_1 \in X$ such that there exist another points $p_2\; (\not=p_1) \in X$ with $f(p_1)=f(p_2)$.  
The {\em double-point loci class} $m_{A_0^2}(f) \in \CH^*(X)$ represents the locus if it has the right dimension, and there is a universal expression of the class in terms of invariants of $X$, $Y$ and $f$ -- that is the {\em double-point formula} \cite{Ronga73, FultonLaksov} (see  \cite[\S9.3]{Fulton}): 
$$m_{A_0^2}(f)=f^*f_*(1)-c_{\kappa}(f)$$
where $c_{\kappa}=c_{\kappa}(f)=c_\kappa(f^*TY-TX)$. 
The next is the {\em triple-point formula} which expresses 
the closure of the set of $p_1 \in X$ such that 
there are distinct $p_2, p_3 \,(\not=p_1) \in X$ mapped to $f(p_1)$. 
That is given by the following formula with a nested structure \cite{Kleiman81, Kleiman90, Ronga84}: 
$$m_{A_0^3}(f)=f^*f_*(m_{A_0^2}(f))-2c_\kappa f^*f_*(1)+2\left(c_\kappa^2+\sum_{i=0}^{\kappa-1}2^ic_{\kappa-i-1}c_{\kappa+i+1}\right).$$ 
Those formulas are obtained by explicitly constructing desingularizations of the double and the triple-point loci by blow-ups.  One may expect that the same method would work inductively for finding higher multiple-point formulas. 
However, it turns out to be very hard in general -- for instance, the general {\em quadruple-point formula} of $m_{A_0^4}(f)$ is already problematic \cite{DL}. 
The reason will be mentioned briefly in the end of this subsection. 

On the other hand, by the projection formula 
$$f_*(\alpha \cdot f^*\beta)=f_*(\alpha)\cdot \beta,$$
the right hand sides of the double and triple-point formulas are written in terms of Chern classes $c_i(f)$ and Landweber-Novikov classes 
$$f^*s_0=f^*f_*(1), \;\; f^*s_{0\cdots01}=f^*f_*(c_\kappa), \;\; f^*s_0^2=f^*f_*f^*f_*(1).$$
A far-extension of this observation finally arrives at Kazarian's conjecture \cite{Kaz03} mentioned in Introduction. 

The conjecture is strongly supported by {\em curvilinear multiple-point theory} -- it was successfully established by Kleiman, Colley and Ran \cite{Colley, Kleiman81, Kleiman90, Ran85} 
for an important particular class of maps, called {\it curvilinear maps}, 
which are morphisms $f: X \to Y$ in $\Sm$ with $\kappa \ge 0$ and corank at most one (i.e., $\dim \ker df_p \le 1$ for any $p \in X$).  
A curvilinear singularity type is denoted by $A_\mu$, i.e., its local ring is isomorphic to $k[x]/\langle x^{\mu+1}\rangle$, and a curvilinear multi-singularity type is denoted by $A_\bullet=(A_{\mu_1}, \cdots, A_{\mu_r})$. 
For generic curvilinear maps, the classes $m_{A_\bullet}(f)$ with some specific types $A_\bullet$ are explicitly computed by an algorithm in \cite{Colley} and closed recursive formulas on $m_{A_0^r}(f)$ have also been studied  \cite{Kleiman81,  Ran85, Kleiman90, Kaz08}. 
The multiple-point theory covers many important applications in classical enumerative geometry, while those formulas may return wrong answers when applying to non-curvilinear maps, i.e., one needs some correction terms supported on singularities of corank greater than one. 
Note that non-curvilinear maps map appear even in low dimensional cases, so an uncovered area has still remained largely, see, e.g., \cite{DL, MarangellRimanyi, Nekarda}.

Another remarkable evidence of the conjecture is recently brought about from a new approach of complex analytic geometry with {\em non-reductive geometric invariant theory} by B\'erczi-Szenes \cite{BercziSzenes, BercziSzenes21} and B\'erczi \cite{Berczi, Berczi23}. 
They have been developing advanced residue techniques for computing  $m_{A_\mu}(f)$ (mono-singularity of type $A_\mu$) and the target multiple-point formula of $n_{A_0^r}(f)$ for generic maps $f$ with $\kappa >0$ via tautological integrals on Hilbert schemes (cf. \S \ref{instanton}, \S \ref{tautological}). 
The latter result also strongly supports part of the conjecture. A particular advantage of their approach is that it is constructive via direct computations, that is in clear contrast to our approach. 

For proving the double-point formula, a key observation is that the scheme-theoretic preimage $(f\times f)^{-1}(\Delta_Y)$ of the diagonal $\Delta_Y \subset Y \times Y$ contains not only the double-point locus of $f$ but also the diagonal $\Delta_X \subset X \times X$. 
To remove the unnecessary contribution from $\Delta_X$,  we may use the blow-up along $\Delta_X$ 
$$\pi: {\rm Bl}(X\times X) \to X \times X.$$ 
The embedding $X \to X \times Y$ of the graph of $f$ naturally induces an embedding 
$$f^{\rm Bl}: {\rm Bl}(X\times X) \to {\rm Bl}((X \times Y)\times (X \times Y))$$
where the target is the blow-up along the diagonal $\Delta_{X\times Y}$. 
Then we have the intersection product of the embedding $f^{\rm Bl}$ with the subscheme 
$$\mbox{Bl}(X\times X) \times \Delta_Y \hookrightarrow 
\mbox{Bl}((X \times Y)\times (X \times Y))$$
in the sense of Fulton-MacPherson, and pushing the product to the first factor $X$ of $\pi$, we define the double-point loci class $m_{A_0^2}(f) \in \CH^\kappa(X)$. 
The contribution from $\Delta_X$ is just the excess intersection along the exceptional divisor, and computing it explicitly,  the desired formula is obtained \cite[\S 9.2]{Fulton}. 
Emphasized is that the formula is valid for {\em arbitrary} morphisms $f$ in $\Sm$ ($\kappa >0$). 
Taking a further blow-up, the triple-point formula is derived in a similar but more involved way. 

A key issue for us is how to generalize this picture of the proof. 
To do this, we employ {\em ordered} Hilbert schemes $X^{[[n]]}$ with the projection 
$$\pi: X^{[[n]]} \to X^n=X\times \cdots \times X.$$
Associated to the graph map of $f$, we introduce the {\it $n$-th Hilbert extension map} 
$$f^{[[n]]}: X^{[[n]]} \to (X\times Y)^{[[n]]}$$ 
(cf. \cite{Gaffney, DL}). 
Note that $f^{[[2]]}=f^{\rm Bl}$ with 
$X^{[[2]]}=\mbox{Bl}(X\times X)$ and $(X\times Y)^{[[2]]}=\mbox{Bl}((X\times Y)\times (X\times Y))$, thus $f^{[[n]]}$ must be a suitable generalization. 
In fact, $f^{[[n]]}$ is a regular embedding (Corollary \ref{Hextension}), and thus we have 
the intersection product of $f^{[[n]]}$ with the subscheme 
$$X^{[[n]]}\times \Delta_Y \hookrightarrow (X\times Y)^{[[n]]},$$
which leads to higher order multiple-point loci classes $m_{A_0^n}(f)$ (Definition \ref{multi-singularity_class}). 

In case of $n\le 3$, the Hilbert scheme $X^{[[n]]}$ is smooth for any dimension of $X$. 
In particular, the preimage of $X^{[[3]]}\times \Delta_Y$ via $f^{[[3]]}$ provides a desingularization of the triple-point locus of generic maps $f$, which coincides with the desingularization by blow-ups used for proving the triple-point formula \cite{Ronga84, Gaffney} (see Remark \ref{tri}). 
In case that $f$ is curvilinear, for any $n$, 
it suffices to consider open subsets of $X^{[[n]]}$ consisting of $0$-dimensional subschemes lying on smooth curve-germs, called curvilinear Hilbert schemes; they are actually smooth and used for desingularizing curvilinear multiple-point loci \cite{Kleiman90}. 
In such nice cases, in principle, everything can be handled by traditional methods within $\Sm$. 

It turn out, however, that the general case including non-curvilinear maps $f$ with $n \ge 4$ and $\dim X \ge 3$ is very problematic.  
It is because the Hilbert scheme $X^{[[n]]}$ itself becomes to be highly singular and  we can not ignore the effect caused by its complicated singularities and various dimensional components. 
For instance, singularities of $X^{[[4]]}$ cause fatal difficulties in constructing a desingularization of the quadruple-point locus for {\em non-curvilinear maps} $f$, as demonstrated by Dias-Le\,Barz \cite{DL}.  

Instead of Hilbert schemes of points, we could have begun with the Fulton-MacPherson configuration space, which is also a compactification of the configuration space of distinct $n$ points and is always smooth \cite{FM}.  Nevertheless, there is a trade-off; the FM configuration space has also its own complicated structure depending on the choice of centers for consecutive blow-ups in its construction, and especially it does not have a functorial property which the Hilbert scheme possesses. 
Anyway, we conclude that it is very hard to reach the final goal, solving the conjecture, by explicit calculations through piling up step-by-step construction of desingularizations of multi-singularity loci. 
In the present paper we explore a new {\em non-constructive} approach in detail.

\subsubsection{\bf Discriminant of ICIS and Legendre maps} \label{icis} 
The case of maps with non-positive codimension ($\kappa \le 0$) belongs to {\em non-curvilinear} enumerative geometry --  known results are quite limited (see, e.g., \cite{Vainsencher, Vainsencher07}). 
In his second paper \cite{Kaz06}, Kazarian extended his conjectures to maps of $\kappa \le 0$. 
Our approach deals with this case as well; indeed, the above Hilbert extension maps are defined without any change. 

In this case, we work with {\em isolated complete intersection singularities} (ICIS), 
e.g.,  $A_\mu$-singularity type defined by $x_1^{\mu+1}+x_2^2+\cdots+x_{|\kappa|+1}^2=0$. 
Namely, given a proper morphism $f: X \to Y$ with $\kappa \le 0$, 
we are concerned with counting singular fibers which have a prescribed combination of several ICIS's. 
In other words, we consider the critical scheme $C(f)=M_{A_1}(f)$ in $X$,  where ${\rm coker}\, df$ has positive dimension, and want to count multiple-points of prescribed type for the singular map 
$$f_C: C(f) \to Y$$ 
obtained by restricting $f$ to $C(f)$ (note that $C(f)$ may not be smooth in general). 
The image of $f_C$ is usually called the {\em discriminant set} of $f$, denoted by $D(f) \subset Y$. 
Of particular interest is the case where $\dim {\rm coker}\, df$ is at most one and $df: X \to \Hom(TX, f^*TY)$ is transverse to the 1st Thom-Boardman stratum $\Sigma^1$. Then 
$C(f)\, (=df^{-1}(\Sigma^1))$ is smooth and $D(f)$ is locally the discriminant of {\em isolated hypersurface singularities}; then, the map  $f_C$ is often referred to as a {\em Legendre map} \cite{Arnold}. 
We will come back to this case in \S \ref{ex_tri2} and \S \ref{leg}. 

\subsubsection{\bf Hilbert's 15th problem and beyond} \label{15th}
In classical literature, a traditional way to study a non-singular projective variety $X$ is to first map it into some projective space $\bP^l$ and count singular points of the image variety -- elementary numerical characters of $X$, named {\it rank, class, class of immersions}, and so on, have been introduced in this way, see, e.g.,  \cite{Salmon, SempleRoth, Kleiman77, Piene}. For instance, Semple-Roth \cite{SempleRoth} noted that for a projective surface there are exactly four basic characters which generate all other characters of Salmon-Cayley-Zeuthen-Baker, often referred to as Enriques' formulas, and thereafter, Roth claimed an analogous statement for a projective $3$-fold that there are seven basic characters which recover dozens of others. That is directly related to the nature of cobordism theory from a modern viewpoint. 
In our language, the target Thom polynomials $n_{\underline{\eta}}$ express those numerical characters in terms of Landweber-Novikov classes $s_I$ of maps $f: X \to \bP^l$, and it turns out that the set of basic characters  exactly correspond to the set of $s_I$'s --  in case of surfaces, there are four Chern monomials $1, c_1, c_1^2, c_2$,  and for $3$-folds, there are seven $1, c_1, c_1^2, c_2, c_1^3, c_1c_2, c_3$, and their pushforward via $f$ generate cobordism invariants  \cite{Sasajima17}. 

As well-known, eminent Hilbert's 15th problem was asking for a rigorous justification of the enumerative geometry of H. Schubert and other pioneers in their era. This is surely more than the quest about the ring structure of cohomology of Grassmannians and flag varieties, because most chapters of Schubert's book \cite{Schubert} are devoted to developing and deepening his methods for solving various kinds of enumerative problems in projective algebraic geometry. After almost a hundred years, Fulton-MacPherson's Intersection Theory \cite{Fulton} successfully provides a firm foundation of enumerative geometry, that is widely recognized nowadays. However, indeed less known, even in a classical context, a number of questions still remains unexplored well, which is consolidated to Kazarian's conjecture.  
Solving the conjecture revives this 19th century enumerative geometry in a modern way with new flavors \cite{Kaz03, Kaz06, MarangellRimanyi, Ohmoto16, Sasajima17, Nekarda, NekardaOhmoto1}. For instance, as an application of Thom polynomials, Nekarda \cite{Nekarda} has expanded a classical work of H. Schubert on contact problems of lines with a projective surface in $\Proj^3$ and $\Proj^4$ with some interests from computer vision. Also one can apply Thom polynomial theory to some questions from mathematical physics such as instanton counting problem, see \S \ref{instanton}.

\subsection{A digest of the Thom-Mather Theory}\label{TM}
Singularity theory of differentiable mappings has been established in both complex analytic category and real $C^\infty$ category (cf. \cite{MatherIV, MatherV, Arnold, MondNuno}). 
Here we briefly present a minimal digest of the local theory for later use. 
In particular, we describe it in algebraic category over an algebraic closed field $k$ of characteristic zero, according to Mather \cite{Mather}. 
In this context, classification of singularities of maps is essentially considered at the level of {\em formal schemes}. 

\subsubsection{\bf Map-germs and equivalence relations} \label{map-germ}
Let $X$ be a scheme and $p\in X$ a closed point. 
Let $\Ost_{X,p}$ denote the stalk at $p$ of the structure sheaf $\Ost_{X}$ and $\mathfrak{m}_{X,p}$ its maximal ideal. 
Also let $\hat{\Ost}_{X,p}$ denote the $\mathfrak{m}_{X,p}$-adic completion of $\Ost_{X,p}$, 
and  $\hat{\mathfrak{m}}_{X,p}$ the maximal ideal generated by  $\mathfrak{m}_{X,p}$. 
For a finite set $S$ of closed points in $X$,  ${\Ost}_{X,S}=\bigoplus_{p \in S} {\Ost}_{X, p}$ ($\hat{\Ost}_{X,S}$ as well). 

Let $f: X \to Y$ be a morphism in $\Sm$ and $q \in Y$ a closed point and $S \subset f^{-1}(q)$ finite closed points. 
The {\em map-germ $f: (X, S) \to (Y, q)$ of $f$ at $S$} means the ring homomorphism $f^*: {\Ost}_{Y, q} \to {\Ost}_{X, S}$ 
(it is called a {\em multi-germ} if $|S|>1$, {\em mono-germ} if $|S|=1$). 

Two mono-germs $f, g: (X, p) \to (Y, q)$ are {\em formally $\K$-equivalent} if their formal quotient algebras are isomorphic:  
$$\hat{\Ost}_{X,p}/f^*(\hat{\mathfrak{m}}_{Y,q})\hat{\Ost}_{X,p} \simeq \hat{\Ost}_{X,p}/g^*(\hat{\mathfrak{m}}_{Y,q})\hat{\Ost}_{X,p}.$$ 
The equivalence class of $f$ is called the (formal) {\em $\K$-type of $f$} throughout the present paper. 
The letter $\K$ refers to the {\em contact equivalence} of map-germs  -- we often use the {\em stably $\K$-equivalence} which is the equivalence among formal quotient algebras associated to map-germs with a fixed codimension $\kappa$ \cite{MatherIV, Arnold, MondNuno}.

We say that $f: (X, p) \to (Y, q)$ is {\em $k$-$\K$-determined} if any mono-germ $g$ satisfying that 
$$g^*({\mathfrak{m}}_{Y, q}){\Ost}_{X, p}\equiv f^*({\mathfrak{m}}_{Y, q}){\Ost}_{X, p} \;\; \mbox{modulo} \;\;  {\mathfrak{m}}_{X, p}^{k+1}$$ 
is formally $\K$-equivalent to $f$. 
Namely, the $\K$-type of $f$ is determined by the Taylor expansion of $f$ at $p$ up to order $k$ in some affine local coordinates of $X$ and $Y$ centered at $p, q$, respectively. Also $f$ is {\em finitely $\K$-determined} if it is $k$-$\K$-determined for some $k$. 
The minimal of such $k$ is called the {\em determinacy degree} of $f$ at $p$; denote it by $k_f(p)$. 
In case of non-positive codimension $\kappa\le 0$ (i.e., $\dim X \ge \dim Y$), 
$f$ is a finitely determined map-germ if and only if $f=0$ defines an isolated complete intersection singularity (ICIS) at $p$. 

There is another natural equivalence among map-germs: we say mono-germs $f, g: (X, p) \to (Y, q)$ are {\em formally $\A$-equivalent} if 
formal ring homomorphisms $f^*, g^*: \hat{\Ost}_{Y, q} \to \hat{\Ost}_{X, p}$ commute with each other via formal local isomorphisms of the source and the target spaces. It is known that formal $\A$-equivalence implies formal $\K$-equivalence, but the converse does not hold  (cf. \S \ref{locally_stable}).

\subsubsection{\bf Sheaf of vector-fields along $f$} 
Let $f: X \to Y$ be a morphism in $\Sm$.
Let $\theta(f)$ denote the sheaf of regular rational vector-fields along $f$, i.e., sections of $f^*TY$. 
Precisely, for open affine sets $U \subset X$ and $V \subset Y$ with $f(U) \subset V$, 
it assigns a $\Ost_X(U)$-module 
$$\theta(f)(U) =\{\, \mbox{derivations $D: \Ost_{Y}(V) \to \Ost_{X}(U)$}\, \}$$
i.e.,  $D$ is $k$-linear and satisfies 
$$D(\varphi \psi)=(f^*\varphi)D\psi+(f^*\psi)D\varphi \quad (\varphi, \psi \in \Ost_{Y}(V)).$$ 
In particular, $\theta(id_X)$ is the sheaf of regular rational vector-fields (derivation) on $X$, denoted also by $\theta_X (=\mathcal{D}er_X)$. 

We associate two sheaf homomorphisms 
$$tf: \theta_X \to \theta(f), \quad \omega f: \theta_Y \to \theta(f)$$
in a natural way as follows. For any open affine sets $U, V$ as above, for $v \in \theta_X(U)$, we define 
$$tf(v): \Ost_Y(V) \to \Ost_X(U), \;\; \varphi \mapsto v(\varphi \circ f).$$
Also for $w \in \theta_Y(V)$, we define 
$$\omega f(w): \Ost_Y(V) \to \Ost_X(U), \;\; \varphi \mapsto w(\varphi) \circ f.$$

We denote the stalk of $\theta(f)$ at $p$ by ${\theta}(f)_p$, 
and the $\mathfrak{m}_{X, p}$-adic completion by $\hat{\theta}(f)_p$. 
Put also $\hat{\theta}(f)_S=\bigoplus_{p \in S} \hat{\theta}(f)_p$ and so on. 
In particular,  we associate 
$$tf: \hat{\theta}_{X, S} \to \hat{\theta}(f)_S, \quad \omega f: \hat{\theta}_{Y, q} \to \hat{\theta}(f)_S.$$

\subsubsection{\bf Locally stable maps} \label{locally_stable} 
A morphism $f: X \to Y$ in $\Sm$ is {\em locally stable at $S$} in the sense of Thom-Mather 
(or say, the germ of $f$ at $S$ is locally stable) if it holds that 
$$\hat{\theta}(f)_S = tf(\hat{\theta}_{X, S})+\omega f(\hat{\theta}_{Y, q})$$
at the formal level \cite{Mather, MatherIV}. 
Namely, any formal infinitesimal deformation of $f$ is naturally recovered by formal vector-fields of the source and the target spaces. 
The local stability of $f$ at $S$ implies 
$${\theta}(f)_S = tf({\theta}_{X, S})+\omega f({\theta}_{Y, q})+\mathfrak{m}_{X,S}^c\theta(f)_S$$ 
for any power $c$. 

In \cite[Thm.A]{MatherIV}, 
Mather also showed that if two locally stable germs are $\K$-equivalent (or more stronger, $\K$-equivalent up to $(\dim Y+1)$-jets), then they are $\A$-equivalent. 
Thus, to classify locally stable map-germs up to $\A$-equivalence, it suffices to classify quotient formal algebras. 
There is a translation of this theorem in algebraic geometry setting, see \cite[\S 5]{Mather}. 

We call $f: X \to Y$ a {\em locally stable map} if for every closed point $q \in Y$ and every finite set of closed points $S \subset f^{-1}(q)$, the germ $f: (X, S) \to (Y, q)$ is locally stable. 

\begin{rem}\upshape
{\bf (Generic linear projection)} 
Locally stable maps form a sufficiently large class of maps in projective algebraic geometry. 
Actually,  the locally stability of maps is an {\em open condition}, provided a certain projectivity condition is assumed \cite[\S 2]{Mather}. 
A typical case is generic linear projection of a smooth scheme $X \subset \Proj^r$ ($m=\dim X$). 
Let $W$ be an open subset of the moduli of linear projections $\pi$ from $\Proj^r$ to $Y=\Proj^l$ ($r>l$) 
such that the center $C_\pi$ intersects $X$ transversely (if $m \le l$, then $C_\pi\cap X=\emptyset$). 
In \cite[\S 3]{Mather}, it is shown that if  $(m, l)$ is a certain pair of dimensions (i.e., it is in the so-called {\em nice range of dimensions}), then 
there is an non-empty Zariski open subset of $W$ consisting of $\pi$ for which $\pi: X-C_\pi \to Y$ is  locally stable. 
\end{rem}

\subsubsection{\bf $\K_e$-codimension} \label{codim} 
By using the Malgrange preparation theorem, the local stability of map-germs is characterized in a more convenient form as follows \cite{MatherIV, Arnold, MondNuno}. 

The {\em $\K_e$-codimension} of $f:(X, p) \to (Y, q)$ is defined by 
$$\ell(f_p):= \dim_k {\theta}(f)_p/(tf({\theta}_{X,p})+f^*({\mathfrak{m}}_{Y,q}){\theta}(f)_p)$$
(often written by $d_e(f, \K)$ in literature). 
This is invariant under (stably) $\K$-equivalence.  
It is known that $f$ is finitely determined if and only if $\ell(f_p)$ is finite \cite{MatherIV, Arnold}. 
In particular, it is shown \cite[Prop.1.6]{MatherIV} that the local stability of a multi-germ $f:(X, S) \to (Y, q)$ is equivalent to the sujectivity of the linear map induced by $\omega f$ (the {\em Kodaira-Spencer map}): 
$$\bar{\omega} f: T_qY \to 
\bigoplus_{p \in S} {\theta}(f)_p/(tf({\theta}_{X,p})+f^*({\mathfrak{m}}_{Y,q}){\theta}(f)_p).$$
For each $p \in S$, let $A_p$ be the kernel of the linear map $\bar{\omega} f_p$ from $T_qY$ to the $p$-component of the direct sum. 
The surjectivity of the above $\bar{\omega} f$ means that $f$ is locally stable at each $p\in S$ (i.e., $\bar{\omega} f_p$ is surjective) and all $A_p$ ($p \in S)$ have {\em regular intersections} in $T_qY$, i.e., 
$${\rm codim} (\bigcap_{p \in S} A_p) = \sum_{p \in S} {\rm codim} A_p = \sum_{p\in S} \ell(f_p)=: \ell(f_S).$$ 
In particular, if  $f:(X, S) \to (Y, q)$ is locally stable, then 
$$\ell(f_S) \le l\; (=\dim Y).$$ 
Namely, $\ell(f_S)$ is the minimal dimension of the target space of a locally stable multi-germ which is stably $\K$-equivalent to $f$. 
Hence, $\ell(f_S)$ is exactly equal to the codimension of the corresponding multi-singularity stratum in $Y$ of the locally stable map $f$ \cite[p.227]{MatherIV}. Later, we will denote by $N_{\underline{\eta}}^\circ(f)$ the stratum in $Y$, where $\underline{\eta}$ is the multi-singularity type of the germ of $f$ at $S$ as defined in Introduction (ii). 

\begin{exam} \upshape 
A multi-germ $f_S$ of $r$-tuple of immersions ($\kappa>0$) with regular intersections 
 is locally stable, and then, $\ell(f_S) (=\ell(A_0^r)) = r\kappa$, that is the codimension of the $r$-th multiple-point locus $N_{A_0^r}^\circ$ in the target space $Y$. 
\end{exam}

\subsection{A digest of Intersection Theory}\label{Intersection}
For later use, we briefly recall a few basics in Intersection Theory of Fulton-MacPherson  
\cite[\S 6]{Fulton}. 

\subsubsection{\bf Refined intersection product \cite[\S 6]{Fulton} } \label{refined_intersection}
Let $\iota: X \to M$ be a regular embedding of codimension $d$ in $\Sch$, and denote the normal bundle by $N_X M$. 
Let $g: V \to M$ be a proper morphism from a purely $k$-dimensional scheme $V$. 
Take the scheme theoretic preimage $W:=g^{-1}(X)=X\times_M V$:  
$$
\xymatrix{
W\ar[r]^{\iota'} \ar[d]_{g'}& V\ar[d]^g \\
X \ar[r]_\iota & M
}
$$ 
We denote the pullback bundle by $N:=g'^*N_X M$ via $g': W \to X$ and the zero section by $s_W: W \to N$. 
Let $\mathcal{J}$ be the ideal sheaf of $W$ in $V$, then the {\em normal cone to $W$ in $V$} is defined by 
$C=C_W V={\rm Spec}\left(\bigoplus_{n=0}^\infty \mathcal{J}^n/\mathcal{J}^{n+1}\right)$ 
with the projection $\pi: C \to W$. 
It is indeed a subcone of the bundle $N$ and determines a $k$-cycle $[C]$ in $N$. 
We then define the {\em intersection product of $V$ by $X$ on $M$}, denoted by $X\cdot V$, to be 
the rational equivalence class 
$$X\cdot V:=s_W^*[C] \;\; \in \CH_{k-d}(W).$$
Moreover, if $[C] = \sum_{i=1}^r m_i[C_i]$, where $C_i$ is an irreducible component of $C$ with the geometric multiplicity $m_i$, we write the support of $C_i$ by $Z_i=\pi(C_i) \subset W$, and 
put $\alpha_i:=(s_W|_{Z_i})^*[C_i]$ in $\CH_{k-d}(Z_i)$.  
It immediately follows from the construction that 
$$X \cdot V = \sum_{i=1}^r m_i \alpha_i \in \CH_{k-d}(W),$$ 
called the {\em canonical decomposition}. 
Furthermore, for any closed subset $Z$ of $W$, set 
$$(X\cdot V)^Z:=\sum_{Z_i \subset Z} m_i \alpha_i \;\; \in \CH_{k-d}(Z),$$
called the {\em part of $X \cdot V$ supported on $Z$}. 

One way to refine the intersection product is to use a section $s_i$ of the bundle $N|_{Z_i}$ other than the zero section. 
Then, $s_i^*[C_i]$ is a well-defined class on $s_i^{-1}(C_i) \subset Z_i$, which refines $\alpha_i$ 
(i.e., $s_i^*[C_i]$ is sent to $\alpha_i$ by the pushforward of the inclusion to $Z_i$). 
If we have enough sections of $N$, for any closed subset $Z \subset W$, 
we may generically choose a section $s_Z$ on $Z$ so that $s_Z^*[C]$ is a well-defined cycle of dimension $k-d$ which represents $(X\cdot V)^Z$. 

\subsubsection{\bf Refined Gysin map \cite[\S 6]{Fulton}} \label{refined_gysin}
The {\em Gysin homomorphism induced from $\iota$} is defined by 
$$\iota^*: \CH_k(M) \to \CH_{k-d}(X), \quad f^*([V]) := X \cdot V$$ 
for any $k$-cycles $[V] \in \CH_k(Y)$. 
More generally, if $g: M' \to M$ is an arbitrary morphism, 
take the fiber square 
$$
\xymatrix{
X' \ar[r]^{\iota'} \ar[d]_{g'}& M' \ar[d]^g \\
X \ar[r]_\iota & M
}
$$
i.e., $X'=g^{-1}(X)=X\times_M M'$. We define the {\em refined Gysin homomorphism} 
$$\iota^!: \CH_k(M') \to \CH_{k-d}(X')$$
by the same formula $f^!([V]) := X \cdot V$. 
It enjoys some natural properties, e.g., 
\begin{itemize}
\item[(a)] If $g$ is proper, and $\alpha \in \CH_k(M')$, then 
$$\iota^*g_*(\alpha)=g'_*\iota^!(\alpha) \;\; \in \CH_{k-d}(X).$$
\item[(b)] If $g$ is flat of relative dimension $n$, and $\alpha \in \CH_k(M)$, then 
$$g'^*\iota^*(\alpha)=\iota^!g^*(\alpha) \;\; \in \CH_{k+n-d}(X').$$
\end{itemize}
See \cite[\S 6.2]{Fulton} for more information.

\subsubsection{\bf Dynamic intersection \cite[\S 11]{Fulton}} \label{dynamic_intersection}
Let $T$ be a smooth curve equipped with a given point $0 \in T$, and let $T^*=T-\{0\}$. 
For a scheme $\mathcal{W}$ over $T$, the fiber of $\mathcal{W}$ over $t \in T$ is denoted by $W_t$. Set 
$\mathcal{W}^*=\mathcal{W}-W_0=\mathcal{W}\times_T T^*$, and 
let $\mathcal{W}'$ be the closure of $\mathcal{W}^*$ in $\mathcal{W}$. 
Then the limit set 
$$\lim_{t\to 0} W_t \subset \mathcal{W}'$$
is defined to be the fiber of $\mathcal{W}'$ over $0$. 
For a cycle $\alpha = \sum n_i [\mathcal{Y}_i] \in \CH_{k+1}(\mathcal{W}^*)$ where $\mathcal{Y}_i$ are irreducible subvarieties of $\mathcal{W}^*$, we also define the limit class 
$$\lim_{t \to 0} \alpha_t :=\sum n_i [(\overline{\mathcal{Y}_i})_0]\;\;\in \CH_{k-d}\left(\lim_{t\to 0} W_t\right)$$ 
where $\overline{\mathcal{Y}_i}$ is the closure of $\mathcal{Y}_i$ in $\mathcal{W}'$, and $(\overline{\mathcal{Y}_i})_0$ is the scheme-theoretic fiber of $\overline{\mathcal{Y}_i}$ over $0$. 

Consider the fiber square of schemes over $T$: 
$$
\xymatrix{
\mathcal{W} \ar[r] \ar[d] & \mathcal{V} \ar[d]^g\\
\mathcal{X} \ar[r]_{\bar{\iota}\;\;\;} & \mathcal{M}
}
$$
Assume that $\bar{\iota}$ itself and $\bar{\iota}_t: X_t \to M_t$ are regular embeddings of codimension $d$, and that 
$\mathcal{V}$ is a scheme with a closed embedding $g$ of relative dimension $k$. 
We then define the {\em limit intersection class} 
$$\lim_{t\to 0} (X_t \cdot V_t):=\lim_{t\to 0}(\mathcal{X}^*\cdot \mathcal{V}^*)_t \;\; \in \CH_{k-d}\left(\lim_{t \to 0} X_t \cap V_t\right).$$
If $\mathcal{V}$ is flat over $T$, the inclusion of $\lim_{t \to 0} X_t \cap V_t$ to $X_0 \cap V_0$ maps 
the limit intersection class to $X_0\cdot V_0$ \cite[Cor.11.1]{Fulton}. 
For instance, the simplest case is $\mathcal{M}=M\times_k T$ and $\mathcal{V}=V\times_k T$. 

\

\section{Hilbert schemes}

As seen above, the $\K$-classification of finitely determined map-germs $(\bA^m, 0) \to (\bA^l, 0)$ 
is nothing but classifying ideals in $\Ost_{\bA^m,0}$ of finite colength having $l$ generators up to isomorphisms of the source space. It is very related to a certain invariant stratification of the Hilbert schemes of $0$-dimensional subschemes.

\subsection{Hilbert scheme of points}
Let $p\in X$ and $I_p \subset \Ost_{X, p}$ an ideal. The Hilbert-Samuel function of $\Ost_{X,p}/I_p$ is defined to be the sequence $T(I):=(t_0, t_1, \cdots)$ where $$t_j=\dim_k \frac{\mm_{X,p}^j}{(I_p \cap \mm_{X,p}^j)+\mm_{X,p}^{j+1}}$$and $\dim_k \Ost_{X,p}/I_p =\sum t_j $ is called the colength of $I_p$. For an ideal sheaf $I$ on $X$ supported at closed points $p_1, \cdots , p_r$, the colength of $I$ means the sum of colength of the stalk $I_{p_i} \subset \Ost_{X, p_i}$. 
Let $X^{[n]}$ denote the Hilbert scheme parameterizing all ideal sheaves on $X$ 
of colength $n$, and $\Hilb^n(\Ost_{X, p})$ the punctual Hilbert scheme which parameterizes the ideals of $\Ost_{X, p}$ of colength $n$ \cite{Iarrobino, Briancon, Gaffney}. 

According to Grothendieck's general framework, $X^{[n]}$ is the projective scheme which represents the contravariant functor  $\Sch \to \mbox{\bf Sets}$ which assigns to every $U \in \Sch$ the set of closed subschemes $Z \subset U \times X$ so that the first factor projection $Z \to U$ is a flat morphism whose fiber are $0$-dimensional subschemes of $X$ of length $n$. In fact, there exists the universal scheme $\Zcal_X^n$ on $X^{[n]}$, i.e., it is a closed subscheme of $X^{[n]}\times X$ such that the first factor projection 
$$\pi:\Zcal_X^n \to X^{[n]}$$ 
is a flat family of subschemes of $X$ of length $n$ which possesses the following universal property: 
if $Z \subset U \times X$ is a flat family parametrized by $U$, then there exists a unique morphism $\varphi_Z: U \to X^{[n]}$ such that $Z=\varphi_Z^*\Zcal_X^n \, (=U \times_{X^{[n]}} \Zcal_X^n)$. 

We denote by $z_I$ an element of $X^{[n]}$ corresponding to the ideal sheaf $I$, i.e., the fiber over $z_I$ in $\Zcal_n$ is $\Spec(\Ost_X/I)$. 
At $z_I \in X^{[n]}$ so that the support of $I$ is one point $p \in X$ (denote by $I_p$ the stalk at $p$), 
the Zariski tangent space is given by 
$$T_{z_I}X^{[n]}=\Hom_{\Ost_{X,p}}(I_p, \Ost_{X,p}/I_p).$$
At $z_I$ supported by distinct $p_1, \cdots, p_r \in X$, the Zariski tangent space $T_{z_I}X^{[n]}$ is defined to be the direct sum of the Zariski tangent space of $X^{[n_i]}$ at $z_{I_{p_i}}$. 
It is known that when $X \in \Sm$ and $m=\dim X\le 2$, $X^{[n]}$ is smooth for any $n$. 

Suppose that we are given a vector bundle $E \to X$ of rank $a$.  
Let $\pi: \Zcal_X^n \subset  X^{[n]}\times X\to X^{[n]}$ be the universal family  
and $\rho: \Zcal_X^n \to X$ the second projection. 
Since $\pi$ is flat, 
$$E^{[n]}:=\pi_*\rho^*E \to X^{[n]}$$
is locally free, so it is a vector bundle of rank $na$. 
We call $E^{[n]}$ the {\em tautological bundle associated to $E$} over $X^{[n]}$. 
By definition, the fiber $E^{[n]}_{z_I}$ at $z_I$ is $E$ restricted to $\Spec(\Ost_X/I)$; 
for the simplicity, in case that $I$ is concentrated at one point $p \in X$, 
we may think of the fiber as 
$$E^{[n]}_{z_I}=\Hom_k(E^*_p, \Ost_{X,p}/I_p), $$ 
where $E^*$ is the dual bundle. 
For instance, for the tangent bundle $TX$ ($m=\dim X$), 
the $mn$-bundle $(TX)^{[n]}$ is defined over $X^{[n]}$. 
The fiber at a point $z_I$ is denoted by $(TX)^{[n]}_{z_I}$, e.g., 
if $I$ is concentrated at $p \in X$, then 
$$(TX)^{[n]}_{z_I}=\Hom_k(\mathfrak{m}_{X,p}/\mathfrak{m}_{X,p}^2, \Ost_{X,p}/I_p)$$
(this notation should be clearly distinguished from $T_{z_I}X^{[n]}$). 

\begin{rem}\label{flattner}\upshape {\bf (Local flattner)} 
A neighborhood structure of $X^{[n]}$ around a point $z_I$ is explicitly described in terms of the {\it local flattener of a map-germ} \cite{Iarrobino, Briancon, Galligo} 
(the notion was originally introduced by Hironaka, Lejeune and Teissier). 
The local flattener of a map-germ $\varphi: (V, z) \to (S, s_0)$ is a flat morphism $\varphi_{Flat}: V_F:=V\times_S S_F \to S_F$ induced by an inclusion $S_F \hookrightarrow S$ such that for any flat morphism $\psi: V' \to S'$ which commutes with $\varphi$ 
by some $S' \to S$, there exists a morphism $S' \to S_F$ by which $\psi$ is obtained from $\varphi_{Flat}$: 
$$
\xymatrix{
V\times_S S_F\ar[r]\ar[d]_{\varphi_{Flat}} & V  \ar[d]^\varphi|(.58)\hole &\\
S_F \ar[r]^{incl} & S & V' \ar[llu] \ar[lu]\ar[d]^\psi\\
& & S'  \ar[llu] \ar[lu]\\
}
$$
For the simplicity, let $z_I \in X^{[n]}$ be concentrated at a single point $p \in X$, and take $h_1, \cdots, h_s \in \Ost_{X,p}$ the standard basis of the ideal $I$ according to some fixed order. Choose $e_1, \cdots, e_n \in \Ost_{X,p}$ which represent a basis of the vector space $\Ost_{X,p}/I$, and set $$H_i=h_i+\sum_{j=1}^n a_{ij} e_j \qquad (a_{ij} \in \bA, \;\; 1\le i \le s).$$ Denote the space of parameters $(a_{ij})$ by $S:=\bA^{sn}$. Define $$V:=V(\langle H_1, \cdots, H_s\rangle) \subset X \times S$$ and let $\varphi: V \to S$ be the projection to the second factor around $z=(z_I, 0) \in V$. Then the local flattener $\varphi_{Flat}: V_F \to S_F$ of $\varphi$ is defined to be the restriction of $\varphi$ on the preimage $V_F$ of the largest subset $S_F$ of $S$ so that $\varphi$ is flat over $S_F$. It is known that $S_F$ at the origin is locally isomorphic to $X^{[n]}$ at $z_I$ and the Zariski tangent space $T_{z_I}X^{[n]}$ is identified with $T_0S=\bA^{sn}$ 
\cite[p.302]{Iarrobino}. \end{rem}

\subsection{Hilbert extension map}\label{relative_hilbert}
Suppose that a map $\iota: X \to M$ in $\Sch$ satisfies that the restriction of the map 
$$(id \times \iota): X^{[n]} \times X \to X^{[n]} \times M$$ 
to the universal family $\Zcal_X^n$ 
is a flat family of points in $M$ parameterized by $X^{[n]}$. 
Then, by the universality of $\Zcal_M^n$, $\iota$ uniquely induces 
a morphism 
$$\iota^{[n]}:X^{[n]}\to  M^{[n]}$$ 
such that $(\iota^{[n]})^*(\Zcal_M^n)=(id \times \iota)(\Zcal_X^n)$ \cite[Lem.\,1.1]{Gaffney}.

\begin{prop}\label{emb}  
If $\iota: X \to M$ is a regular embedding, then $\iota^{[n]}:X^{[n]}\to  M^{[n]}$ is also a regular embedding. 
Let $\nu$ be the normal bundle for $\iota$, then 
the normal bundle for $\iota^{[n]}$ is the tautological bundle $\nu^{[n]}$. 
\end{prop}

\proof 
Assume that $\iota: X \hookrightarrow M$ is a regular embedding of codimension $l$.  
Let $\Ical_X$ denote the ideal sheaf defining $X$,  $\pi: \Zcal_M^n \to M^{[n]}$ the universal flat family  
and $\rho: \Zcal_M^n \to M$ the second projection. 
Then $\iota^{[n]}(X^{[n]})$ is defined as the zero locus of the composed sheaf homomorphism between locally free sheaves over $M^{[n]}$: 
$$\pi_*\rho^*\Ical_X \hookrightarrow \pi_*\rho^*\Ost_M \to \pi_*\Ost_{\Zcal_M^n}.$$ 
Locally we may take $l$ generators of $\Ical_X$. 
Since the dimension of every stalk of $\pi_*\Ost_{\Zcal_M^n}$ constantly equals $n$,  $\iota^{[n]}(X^{[n]})$ is locally defined by a regular sequence of exactly $ln$ generators, i.e., $\iota^{[n]}$ is a regular embedding (cf. \cite[B.7]{Fulton}). 
The second statement immediately follows,  
e.g., the normal at $z_I \in X^{[n]}$ is 
$$N_{z_I}=\Hom_k(\Ical_{X, p}/\Ical_{X, p}^2, \Ost_{X, p}/I_p),$$
when $I$ is concentrated at a point $p \in X$. 
\qed

\

\begin{rem}\upshape Here is a closer look at the proof of Proposition \ref{emb} in a down-to-earth way using local flattner: let $z_I \in X^{[n]}$ be concentrated at a single point $p \in X\, (\subset M)$ and notations are the same as in Remark \ref{flattner}. Let $u_1, \cdots, u_l \in \Ost_{M, p}$ be a regular sequence defining $X$ around $p$. Then $\iota^{[n]}(z_I)$ corresponds to the ideal $$I':=\langle h_1, \cdots, h_s, u_1, \cdots, u_l\rangle \subset \Ost_{M, p}$$ Obviously, $\Ost_{M, p}/I'=\Ost_{X, p}/I$.   According to Remark \ref{flattner}, from those $(s+l)$ generators of $I'$, we obtain the local flattener $\tilde{\varphi}: \tilde{V}_F \to \tilde{S}_F$, where $\tilde{S}_F \subset \bA^{(s+l)n}$ expresses the neighborhood structure of $M^{[n]}$ at $z_{I'}$. By the construction in Remark \ref{flattner}, it is easy to find an automorphism of the product space $\bA^{sn} \times \bA^{ln}$ so that $$\tilde{S}_F\simeq S_F \times \bA^{ln}$$ where $S_F$ is locally isomorphic to $X^{[n]}$ at $z_I$ and the second factor $\bA^{ln}$ corresponds to the normal to the embedding $\iota^{[n]}$. \end{rem}

Given a morphism $f: X \to Y$, 
consider the embedding $\iota: X \to M=X \times Y$ of the graph of $f$ . 
The scheme-theoretic image of $\Spec(\Ost_X/I)$ by the graph map is $\Spec(\Ost_{X\times Y}/I')$ 
with $I'=\langle I, I_{\Gamma}\rangle$ where $I_{\Gamma}$ is the ideal sheaf defining the graph. 
Then $\iota^{[n]}(z_I)=z_{I'}$. We denote it by $f^{[n]}$. 

\begin{definition}\upshape \label{extension}
{\bf (Hilbert extension map \cite{Gaffney, DL})} 
For a morphism $f: X \to Y$ in $\Sch$, 
the graph map $X \to X \times Y$ induces an embedding $f^{[n]}: X^{[n]}\to (X\times Y)^{[n]}$.  
We call them the {\em $n$-th Hilbert extension map of $f$}. 
\end{definition}

It is notable that the Hilbert extension map is defined even for morphisms between singular varieties. 
Further, if the target space $Y$ is smooth and of dimension $p$ (while $X$ can be singular), the graph map is a regular embedding with codimension $p$, therefore Proposition \ref{emb}, implies  

\begin{cor}\label{Hextension}
For a morphism $f: X \to Y$ in $\Sch$ with $Y \in \Sm$, 
the $n$-th Hilbert extension map $f^{[n]}$ is a regular embedding. 
\end{cor}

\begin{rem}\upshape
One may wonder why nobody paid attention so far to the simple fact of Corollary \ref{Hextension}. Perhaps, the reason would be that in that time there had been no tool to explore further on highly singular spaces, Hilbert schemes, and the fact did not help at all. 
\end{rem}

\subsection{Relative Hilbert scheme}
Assume $f: X \to Y$ to be in $\Sm$. 
As seen above, although Hilbert schemes $X^{[n]}$ and $(X\times Y)^{[n]}$ are singular in general,  $f^{[n]}$ is always  normally non-singular;   
the Zariski tangent space of $(X\times Y)^{[n]}$ at $z_{I'}=f^{[n]}(z_I)$ canonically splits:
$$T_{x_{I'}} (X\times Y)^{[n]} = T_{z_I}X^{[n]} \oplus N_{z_{I'}}$$
where the second factor is the normal vector space 
$$N_{z_{I'}}=(TY)^{[n]}_{z_{I'}}=\bigoplus_{i=1}^r \Hom_{k}\left(\frac{\mathfrak{m}_{Y, q_i}}{\mathfrak{m}_{Y, q_i}^2}, \frac{\Ost_{X, p_i}}{I_{p_i}}\right) \simeq  \bA^{ln},$$ 
where $I$ is supported by distinct closed points $p_1, \cdots, p_r \in X$ and $q_1, \cdots, q_r \in Y$ ($q_i=f(p_i)$), and $l=\dim Y$, $n=\sum_{i=1}^r \dim\Ost_{X, p_i}/I_{p_i}$. 

In more detail, 
the Hilbert extension map is a section of a smooth morphism with relative dimension $pn$ over $X^{[n]}$. 
Let $M=X\times Y$ and 
$$\pi_X: \Zcal_X\subset M^{[n]}\times M \to M^{[n]}\times X$$
the universal scheme of $M^{[n]}$ with the natural projection of the second factor. 
Then, there is a maximal Zariski open subset $V_X$ of $M^{[n]}$ such that $\pi_X(\Zcal_M)$ is flat over $V_X$. 
By the universality, there is a morphism 
$$\Phi: V_X \to X^{[n]}$$ 
such that $\pi_X(\Zcal_M)$ is the pullback of $\Zcal_X$ by $\Phi$. 
It is clear to see that 

\begin{prop}\label{V_X}
The map $\Phi$ is a smooth morphism in $\Sm$ whose relative tangent bundle is the tautological bundle $(TY)^{[[n]]}|_{V_X}$, and the Hilbert extension map $f^{[n]}$ is a section of $\Phi$, i.e.,  
$$f^{[n]}: X^{[n]} \to V_X \subset (X\times Y)^{[n]} \quad (\Phi\circ f^{[n]}=id).$$ 
\end{prop}

There is a subvariety of $ (X\times Y)^{[n]}$ (actually, a subvariety of $V_X$), denoted by $\mDelta^{[n]}$, consisting of ideal sheaves on $X \times Y$ which are supported on finite sets of the form $S \times q$ 
and contains $\mm_{Y,q}$ for this $q$, i.e., each stalk is of the form $\langle I_{p_i}, \mm_{Y,q}\rangle$. 

\begin{prop}\label{diag}
{\rm  \cite[Prop.1.4]{Gaffney}}
The morphism 
$$\iota_\mDelta: X^{[n]} \times Y \to (X\times Y)^{[n]}, \;\; \iota_\mDelta(z_I, q):= z_{\langle I,\mm_{Y,q}\rangle}$$ 
is isomorphic onto $\mDelta^{[n]}$. Further, for a morphism $f: X \to Y$ in $\Sm$,  
$$z_{I'}=f^{[n]}(z_I) \in \mDelta^{[n]}\;\; \Longleftrightarrow \;\; 
f^*(\mm_{Y,q}) \subset I_{p_i} \;\;\; (\forall\, p_i \in S)$$ 
where $I$ is supported at $S$ with $q=f(p_i)$ for any $p_i \in S$. 
\end{prop}

\proof 
The first assertion repeats the definition of $\mDelta^{[n]}$. 
The condition that $z_{I'}=f^{[n]}(z_I)\in \mDelta^{[n]}$ is equivalent to that 
in $\Ost_{X\times Y, (p_i, q)}$, 
\begin{eqnarray*}
I'_{p_i, q}&=&\langle I_{p_i}, y_1-f_1(x), \cdots, y_l-f_l(x)\rangle\\
&=&\langle I_{p_i}, y_1, \cdots, y_l\rangle
\end{eqnarray*}
where $(y_1, \cdots, y_l)$ are local coordinates of $Y$ centered at $q$ and $f_j(x):=y_j\circ f(x)$. 
\qed 

\

\begin{definition}\label{relativeHilbert}\upshape
{\bf (Relative Hilbert scheme)} 
For a morphism $f: X \to Y$ in $\Sm$, 
we define the {\em relative Hilbert scheme $X^{[n]}(f)$} 
to be the scheme-theoretic preimage $(f^{[n]})^{-1}(\mDelta^{[n]})$, 
i.e., the reduced scheme of the fiber product of morphisms $\iota_\mDelta$ and $f^{[n]}$. 
$$
\xymatrix{
X^{[n]}(f) \ar[r]\ar[d] & \;\; \mDelta^{[n]}  \ar[d] \\
X^{[n]} \ar[r]_{f^{[n]}\;\;\;} & (X\times Y)^{[n]} 
}
$$
\end{definition}
By definition, the relative Hilbert scheme $X^{[n]}(f)$ parametrizes $n$ points of the fiber of $f: X \to Y$. 
Namely, let $$\bar{f}: X^{[n]}(f) \to Y$$ 
denote the projection to the second factor of $\mDelta^{[n]}=X^{[n]} \times Y$, 
then the fiber of $\bar{f}$ over $q \in Y$ is the modulus of $0$-dimensional subschemes lying on $f^{-1}(q)$.

\begin{exam}\label{tri} 
\upshape 
{\bf (Desingularization of triple-point locus)} 
In case of $\kappa > 0$, the relative Hilbert scheme $X^{[n]}(f)$ is regarded as a sort of `partial' desingularization of $n$-th multiple-point locus of $f$. 
In fact, it is known \cite{DL, Gaffney} that for locally stable $f: X \to Y$ in $\Sm$, 
$X^{[3]}(f)$ gives a resolution of singularities of the triple-point locus in the target $Y$, 
which is isomorphic to the resolution constructed in a different way using blow-ups \cite{Ronga84}.  
If $f$ is locally stable and $\dim \ker df_x \le 2$, i.e.,  
$\Sigma^3(f)=\emptyset$,  then $X^{[n]}(f)$ for any $n$ gives rise to 
a resolution of singularities of the $n$-th multiple-point locus of $f$ \cite[Cor.3.7]{Gaffney}. 
\end{exam}

\subsection{Differential of Hilbert extension map} \label{DH}
Due to Grothendieck, the differential of a morphism in $\Sch$ at a point is defined as a linear operator between the Zariski tangent spaces. We apply this to the map $f^{[n]}$ and it is easy to see 

\begin{prop} \label{differential}
{\rm \cite[Thm.1.3]{Gaffney}} 
Let $f: X \to Y$ be in $\Sm$, and suppose that $I$ is concentrated at a point $p \in X$ such that $\Spec(\Ost_{X,p})/I_p$ lies on the fiber of $f$. 
Then, the differential of $f^{[n]}$ at $z_I$ between the Zariski tangent spaces 
$$Df^{[n]}: T_{z_I} X^{[n]} \to T_{z_{I'}}(X\times Y)^{[n]}$$
is given by the formula
$$(Df^{[n]}\varphi)(h)=\varphi(h \circ \iota)$$
where $\varphi \in \Hom_{\Ost_{X, p}}(I_p, \Ost_{X, p}/I_p)$, $h \in I'=\langle I, I_\Gamma \rangle$, 
and $\iota: X \to X\times Y$ is the graph map. 
\end{prop}

In the above setting, a vector field-germ (derivation) $v$ on $X$ at $p$ defines a particular element of 
$T_{z_I} X^{[n]}=\Hom_{\Ost_{X, p}}(I_p, \Ost_{X, p}/I_p)$, i.e., 
$v: \varphi\mapsto v(\varphi)$ modulo $I$. 
That yields a canonical inclusion
$$T_pX \hookrightarrow T_{z_I} X^{[n]}.$$
Now, take local coordinates $(x_1, \cdots, x_m)$ and $(y_1, \cdots, y_l)$ centered at $p$ and $q=f(p)$, respectively, and put $f_j(x)=y_j\circ f(x)$. 
Then, the image of $v=a(x) \frac{\rd}{\rd x_i}$ with $a(x) \in \Ost_{X,p}$ via the differential $Df^{[n]}$ has the factor in the normal vector space 
$$(TY)^{[n]}_{z_{I'}}=\Hom_{k}\left({\mathfrak{m}_{Y,q}}/{\mathfrak{m}_{Y,q}^2}, {\Ost_{X, p}}/I_p\right)$$ 
which is determined by 
$$
Df^{[n]}\left(a(x) \frac{\rd}{\rd x_i}\right)(y_j) = a(x)\frac{\rd f_j}{\rd x_i}(x)\;\; \mbox{modulo $I_p$}.
$$
This is immediately extended to multi-singularity case. In other words, via the tautological construction, 
we have the vector bundle map 
$$(df)^{[n]}: (TX)^{[n]} \to (f^*TY)^{[n]} \;\; \mbox{over}\;\;  X^{[n]}$$ 
induced from the differential $df: TX \to f^*TY$, and actually 
$Df^{[n]}$ is nothing but the restriction to $\bigoplus_{p \in S} T_{p} X$ of $(df)^{[n]}$ 
at $z_I$ where $I$ is supported at $S$. 
Also it is interpreted in terms of the Thom-Mather theory (\S \ref{TM}) as follows (cf. \cite[p.284]{Gaffney}). 
Note that germs of sections of $(f^*TY)^{[n]}$ at $z_I$ form the $\Ost_{X,S}$-module 
$$\frac{\theta(f)_S}{I.\theta(f)_S}=(\Ost_{X,S}/I)\otimes_{\Ost_{X,S}} \theta(f)_S,$$ 
and also that 
$f^*(\mathfrak{m}_{Y,q}) \subset I$ and 
$\omega f(\mathfrak{m}_{Y,q} \theta_{Y,q}) \subset f^*(\mathfrak{m}_{Y,q})\theta(f)_S$.

\begin{lem}\label{tf}
The restriction of the differential 
$$Df^{[n]}: \bigoplus_{p \in S} T_{p} X\to (TY)^{[n]}_{z_{I'}}$$
is equivalent to the $\Ost_{X, S}$-module homomorphism $tf$ modulo $I.\theta(f)$ 
$$
tf: \theta_{X, S} \to \frac{\theta(f)_S}{I.\theta(f)_S}= \frac{\theta(f)_S}{(I+f^*(\mathfrak{m}_{Y,q}))\theta(f)_S+\omega f(\mathfrak{m}_{Y,q} \theta_{Y,q})}. 
$$ 
\end{lem}


\subsection{Ordered Hilbert scheme}\label{ordered_H}
For our purpose, more suitable is the {\it ordered Hilbert scheme} defined as follows (cf. Rennemo \cite{Rennemo}). 

Let $S^nX$ be the $n$-th symmetric product of $X$, i.e., 
$$S^nX:=
X^n/\mathfrak{S}_n,  \quad X^n=\overbrace{X\times \cdots \times X}^{n}$$
where $\mathfrak{S}_n$ is the $n$-th symmetric group. 
It is also the moduli scheme parameterizing effective $0$-dimensional cycles of degree $n$ in $X$; every element of $S^nX$ is written as $\sum_{i=1}^r n_i [p_i]$, where $p_i \in X$,  $n_i \in \N$ and $\sum_{i=1}^r n_i=n$. 
The Hilbert scheme $X^{[n]}$ is regarded as a `partial desingularization' of $S^nX$; 
there exists a morphism 
$$\pi^{HC}: X^{[n]} \to S^nX, \;\; z_I \mapsto \sum_{p \in X} (\dim \Ost_{X,p}/I_p)\, [p],$$ 
called the {\it Hilbert-Chow morphism}. 
We see that the fiber of $\pi^{HC}$ at $n[p] \in S^nX$ is the punctual Hilbert scheme $\Hilb^n(\Ost_{X,p})$, and the fiber at $\sum_{i=1}^r n_i [p_i]$ is isomorphic to the Cartesian product of $\Hilb^{n_i}(\Ost_{X,p_i})$. 
Especially, the general fiber consists of one point, so $\pi^{HC}$ is birational. 
It is known that for a smooth surface $X$ ($m=\dim X=2$),  $\pi$ is a crepant resolution of $S^nX$ for any $n$. 

The {\em ordered} Hilbert scheme $X^{[[n]]}$ is defined by the fiber product of the Hilbert-Chow morphism and the canonical projection $X^n \to S^nX$: 
$$
\xymatrix{
X^{[[n]]} \ar[r]\ar[d] & X^{[n]} \ar[d]\\
X^n \ar[r] & S^nX 
}
$$
We denote a point of $X^{[[n]]}$ by $(x, z_{I})$ with 
$$x=(p_1, \cdots, p_n) \in X^n, \quad z_I \in X^{[n]}$$ 
where $p_i$ may be duplicated and $\pi^{HC}(z_I)=\sum [p_i]$.   
Note that $X^{[[n]]}$ is a compactification of the configuration space, i.e., 
it has an open dense subset isomorphic to the configuration space of ordered $n$-tuples of points on $X$ 
$$F(X, n):=\{\,(p_1, \cdots, p_n) \in X^n\, |\, p_i\not=p_j \;(i\not=j)\, \}.$$ 
In particular, when $n=2$,  $X^{[[2]]}$ coincide with the blow-up $\mbox{Bl}(X\times X)$ of $X \times X$ along the diagonal  (cf. \cite{FM}). 

Basic properties of non-ordered Hilbert schemes are naturally transported to ordered ones as follows. 

\begin{itemize}
\item {\bf (Regular embedding)} 
Given a regular embedding $\iota: X \hookrightarrow M$, 
we set  
$$\iota^{[[n]]}:X^{[[n]]}\to  M^{[[n]]}, \;\;\;\;
\iota^{[[n]]}(p_1, \cdots, p_n, z_I):= (\iota(p_1), \cdots, \iota(p_n), \iota^{[n]}(z_I)).$$ 
This is also a regular embedding as well as Lemma \ref{emb}. In fact,  the sheaf homomorphism over $M^{[n]}$ in the proof of Lemma \ref{emb} is pulled back via the projection $M^{[[n]]}\to M^{[n]}$ and the zero locus in $M^{[[n]]}$ exactly coincides with $X^{[[n]]}$. 
\item  {\bf (Hilbert extension map)} 
For a morphism $f: X \to Y$ in $\Sm$, we have the ordered Hilbert extension map 
$$f^{[[n]]}: X^{[[n]]} \to (X\times Y)^{[[n]]}$$
and the ordered relative Hilbert scheme 
$$\bar{f}: X^{[[n]]}(f) \to Y$$ 
as well. 
The normal bundle of the embedding $f^{[[n]]}$ is isomorphic to 
the pullback of that of non-ordered extension map $f^{[n]}$, and thus 
the normal at $(x, z_{I})$ is identified with $(TY)^{[n]}_{z_{I'}}$ 
(notations are the same as above). 
Also there is a canonical inclusion 
$$T_{p_i}X \hookrightarrow T_{(x, z_I)} X^{[[n]]}.$$ 
In particular, Lemma \ref{tf} can be read off as for the differential of $f^{[[n]]}$, 
i.e., $Df^{[[n]]}$ with the restriction to $T_{p_i}X$ and the projection to the normal space 
is entirely the same as that of $Df^{[n]}$.
\item  {\bf (Disjoint union)} 
Note that for the disjoint union, there is a canonical decomposition of the non-ordered Hilbert scheme: 
$$(X_1\sqcup X_2)^{[n]}=\bigsqcup_{i+j=n} X_1^{[i]}\times X_2^{[j]}.$$
For the ordered version, we have the decomposition 
$$(X_1\sqcup X_2)^{[[n]]}=\bigsqcup_{I\sqcup J} X_1^{[[I]]}\times X_2^{[[J]]}$$
where the summand runs over all partitions $I \sqcup J=\{1, \cdots, n\}$. 
\end{itemize}

\section{Thom-Kazarian Principle}

From now on, we always consider the {\em ordered} Hilbert schemes. 
Here is our basic diagram which we will frequently use in this paper: 
$$
\xymatrix{
X^{[[n]]}(f) \ar@/^22pt/[rr]^{\bar{f}} \ar[r]\ar[d] & X^{[[n]]} \times Y \ar[r] \ar[d]^\iota & Y\\
X^{[[n]]} \ar[r]_{f^{[[n]]}\;\;\;} & (X\times Y)^{[[n]]} &\\
}
\eqno{(1)}
$$
For our purpose, 
the first factor projection takes an important role, denoted by 
$${\rm pr}_1: X^{[[n]]} \to X, \qquad (p_1, \cdots, p_n, z_I) \mapsto p_1.$$
The map $\bar{f}$ is naturally factored as follows: 
$$
\xymatrix{
&X^{[[n]]}(f) \ar[ld]_{{\rm pr}_1} \ar[rd]^{\bar{f}}&\\
X \ar[rr]_f && Y
}
$$

\subsection{Geometric subsets} \label{gs}
Let $\eta$ be a mono-singularity type of finitely determined germs $(\bA^m, 0) \to (\bA^l, 0)$. 
We often abuse the notation $\eta$ to mean a map-germ of that type.  
As in Introduction (i), put 
$$
I_{\eta}:=\eta^*(\mathfrak{m}_{\bA^{l},0})\Ost_{\bA^m,0}+\mathfrak{m}_{\bA^m,0}^{k_f+1} 
$$
and denote the colength by 
$$n(\eta):=\dim_k \Ost_{\bA^m,0}/I_{f} < \infty.$$
Note that adding the power of the maximal ideal is needed only in case of $\kappa < 0$ (indeed, if $\kappa \ge 0$, finite determinacy implies $I_{\eta}=\eta^*(\mathfrak{m}_{\bA^{l},0})\Ost_{\bA^m,0}$ by Nakayama's lemma). 

We denote by  
$$K_\eta^\circ \subset \Hilb^{n(\eta)}(\Ost_{\bA^m, 0})$$
the subset consisting of all ideals $I \subset \Ost_{\bA^m, 0}$ isomorphic to $I_\eta$. Note that $K_\eta^\circ$ is a locally closed smooth subscheme of the punctual Hilbert scheme; actually it is an orbit under the actions of automorphism-germs $\tau$ of $(\bA^m, 0)$: $\tau. I := \tau^*I$. 
Additionally, if requested, we may take $\eta$ to be an {\em equisingularity type}, which is a modulus of $\K$-type, 
and then $K_\eta^\circ$ is a locally closed smooth subscheme formed by the family of such orbits. 

Let $\underline{\eta}=(\eta_1, \cdots, \eta_r)$ be a multi-singularity type of germs $(\bA^m, 0) \to (\bA^l, 0)$ with $l=m+\kappa$,  and let $n=n(\underline{\eta}):=\sum_{i=1}^r n(\eta_i)$ and $\ell=\ell(\underline{\eta})=\sum_{i=1}^r \ell(\eta_i)$. Let $X$ be a smooth scheme of dimension $m$. 
We set 
$$\Xi^\circ(\underline{\eta})\, (=\Xi^\circ(X;\underline{\eta}))\subset X^{[[n]]}$$ 
to be the constructible subset consisting of $(x, z_I)$ which satisfy the following conditions: 
\begin{enumerate}
\item $I$ is supported at  $r$ distinct closed points $p_1, \cdots,  p_r$ in $X$; 
\item $\Ost_{X, p_i}/I_{p_i}$ is isomorphic to the $\K$-type $\eta_i$  ($1\le i\le r$); 
\item $x=(\delta^{n_1}(p_1), \delta^{n_2}(p_2), \cdots, \delta^{n_r}(p_r)) \in X^n$ where  $\delta^r: X \to X^r$ is the diagonal map $p \mapsto (p, \cdots, p)$ and $n_i=n(\eta_i)$. 
\end{enumerate}
Note that ${\rm pr}_1(x, z_I)=p_1$. By definition, 
$\Xi^\circ(\underline{\eta})$ has the natural projection to the configuration space 
$$\Xi^\circ(\underline{\eta}) \to F(X, r), \quad (x, z_I) \mapsto (p_1, \cdots, p_r).$$
The preimage of a point of $F(X, r)$ is isomorphic to the Cartesian product $\prod_{i=1}^r K_{\eta_i}^\circ$;  especially, $\Xi^\circ(\underline{\eta})$ itself is smooth and this projection is a smooth morphism in $\Sm$. 

\begin{definition}\label{geom_subset}\upshape
We define the {\em geometric subset $\Xi(\underline{\eta})\, (=\Xi(X;\underline{\eta}))$ associated to a multi-singularity type $\underline{\eta}$} to be the Zariski closure of the locally closed smooth subscheme $\Xi^\circ(\underline{\eta})$ in $X^{[[n]]}$. 
\end{definition}

Precisely saying, a multi-singularity type $\underline{\eta}$ is an isomorphism class of quotient algebras associated to multi-germs, 
so it corresponds to the {\em stably} $\K$-equivalence, and we do not need to fix $m, l$, but only $\kappa (=l-m)$. 
Therefore, for any smooth scheme $X'$ of dimension $m'$, $\Xi(X'; \underline{\eta})$ is defined whenever it makes sense, i.e, each entry $\eta_i$ is realized as the $\K$-type of germs $(\bA^{m'}, 0) \to (\bA^{m'+\kappa}, 0)$. 
In particular, we see 

\begin{lem}\label{geom_subset_K}
Given an embedding $X' \subset X$ in $\Sm$,  the preimage of $\Xi(X; \underline{\eta})$ via the embedding $X'^{[[n]]} \hookrightarrow X^{[[n]]}$ coincides with the geometric subset $\Xi(X'; \underline{\eta})$, if it exists. 
\end{lem}

Let $\underline{\eta}=(\eta_1, \cdots, \eta_r)$ be a multi-singularity type. 
For any $I =\{i_1, \cdots, i_k\}$ with $1\le i_1<\cdots<i_k\le r$, we set 
$$\underline{\eta}_I=(\eta_{i})_{i \in I}=(\eta_{i_1}, \cdots, \eta_{i_k}).$$ 
As seen in \S \ref{ordered_H}, 
the (ordered) Hilbert scheme of disjoint union $X_1\sqcup X_2$ has a canonical decomposition. 
Then, the corresponding geometric subset in $(X_1\sqcup X_2)^{[[n]]}$ is naturally decomposed as follows: 

\begin{lem} \label{H_disjoint} 
It holds that 
$$\Xi(X_1\sqcup X_2;\underline{\eta}) = 
\bigsqcup_{I\sqcup J} \Xi(X_1; \underline{\eta}_{I}) \times \Xi(X_2; \underline{\eta}_{J}),$$
where $I$ runs over all subsets of $\{1, \cdots, r\}$ and $J$ is the complement ($I$ and $J$ can be the empty set). 
\end{lem}

\subsection{Multi-singularity loci} 
In the literature, `multi-singularity loci' have been treated using multi-jet extension maps introduced by J. Mather (see Remark  \ref{jet-extension} below). However the meaning is unclear when one treats degenerate maps. Here we present a right definition using Hilbert extension maps. 

There is a natural inclusion: 
$$\Xi(\underline{\eta})\times Y \subset X^{[[n]]} \times Y \hookrightarrow  (X\times Y)^{[[n]]}$$ 
and the composition is also denoted by $\iota_\mDelta$. 

\begin{definition}\label{multi-singularity_loci}\upshape
{\bf (Multi-singularity loci)} 
For $\underline{\eta}=(\eta_1, \cdots, \eta_r)$ 
with $n(\underline{\eta})=n$ and for an arbitrary proper morphism $f: X \to Y$ in $\Sm$, 
we define $C_{\underline{\eta}}(f) \subset X^{[[n]]}(f)$  
to be the reduced scheme of the fiber product 
$$
\xymatrix{
C_{\underline{\eta}}(f) \ar[r]\ar[d] & \Xi(\underline{\eta}) \times Y \ar[d]^{\iota_\mDelta}\\
X^{[[n]]} \ar[r]_{f^{[[n]]}\;\;\;} & (X\times Y)^{[[n]]}
}
\eqno{(1')}
$$
The {\em $\underline{\eta}$-singularity locus} in the source is introduced as the reduced scheme of the image
$$M_{\underline{\eta}}(f):={\rm pr}_1(C_{\underline{\eta}}(f))_{red} \subset X,$$
and the locus in the target is also given by  
$$N_{\underline{\eta}}(f):=\bar{f}(C_{\underline{\eta}}(f))_{red} \subset Y.$$
\end{definition}

\

In Definition \ref{multi-singularity_loci}, it is emphasized that $f$ is arbitrary.  
Even if $f$ is very degenerate and has no $\underline{\eta}$-singular point (i.e., the multi-germ of $f$ ay any finite set $S \subset X$ is not of type $\underline{\eta}$), for instance in case that $f$ is a constant map,   the $\underline{\eta}$-singularity locus $M_{\underline{\eta}}(f)$ still makes sense and non-empty. So, in general, the locus $M_{\underline{\eta}}(f)$ differs from the Zariski closure of the set of $\underline{\eta}$-singular points of $f$. They coincide if $f$ is appropriately generic, see Remark \ref{jet-extension}. 

Note that we do not describe specifically a scheme structure of the multi-singularity locus of a given type. It is a different meaningful question (e.g., Mond-Pellikaan  \cite{MondPellikaan} give the scheme structure of multiple-point loci $M_{A_0^r}(f)$ using Fitting ideals).

\begin{rem}\upshape \label{jet-extension}
{\bf (Multi-jet extension map)} 
For a locally stable map $f: X\to Y$, 
the geometric meaning of the locus $M_{\underline{\eta}}(f)$ is very clear  \cite{MatherIV, MatherV}. 
If $f$ is locally stable and admits $\underline{\eta}$-singular points, 
there is an non-empty locally closed smooth subscheme of $X$, denoted by $M_{\underline{\eta}}^\circ\, (=M_{\underline{\eta}}^\circ(f))$,  consisting of $p_1 \in X$ such that 
there are distinct closed points $p_2, \cdots, p_r \, (\not=p_1)\in X$ satisfying 
\begin{itemize}
\item  $f(p_1)=\cdots =f(p_r)=:q \in Y$;  
\item the germ of $f$ at $S=\{p_1, \cdots, p_r\}$ is locally stable; 
\item the germ of $f$ at $p_i$ is of type $\eta_i$ for each $1\le i\le r$. 
\end{itemize}
Also we have the locus $N_{\underline{\eta}}^\circ$ consisting of such points $q \in Y$; it is a locally closed smooth subscheme of $Y$ and $f$ sends $M_{\underline{\eta}}^\circ$ onto $N_{\underline{\eta}}^\circ$ locally isomorphically. 
By the local stability of $f$ at $S$, it holds that 
$$ \dim M_{\underline{\eta}}^\circ=\dim N_{\underline{\eta}}^\circ=\dim Y-\ell$$
where $\ell=\ell(\underline{\eta})=\sum \ell(\eta_i)$ is the codimension of the type $\underline{\eta}$ (\S \ref{locally_stable}). 
In fact, $M_{\underline{\eta}}^\circ$ and $N_{\underline{\eta}}^\circ$ are precisely defined using the multi-jet extension map 
$${}_rj^kf: F(X, r) \to {}_rJ^k(X, Y)$$
over the configuration space $F(X, r)$ of ordered $r$ points on $X$ \cite{MatherV}. 
The local stablility of $f$ is equivalent to the condition that the map ${}_rj^kf$ is transverse to any multi-$\K$-orbits \cite[Thm.4.1]{MatherV}, and thus the preimage of the $\underline{\eta}$-orbit is smooth and is projected to  $M_{\underline{\eta}}^\circ$ and $N_{\underline{\eta}}^\circ$. 
Then the Zariski closure of $M_{\underline{\eta}}^\circ \subset X$ and $N_{\underline{\eta}}^\circ \subset Y$, respectively, coincide with the loci $M_{\underline{\eta}}(f)$ and $N_{\underline{\eta}}(f)$ in Definition \ref{multi-singularity_loci}. 
A more rigorous treatment will be discussed below. 
\end{rem}

\subsection{Invariance and transversality}\label{invariance}
The smooth subscheme $\Xi^\circ(\underline{\eta})$ is invariant under the action of isomorphism-germs of $X$ at $r$ distinct points\footnote{Actually, in \cite{Gaffney}, the non-ordered version of $\Xi^\circ(\underline{\eta})$ is treated, and called an {\em invariant submanifold} of $X^{[n]}$.}.   
That means the following. 
Pick $(x, z_I) \in \Xi^\circ(\underline{\eta})$, where $I$ is supported on $S=\{p_1, \cdots, p_r\}$ and $x=(\delta^{n_1}(p_1), \cdots, \delta^{n_r}(p_r)) \in X^n$. 
Take any isomorphism-germ 
$$\tau: (X, S') \to (X, S) \quad \mbox{with $\tau(S')=S$},$$ 
and put $x'=(\delta^{n_1}(p'_1), \cdots, \delta^{n_r}(p'_r)) \in X^n$, $\tau(p'_i)=p_i$. 
Then $(x', \tau^*z_I)$ belongs to  $\Xi^\circ(\underline{\eta})$. 
The geometric subset $\Xi(\underline{\eta})$ is also invariant under the action of isomorphism-germs of $X$ at distinct points.

Recall the canonical inclusion $T_pX \hookrightarrow T_{z_I}X^{[n]}$ mentioned just after Proposition \ref{differential}. 
At $(x, z_I) \in \Xi^\circ(\underline{\eta})$, 
any vector field-germ of $X$ at $p_i$ creates an infinitesimal deformation of $I$ at $p_i$ in the corresponding punctual Hilbert scheme, which is now tangent to $K_{\eta_i}^\circ$ by the invariance \cite[Lem.3.3]{Gaffney}. Thus we have

\begin{lem}\label{inv_mfd}
Let $\Xi^\circ(\underline{\eta})$ be as above. At every point $(x, z_I) \in \Xi^\circ(\underline{\eta})$, 
there is a canonical inclusion 
$$\bigoplus_{i=1}^r T_{p_i}X \hookrightarrow T_{(x,z_{I})} \Xi^\circ(\underline{\eta}).$$ 
\end{lem}

Next, let us recall $\Phi: V_X \to X^{[[n]]}$ in Proposition \ref{V_X}  
(here we use ordered Hilbert schmes). 
Write $\Xi^\circ:=\Xi^\circ(\underline{\eta}) \subset X^{[[n]]}$ for short. 
Then 
$$\Phi: V_{\Xi^\circ}:=\Phi^{-1}(\Xi^\circ)  \to \Xi^\circ$$ 
is a smooth morphism in $\Sm$. 
For a given $f:X \to Y$, the Hilbert extension map defines a section 
$$f_{\Xi^\circ}:=f^{[[n]]}|_{\Xi^\circ}: \Xi^\circ \to V_{\Xi^\circ}.$$
Also $\iota_\Delta: \Xi^\circ \times Y \hookrightarrow V_{\Xi^\circ}$ is an embedding in $\Sm$.

We say that $f: X \to Y$ is {\em locally stable along the $\underline{\eta}$-locus} if  
the germ $f: (X, S) \to (Y, q)$ is locally stable for every $q \in N_{\underline{\eta}}(f)$ and finite $S \subset  f^{-1}(q)$.

We have the following key property (cf. {\rm \cite[Thm.3.5]{Gaffney}}):

\begin{prop}\label{transversality}
Suppose that $f:X \to Y$ is locally stable along the $\underline{\eta}$-locus. 
Then the map $f_{\Xi^\circ}$ is transverse to $\iota_\mDelta(\Xi^\circ \times Y)$ in $V_{\Xi^\circ}$. 
In particular, the intersection, denoted by $C_{\Xi^\circ}$, is a smooth subscheme of $\Xi^\circ$. 
\end{prop}

\proof 
As mentioned in Remark \ref{ordered_H}, 
the normal to the embedding $f_{\Xi^\circ}$ at $(x, z_{I})\in C_{\Xi^\circ}$ is identified with 
$$N_{z_{I'}}:=(TY)^{[n]}_{z_{I'}}\;\; \mbox{ (=the relative tangent of $\Phi$ at $(x, z_{I'})$)}$$ 
where $I$ is supported at $S=\{p_1, \cdots, p_r\}$ with $f(p_i)=q \in Y$ and 
$f^*(\mathfrak{m}_{Y,q})\subset I_{p_i}$  for all $i$. 
Note that this normal space contains all tangent vectors of the second factor $Y$ of $\Xi^\circ \times Y$. 
Thus the required transversality condition is fulfilled if the composed linear map of $Df_{\Xi^\circ} (=Df^{[[n]]})$ 
with the canonical inclusion (Lemma \ref{inv_mfd}) and the projection 
$$
\xymatrix{
\bigoplus_{i=1}^r T_{p_i}X\ar[r] & T_{(x,z_{I})} \Xi^\circ   \ar[r]^{\mbox{\tiny $Df_{\Xi^\circ}$}}& N_{z_{I'}}/T_qY
}
$$ 
is surjective. 
Using Lemma \ref{tf} and 
$$\omega f(\theta_{Y,q}) = T_qY \oplus \omega f(\mathfrak{m}_{Y,q}\theta_{Y,q}),$$ 
this is equivalent to the surjectivity of the $\Ost_{X,S}$-homomorphism 
$$
\xymatrix{
tf: \theta_{X, S} \ar[r] & \displaystyle \frac{\theta(f)_S}{(I+f^*(\mathfrak{m}_{Y,q}))\theta(f)_S+\omega f(\theta_{Y,q})} 
}
$$ 
and it follows from the assumption, that is the local stability of $f$ at $S$ (\S \ref{locally_stable}). \qed

\

Let $f: X\to Y$ be locally stable, and put 
$$l=\dim Y, \;\; n=n(\underline{\eta}), \;\; \ell=\ell(\underline{\eta})$$ 
and also $m=\dim X$ and $\kappa=l-m$. 

Take $C_{\Xi^\circ} \subset \Xi^\circ\, (:=\Xi^\circ(\underline{\eta}))$ and $(x, z_{I})\in C_{\Xi^\circ}$ 
as in Proposition \ref{transversality} and its proof. 
Since $f$ is locally stable, there is an open dense subset of $C_{\Xi^\circ}$, say $W$,  consisting of $(x, z_I)$ such that $z_I$ is uniquely determined such as 
$$I_{p_i}=f^*(\mathfrak{m}_{Y, q})\Ost_{X, p_i}+\mathfrak{m}_{X, p_i}^{k_{\eta_i}+1} \quad (1\le i\le r),$$ 
i.e., the germ of $f$ at $p_i$ is of type $\eta_i$. 
Put 
$$M_{\underline{\eta}}^\circ:={\rm pr}_1(W)_{red} \subset X, \quad 
N_{\underline{\eta}}^\circ:=\bar{f}(W)_{red} \subset Y.$$
(see Remark \ref{jet-extension}). Note that  
$$\dim C_{\Xi^\circ}=\dim W=\dim M_{\underline{\eta}}^\circ = \dim N_{\underline{\eta}}^\circ=l-\ell.$$

Although $\Xi(\underline{\eta})$ is determined only by isomorphic types $\eta_i$ of local algebras 
(so the information of $Y$ does not matter), we have the following relation with $l$ and $\ell$: 

\begin{cor}\label{dim_geometric_subset}
The geometric subset $\Xi(\underline{\eta})$ has dimension $ln-\ell$. 
\end{cor}

\proof 
Put $s:=\dim \Xi(\underline{\eta})$. 
Suppose that the germ $f: (X, S)\to (Y,q)$ is of type $\underline{\eta}$ and locally stable, 
and let $I$ denote the corresponding ideal sheaf on $X$ supported at $S$. 
Since $f_{\Xi^\circ}: \Xi^\circ \to V_{\Xi^\circ}$ is transverse to $\Xi^\circ \times Y$ at $f_{\Xi^\circ}(z_I)$, 
we see 
\begin{eqnarray*}
\dim C_{\Xi^\circ}&=&\dim \Xi^\circ + \dim (\Xi^\circ\times Y) - \dim V_{\Xi^\circ}\\
&=&2s+l-(s+ln)\\
&=&s+l-ln.
\end{eqnarray*}
This is equal to $l-\ell$, thus we get $s=ln-\ell$. 
\qed

\subsection{Multi-singularity loci classes} 
Let $\underline{\eta}=(\eta_1, \cdots, \eta_r)$ be a multi-singularity type of germs $(\bA^m,0) \to (\bA^l,0)$. 
Given an arbitrary proper morphism $f: X \to Y$ with $m=\dim X$ and $l=\dim Y$, 
take the fiber square $(1')$. 
Then we can define the intersection product of the regular embedding $f^{[[n]]}$ and the subscheme $\Xi(\underline{\eta})\times Y$ in the sense of Fulton-MacPherson \cite{Fulton}. 
The expected dimension of the intersection is $l-\ell$.  
Thus, we have 
$$X^{[[n]]}\cdot (\Xi(\underline{\eta})\times Y) \in \CH_{l-\ell}(X^{[[n]]}(f)).$$
This is nothing but the image of the class $[\Xi(\underline{\eta})\times Y]$ 
via the {\it refined Gysin map} \cite[\S 6.2]{Fulton}
$$(f^{[[n]]})^!: \CH_*(X^{[[n]]}\times Y) \to \CH_{*-ln}(X^{[[n]]}(f)).$$

\begin{definition}\label{multi-singularity_class}\upshape
{\bf (Source/target multi-singularity loci classes)} 
For an arbitrary proper morphism $f: X \to Y$ of codimension $\kappa$ in $\Sm$, 
we define the {\em source and the target $\underline{\eta}$-singularity loci classes}, respectively, by 
the pushforward of the intersection product of the Hilbert extension map with the geometric subset $\Xi(\underline{\eta})\times Y$: 
$$m_{\underline{\eta}}(f):={\rm pr}_{1*}(X^{[[n]]}\cdot (\Xi(\underline{\eta})\times Y)) \in \CH_{l-\ell}(X)=\CH^{\ell-\kappa}(X),$$
$$n_{\underline{\eta}}(f):=\bar{f}_*(X^{[[n]]}\cdot (\Xi(\underline{\eta})\times Y)) \in \CH_{l-\ell}(Y)=\CH^\ell(Y). $$
\end{definition}

In a finer sense, those classes are actually localized on the $\underline{\eta}$-singularity loci of the source and the target, respectively, i.e., 
$$X^{[[n]]}\cdot (\Xi(\underline{\eta})\times Y) \in \CH_{l-\ell}(C_{\underline{\eta}}(f)),$$  
$$m_{\underline{\eta}}(f) \in \CH_{l-\ell}(M_{\underline{\eta}}(f)), 
\quad 
n_{\underline{\eta}}(f) \in \CH_{l-\ell}(N_{\underline{\eta}}(f)).$$

\begin{rem}\upshape \label{lcim}
{\bf (Local complete intersection morphism)} 
Definition \ref{multi-singularity_class} makes sense even in case that $X$ is singular but equidimension and a proper morphism $f: X \to Y$ in $\Sch$ is a {\em local complete intersection morphism} with smooth $Y$, i.e., $f$ is a composition of a map $M \to Y$ in $\Sm$ and a regular embedding $\iota: X \hookrightarrow M$. 
We define the geometric subset $\Xi(X; \underline{\eta}) \subset X^{[[n]]}$ to be the schematic preimage of $\Xi(M; \underline{\eta}) \subset M^{[[n]]}$ via $\iota^{[[n]]}$ and it yields source and target multi-singularity loci classes for $f: X \to Y$. 
\end{rem}

\begin{rem}\upshape \label{FL}
{\bf (Mono-singularity loci class)} 
Consider the case of $r=1$: let $\eta_1$ be a mono-singularity type of germs $(\bA^m, 0)\to (\bA^l,0)$. Like as degeneracy loci class mentioned in \S \ref{TP}, 
the $\eta_1$-singular loci class of a map $f: X \to Y$ in $\CH^*(X)$ is defined as the intersection product of the jet-extension map $d^kf: X \to J^k(X, Y)$ with the total space of the subbundle with fiber being the closure of the $\K$-orbit of type $\eta_1$ ($k \ge k_{\eta_1}$). 
Comparing the jet space $J^k(m, l)$ with the punctual Hilbert scheme $\Hilb^n(\Ost_{\bA^m,0})$, it is not difficult to see that the class coincides with $m_{\eta_1}(f)$ in Definition \ref{multi-singularity_class}  as a particular case. 
\end{rem}

The following properties are obvious from the definition. 

\begin{lem}\label{iso}
(1) It holds that 
$$f_*m_{\underline{\eta}}(f)=n_{\underline{\eta}}(f).$$
(2) For any isomorphism $\sigma: X' \to X$,  
$$\sigma_*m_{\underline{\eta}}(f\circ \sigma)=m_{\underline{\eta}}(f), 
\quad n_{\underline{\eta}}(f\circ \sigma)=n_{\underline{\eta}}(f).$$  
\end{lem}

The $\underline{\eta}$-singularity loci class satisfies the base-change property in the following sense. 
Suppose that a proper morphism $f: X \to Y$ and 
a morphism $\varphi: Y' \to Y$ are transverse to each other within $\Sm$. 
Then we have the fiber square $X':=X\times_Y Y'$ in $\Sm$ with the same codimension $\kappa(f')=\kappa(f)=\kappa$: 

$$
\xymatrix{
X' \ar[r]^{\varphi'} \ar[d]_{f'} & X \ar[d]^f\\
Y' \ar[r]_{\varphi} & Y
}
$$

\begin{prop}\label{pullback}
Let $\varphi$ and $f$ be transverse to each other, and consider the above diagram.  
Then, the pullback operations $(\varphi')^*$ and $\varphi^*$ preserves the $\underline{\eta}$-singularity loci classes: 
$$m_{\underline{\eta}}(f')=(\varphi')^*(m_{\underline{\eta}}(f)), \quad n_{\underline{\eta}}(f')=\varphi^*(n_{\underline{\eta}}(f)).$$ 
\end{prop}

 \proof 
 The relative Hilbert scheme $X'^{[[n]]}(f')$ associated to $f': X' \to Y'$ parametrizes $0$-dimensional subschemes of $X$ lying on the scheme theoretic fiber $f^{-1}(\varphi(y'))$ where $y' \in Y'$ varies.
Actually, it is identified with the fiber product of $\bar{f}: X^{[[n]]}(f) \to Y$ and  $\varphi: Y' \to Y$, 
and we have the following commutative diagram:  
$$
\xymatrix{
X' \ar@/^18pt/[rrr]^{f'}\ar[d]_{\varphi'}&X'^{[[n]]}(f') \ar[l]^{{\rm pr}_{1}\;\;\;} \ar[r] \ar[d]^{\bar{\varphi}} & X^{[[n]]} \times Y' \ar[d]_{id\times \varphi} \ar[r] & Y'\ar[d]^{\varphi}\\
X \ar@/_18pt/[rrr]_{f}& X^{[[n]]}(f) \ar[l]_{{\rm pr}_{1}\;\;\;} \ar[r] & X^{[[n]]}\times Y \ar[r] & Y
}
$$
In other words, we may consider the regular embedding $\iota: X' \to X \times Y'$ and its Hilbert extension map 
$$\iota^{[[n]]}:X'^{[[n]} \to (X\times Y')^{[[n]]},$$ 
then $X'^{[[n]]}(f')$ is identified with the fiber product of $\iota^{[[n]]}$ and the subscheme 
$$X^{[[n]]} \times Y' \hookrightarrow (X\times Y')^{[[n]]}.$$ 
Note that the normal bundle of $\iota$ is naturally isomorphic to the pullback of $TY$ via $\varphi\circ f'\, (=f \circ \varphi')$. 
Therefore, the normal bundle of $\iota^{[[n]]}$ along $X'^{[[n]]}(f')$  is isomorphic to  
the normal bundle $TY^{[[n]]}$ of $(f)^{[[n]]}$ restricted to $X^{[[n]]}(f)$. 
The base-change $\varphi$ induces the refined Gysin map 
$$\varphi^{!}: \CH_*(X^{[[n]]}(f) ) \to \CH_*(X'^{[[n]]}(f')).$$ 
By the above construction, 
given a geometric subset $\Xi(\underline{\eta}) \subset X^{[[n]]}$, we have 
$$\varphi^{!}(X^{[[n]]}\cdot (\Xi(\underline{\eta})\times Y))=X'^{[[n]]}\cdot (\Xi(\underline{\eta})\times Y').$$ 
Since $\varphi^{!}$ enjoys $\varphi^*\circ \bar{f}_*=\bar{f'}_*\circ \varphi^{!}$ ((a) in \S \ref{refined_intersection}), 
we see 
\begin{eqnarray*}
n_{\underline{\eta}}(f')
&=&\bar{f'}_*(X'^{[[n]]}\cdot (\Xi(\underline{\eta})\times Y'))\\
&=& \bar{f'}_*\circ \varphi^{!}(X^{[[n]]}\cdot (\Xi(\underline{\eta})\times Y))\\
&=& \varphi^*\circ \bar{f}_*(X^{[[n]]}\cdot (\Xi(\underline{\eta})\times Y))\\
&=&\varphi^*(n_{\underline{\eta}}(f)).
\end{eqnarray*}
Also since $(\varphi')^*\circ {\rm pr}_{1*}={\rm pr}_{1*}\circ \varphi^{!}$,   
the naturality of $m_{\underline{\eta}}$ for pullback follows in entirely the same way. 
This completes the proof. 
\qed

\begin{rem}\label{n_pullback}\upshape 
In Proposition \ref{pullback}, we can weaken the assumption. Even in case that $\varphi$ is not transverse to $f$ and the fiber product $X'$ is singular (while $X, Y, Y'$ are smooth, $l'=\dim Y'$), we may consider the Hilbert extension map associated to the regular embedding $\iota: X' \to X\times Y'$ by Proposition \ref{emb}. Then we have the intersection product of $\iota^{[[n]]}$ and the subscheme of $(X\times Y')^{[[n]]}$ isomorphic to $\Xi(\underline{\eta})\times Y'$, and we may define $$m_{\underline{\eta}}(f'):={\rm pr}_{1*}(X'^{[[n]]}\cdot (\Xi(\underline{\eta})\times Y')) \in \CH_{l' -\ell}(X'), $$ $$n_{\underline{\eta}}(f'):=\bar{f'}_*(X'^{[[n]]}\cdot (\Xi(\underline{\eta})\times Y')) \in \CH_{l'-\ell}(Y'). $$ For these classes, the base-change property (Proposition \ref{pullback}) follows from entirely the same proof as above.  \end{rem}

If $f$ satisfies an appropriate genericity condition, 
the subschemes $M_{\underline{\eta}}(f)$ and $N_{\underline{\eta}}(f)$ have the expected dimension 
and represent the multi-singularity loci classes up to certain multiplicities. 
To be precise, we present the following general statement. 

\begin{prop}\label{multi-sing}
If a proper morphism $f: X \to Y$ is  locally stable along the $\underline{\eta}$-locus, 
it holds that 
$$m_{\underline{\eta}}(f)= \#\Aut(\underline{\eta}/\eta_1) [M_{\underline{\eta}}(f)], \quad 
n_{\underline{\eta}}(f)= \#\Aut(\underline{\eta}) [N_{\underline{\eta}}(f)], $$
where $\underline{\eta}/\eta_1:=(\eta_2, \cdots, \eta_r)$ 
and $\#\Aut(\underline{\eta})$ is the order of the subgroup $\Aut(\underline{\eta})$ of $\mathfrak{S}_r$ which consists of permutations preserving the collection $\underline{\eta}$. 
\end{prop}

In this case, the map $f$ sends $M_{\underline{\eta}}(f)$ onto $N_{\underline{\eta}}(f)$ in finite-to-one; 
the degree is given by 
$$\deg_1 \underline{\eta}:=\frac{\#\Aut(\underline{\eta})}{\#\Aut(\underline{\eta}/\eta_1)}=\# \, \mbox{entries of type $\eta_1$ in $\underline{\eta}$}.$$
For instance, for the type $\underline{\eta}=(A_0, A_0, A_0, A_1, A_1)$, i.e., $\eta_1=\eta_2=\eta_3=A_0$ and $\eta_4=\eta_5=A_1$,  $\#\Aut(\underline{\eta})=3!2!$, $\#\Aut(\underline{\eta}/\eta_1)=2!2!$ and $\deg_1=3$. 

\

\proof 
Let $f: X\to Y$ be proper and locally stable. As in subsection \ref{invariance}, 
let $W$ be the subset of $C_{\Xi^\circ}$ consisting of $(x, z_I)$ such that $z_I$ is uniquely determined such that the germ of $f$ at $p_i$ is of type $\eta_i$. 
Then, ${\rm pr}_1: W \to M_{\underline{\eta}}^\circ$ and $\bar{f}: W \to N_{\underline{\eta}}^\circ$ are finite-to-one. In fact, at any $p_1\in M_{\underline{\eta}}^\circ$,  every point of ${\rm pr}_1^{-1}(p_1)$ in $W$ is of the form $(x', z_I)\in \Xi^\circ$, 
where $x'$ is obtained from $x$ by permuting $p_2, \cdots, p_r$ such that types $\eta_2, \cdots, \eta_r$ are preserved. 
Thus we have 
$$\deg \left[{\rm pr}_1: W \to M_{\underline{\eta}}^\circ\right] = \#\Aut(\eta_2, \cdots, \eta_r).$$
Also it is immediate to see 
$$\deg \left[\bar{f}: W \to N_{\underline{\eta}}^\circ\right] = \#\Aut(\eta_1, \cdots, \eta_r).$$

Recall that $C_{\underline{\eta}}:=C_{\underline{\eta}}(f)$ is obtained by intersecting $f^{[[n]]}: X^{[n]]} \to V_X$ with $\Xi(\underline{\eta})\times Y \subset V_X$ and it has dimension $l-\ell$. 
It is not difficult to see that 
$\CH_{l-\ell}(C_{\underline{\eta}}) \simeq \CH_{l-\ell}(C_{\Xi^\circ})$, and hence,   
$$[C_{\underline{\eta}}]=X^{[[n]]}\cdot (\Xi(\underline{\eta})\times Y) \in \CH_{l-\ell}(C_{\underline{\eta}}),$$ 
for it is so on $C_{\Xi^\circ}$.  
Then it follows that 
$${\rm pr}_{1*}[C_{\underline{\eta}}]=\# \Aut(\eta_2, \cdots, \eta_r)\, [M_{\underline{\eta}}(f)]$$ 
and also that 
$$\bar{f}_{*}[C_{\underline{\eta}}]=\# \Aut(\eta_1, \cdots, \eta_r)\, [N_{\underline{\eta}}(f)].$$ 
This completes the proof. 
\qed

\begin{rem}\upshape
When one takes a locally stable map $f:X \to Y$ and a base change $\varphi: Y' \to Y$ being transverse to $f$,  the obtained map $f': X' \to Y'$ is also locally stable 
and the $\underline{\eta}$-singularity loci of $f$ are pulled back to those of $f'$. 
Proposition \ref{multi-sing} implies a special case of the base-change property of the $\underline{\eta}$-singularity loci classes (Proposition \ref{pullback}) when applied to locally stable maps. 
\end{rem}

\subsection{Main Theorems}\label{thms}
Our main result is an affirmative solution to Kazarian's conjecture in algebraic geometry over 
an algebraically closed field $k$ of characteristic zero. 
That is the existence theorems of universal polynomials expressing $m_{\underline{\eta}}(f)$ and $n_{\underline{\eta}}(f)$. 
We call them the {\em  source and target multi-singularity Thom polynomials of type $\underline{\eta}$}, respectively\footnote{We also use this terminology for the polynomials up to constant multiplication (precisely saying, multiplied by $\frac{1}{\#\Aut(\underline{\eta}/\eta_1)}$ and $\frac{1}{\#\Aut(\underline{\eta})}$) like as in Table 1 in Introduction (cf. Proposition \ref{multi-sing}). }. 
As seen below, those Thom polynomials are mutually {\em well structured} subject to the hierarchy of adjacency relations of multi-singularity types. 
Here we use Chow rings $\CH^*(-)_\Q$ with rational coefficients.

\begin{thm} \label{main_thm1}
{\rm \bf (Target multi-singularity Thom polynomials)} 
Let $\underline{\eta}=(\eta_1, \cdots, \eta_r)$ be a multi-singularity type 
of map-germs of relative codimension $\kappa$ with $\codim \underline{\eta}=\ell$. 
Then, one can uniquely assign an abstract Chern polynomial 
with homogeneous degree $\ell$ 
$$R_{\underline{\eta}}=\sum a_I c^I \in \Q[c_1, c_2, \cdots],$$ 
called the residual polynomial of type $\underline{\eta}$, 
such that $R_{\underline{\eta}}$ does not depend on the order of entries in $\underline{\eta}$ 
and recursively admits the following property: 
For an arbitrary proper morphism $f: X \to Y$ of codimension $\kappa$ in $\Sm$, 
the target $\underline{\eta}$-singularity loci class $n_{\underline{\eta}}(f)$ 
is expressed by 
$$n_{\underline{\eta}}(f) =\sum f_*(R_{J_1})\cdots f_*(R_{J_s})$$
in $\CH^\ell(Y)_\Q$, where the summand runs over all possible partitions of the set $\{1, \cdots, r\}$ 
into a disjoint union of non-empty unordered subsets, 
$\{1, \cdots, r\}= J_1\sqcup \cdots \sqcup J_s \; (s\ge 1)$, and 
for each $J=\{j_1, \cdots, j_k\}$,  
$R_J$ stands for the residual polynomial $R_{(\eta_{j_1}, \cdots, \eta_{j_k})}$ 
evaluated by $c_i=c_i(f)$, the quotient Chern classes associated to $f$. 
\end{thm}

For the source multi-singularity loci class,  due to some technical reason, we consider only the case of $\kappa \ge 0$ and appropriately good maps $f: X \to Y$.  Note that a proper locally stable map $f: X \to Y$ introduces a stratification of $X$ into finitely many strata $M_{\underline{\xi}}(f)^\circ$ indexed with multi-singularity types $\underline{\xi}$ (here, strata are not assumed to be connected). In the present paper, we say that $f: X \to Y$ is {\em admissible} if $f$ is a proper locally stable map satisfying the following properties\footnote{This unfavorable condition (2) is due to the lack of classifying space of multi-singularity types.  The condition should be satisfied by  `universal map' on the expected classifying spaces  (see \S \ref{classifying_stack}). }: 
\begin{itemize}
\item[(1)] $f$ has only local singularities of quasi-homogeneous type (see Remark \ref{qh}); 
\item[(2)] for every non-empty stratum $M_{\underline{\xi}}(f)^\circ$, the top Chern class of its normal bundle is a non-zero divisor in $\CH^*(M_{\underline{\xi}}(f)^\circ)_\Q$. 
\end{itemize}

\begin{thm} \label{main_thm2}
{\rm \bf (Source multi-singularity Thom polynomials)} 
Let $\underline{\eta}=(\eta_1, \cdots, \eta_r)$ be a multi-singularity type of map-germs of codimension $\kappa \ge 0$. Then, for any admissible map $f: X \to Y$ of codimension $\kappa$,  the source $\underline{\eta}$-singularity loci class $m_{\underline{\eta}}(f)$ is expressed by 
$$m_{\underline{\eta}}(f) = \sum R_{J_1}\cdot f^*f_*(R_{J_2}) \cdots f^*f_*(R_{J_s})$$
in $\CH^{\ell-\kappa}(X)_\Q$, where the summand runs over all possible partitions of the set $\{1, \cdots, r\}$ 
into a disjoint union of non-empty unordered subsets 
so that the subset containing the entry $1 \in \{1, \cdots, r\}$ is denoted by $J_1$. 
\end{thm}

In particular,  {\em non-curvilinear} multiple-point theory exists -- Theorem \ref{main_thm2} guarantees the existence of a unique universal expression for the source multiple-point class $m_{A_0^r}(f)$ ($r\ge 4$) for certainly nice maps $f$ without any corank condition. 

Note that when applying the projection formula to the source Thom polynomials (the right hand side of the above equality), one recovers the known relation in 
Lemma \ref{iso} (1) 
$$f_*m_{\underline{\eta}}(f)=n_{\underline{\eta}}(f).$$

It is highly expected that the source multi-singularity Thom polynomials exist for arbitrary maps in any case of codimension $\kappa \in \Z$.  This part of Kazarian's conjecture still remains open, see \S \ref{perspective} for some detail. 

\begin{rem}\upshape 
Theorems \ref{main_thm1} and \ref{main_thm2} can be stated also for the $K$-theory $K^0$, instead of $\CH^*$, where residual polynomials are polynomials with $\Q[\beta, \beta^{-1}]$-coefficients, see \S \ref{multi_cob}. In fact, much more generally, Theorem \ref{main_thm1} will be proven for multi-singularity loci classes valued in a universal cohomology theory called the {\em algebraic cobordism} $\Omega^*$ \cite{LevineMorel, LevinePandhari} (Theorem \ref{main_thm}). 
\end{rem}

\begin{rem}\label{lcim_2}\upshape 
As mentioned in Remark \ref{lcim}, multi-singularity loci classes $m_{\underline{\eta}}(f)$ and $n_{\underline{\eta}}(f)$ can be defined for a proper local complete intersection (l.c.i.) morphism $f=F \circ \iota: X \to M \to Y$ with a regular embedding $\iota$ (with the normal bundle $\nu_{M, X}$) and $F: M \to Y$ in $\Sm$. Since such $f$ has the virtual normal bundle $\nu_f:=f^*TY-\iota^*TM+\nu_{M, X} \in K_0(X)$, the Chern classes $c_i(f)$ and Landweber-Novikov classes are defined. 
Then the above Theorem \ref{main_thm1} can be generalized for such proper l.c.i. morphisms, that is, 
the target class $n_{\underline{\eta}}(f) \in \CH^*(Y)$ is expressed as the specialization of the target Thom polynomial of $f_t: X_t \to Y\, (t \in T)$ via dynamic intersection (\S \ref{dynamic_intersection}). 
\end{rem}

\subsection{Examples} \label{appl}
There are some computations of residual polynomials by the restriction method  (see Remark \ref{restriction}). 
Here we pick up a few samples. 

\subsubsection{\bf Triple-point formula, revisited ($\kappa>0$)} \label{ex_tri}
The double and triple-point formulas reviewed in \S \ref{multiple} are reproduced from our theorems with residual polynomials
$$R_{A_0}=1, \quad R_{A_0^2}=-c_\kappa, \quad 
R_{A_0^3}=2\left(c_\kappa^2+\sum_{i=0}^{\kappa-1}2^ic_{\kappa-i-1}c_{\kappa+i+1}\right)$$
which are computed by the restriction method \cite{Kaz03}. 
Also the quadruple-point formula is obtained using $R_{A_0^4}$ computed in \cite{MarangellRimanyi}. 
For the type of $A_0^2=(A_0, A_0)$, there are two possible partitions 
$$\{\{1\}, \{2\}\}, \quad \{1,2\},$$ 
and hence, the double-point formula is recovered: 
$$m_{A_0^2}=R_{A_0}\cdot f^*f_*R_{A_0} + R_{A_0^2} = f^*f_*(1)-c_\kappa.$$
For the type of $A_0^3=(A_0, A_0, A_0)$, there are five partitions 
$$\{\{1\}, \{2\}, \{3\}\}, \; \{\{1\}, \{2, 3\}\}, \; \{\{1, 2\}, \{3\}\}, \; \{\{1, 3\}, \{2\}\}, \; \{1, 2, 3\},$$ 
and hence 
\begin{eqnarray*}
m_{A_0^3}
&=&R_{A_0}\cdot (f^*f_*R_{A_0})^2 + R_{A_0}\cdot f^*f_*R_{A_0^2} + 2 R_{A_0^2}\cdot f^*f_*R_{A_0} + R_{A_0^3}\\
&=& f^*f_*(f^*f_*(1)+R_{A_0^2})+ 2 R_{A_0^2}\cdot f^*f_*(1) + R_{A_0^3}\\
&=& f^*f_*(m_{A_0^2})-2c_\kappa f^*f_*(1) + R_{A_0^3}. 
\end{eqnarray*}
This is just the triple-point formula mentioned before. 

In the target space, especially in case of $\kappa=1$, we have 
$$n_{A_0^2}= s_0^2-s_1, \quad n_{A_0^3}=s_0^3-3s_0s_1+2(s_2+s_{01}).$$

\subsubsection{\bf Triple-point formula for the discriminant  ($\kappa\le 0$)} \label{ex_tri2}
In the negative codimensional case ($\kappa <0$),  the source multi-singularity Thom polynomial for $m_{\underline{\eta}}(f)$ has been missing yet (Theorem \ref{main_thm2} does not cover this case). Even though, in most of applications to classical enumerative problems, it suffices to compute the target Thom polynomials for $n_{\underline{\eta}}(f)$ in Theorem \ref{main_thm1}. 
For example, consider the case of $\kappa=-1$, and 
let $f: X^{l+1}\to Y^l$ be a proper map in $\Sm$. 
The critical locus $C(f)=\overline{A_1(f)}$ has the expected dimension $l-1$ 
(i.e., of codimension $2$), and the Thom polynomial is 
$m_{A_1}=R_{A_1}=c_1^2-c_2$, thus 
$$n_{A_1}=s_2-s_{01}.$$ 
The double and triple-point formulas of the discriminant set $D(f)$ in $Y$ are computed in \cite{Kaz06}: 
\begin{eqnarray*}
n_{A_1^2}&=&f_*(R_{A_1^2})+(f_*(R_{A_1}))^2\\
&=&-7s_3+8s_{11}-s_{001}+(s_2-s_{01})^2,\\
n_{A_1^3}&=&f_*(R_{A_1^3})+3f_*(R_{A_1})(f_*(R_{A_1}))^2+(f_*(R_{A_1}))^3\\
&=&138s_4-158s_{21}+2s_{02}+20s_{101}-2s_{0001}\\
&&\textstyle +3(s_2-s_{01})(-7s_3+8s_{11}-s_{001})+(s_2-s_{01})^3. 
\end{eqnarray*}

Using the above $n_{A_1^3}$ we can immediately reproduce classically known formulas in the following examples (i) and (ii). In old literature, $D(f)$ in $Y=\Proj^3$ is often called a {\em reciprocal surface}. 

\vspace{5pt}

\t
(i) We count the number $T$ of triple-points of the dual variety $S^* \subset \Proj^{3*}=:Y$ for a surface $S \subset \Proj^3$. 
Apply the Thom polynomial $n_{A_1^3}$ to the projection $f: X \to Y$, 
where $X \subset S \times Y$ is the incident variety of dimension $4$, 
that recovers an old result of Salmon \cite{Salmon, SempleRoth}: 
$$T=\textstyle \frac{1}{6}d(d-2)(d^7-4d^6+7d^5-45d^4+114d^3-111d^2+548d-960).$$

\vspace{5pt}

\t
(ii) We find the number $N_3(d)$ of plane curves of degree $d$ with three nodes passing through some general points. 
Apply the same Thom polynomial $n_{A_1^3}$ to the projection $f: X \to Y$ of the universal family, where $Y$ is a generic subspace $\Proj^3$ in the complete linear system $H^0(\Proj^2; \Ost(d))$. 
That leads to an old formula of Roberts (cf. \cite{KP}): 
$$N_3(d)=\textstyle \frac{1}{2}(9d^6-54d^5+9d^4+423d^3-458d^2-829d+1050).$$

\vspace{5pt}

Theorem \ref{main_thm1} guarantees the existence of $r$-tuple-point formula for $n_{A_1^r}(f)$ of discriminants, which provides a possible common generalization of both (i) and (ii) above. 
As a remark,  to count quadruple-points of $V^*$ for a $3$-fold $V \subset \Proj^4$, we can apply the Thom polynomial for $n_{A_1^4}$ with $\kappa=-2$ (because the incident variety is of dimension $6$),  while to find $N_4(d)$, we use the Thom polynomial for $n_{A_1^4}$ with $\kappa=-1$. Of course, these two Thom polynomials are different, because of different codimension $\kappa$, but when considering the restriction of $f$ to the critical locus $C(f)$ (see \S \ref{icis}),  they carry the same information about the `quadruple-point formula' for Legendre maps $f_C: C(f) \to Y$. 
This leads us to a specific version of Thom polynomials, called {\em Legendrian Thom polynmomials}, which will be discussed briefly in \S \ref{critical} and \S \ref{leg}. 

\subsubsection{\bf Instanton counting} \label{instanton}
In relation with the above example (ii), a modern typical question raises as {\it instanton counting} in string theory -- 
given a projective surface $S$ with an ample line bundle $L$, it asks how many curves with $r$ nodes ($A_1$-singularities) do exist in a generic subspace $\bP^r$ of the linear system $|L|$. The number should be counted in a universal way in terms of Chern classes of $S$ and $L$, that is known as part of the  {\em G\"ottsche conjecture}, which is now a theorem proven by several different techniques \cite{LT14, Rennemo, Berczi, Berczi23}. A key issue common to these approaches is the computation of a certain tautological integral over the Hilbert scheme $S^{[[n]]}$. 

On the other hand, our universal polynomial approach to this problem is as follows. Consider the universal family of nodal curves 
$$f: X \to Y:=\bP^r \subset \Proj(H^0(S; L))\qquad  (X \subset S \times Y)$$ 
where $\Proj^r$ is a generic subspace. 
Each $r$-nodal curve corresponds to a target $r$ multiple-point of type $A_1$ of $f$ in $Y$, and thus the desired enumeration is given by the degree $\frac{1}{r!} \int_Y n_{A_1^r}(f)$ in our theory. This leads to another proof of the G\"ottsche conjecture and its generalization, see  \cite[\S 10.1]{Kaz03}. 

Actually, these two different approaches deal with the same thing. 
By our definition, the target multi-singularity loci class $n_{A_1^r}(f)$ is determined by the intersection product of the Hilbert extension map $f^{[[n]]}: X^{[[n]]} \to (X\times Y)^{[[n]]}$ with the geometric subset $\Xi(X; A_1^r) \times Y$. 
Like as in the proof of Proposition \ref{pullback},  we can see that the intersection product is interpreted as the intersection of the relative Hilbert scheme $X(f)^{[[n]]}$ with the subscheme $\Xi(S; A_1^r) \times Y$ in $S^{[[n]]} \times Y$. It then leads to essentially the same quantity as computed in the tautological integral approach. 
The same enumerative problem can be considered for not only $A_1^r$ 
but also other multi-singularity type $\underline{\eta}$ of plane curves or of divisors/complete intersections as a higher dimensional generalization (cf. \cite{Rennemo, Tzeng, Vainsencher, Vainsencher07}).

\section{Algebraic Cobordism} \label{alg_cob}
\subsection{Oriented cohomology theory}
As an algebro-geometric counterpart to Quillen's program on complex cobordism theory, Levine-Morel \cite{LevineMorel} have introduced the notion of an oriented cohomology theory on $\Sm$, or more generally  an oriented Borel-Moore homology theory on $\Sch$. 
Let $R^*$ be the category of graded commutative rings. 
An {\em oriented cohomology theory} on $\Sm$ is an additive (pull-back) functor $A^*: \Sm^{op} \to R^*$ which enjoys   
the following axioms: 
\begin{itemize}
\item[-]
{\it pushforward} is associated to every proper morphism\footnote{In the original definition in \cite{LevineMorel}, morphisms are assumed to be projective, but this assumption is weakened by Gonzalez-Karu \cite{GonzalezKaru}.} 
so that it satisfies expected  properties (covariant functoriality, projection formula, compatibility for transverse cartesian square); 
\item[-] 
for every vector bundle $E \to X$ of rank $n$, let $P:=\bP(E^*)$ the scheme of hyperplanes of $E$ and let $\Ost(1)$ be the tautological quotient line bundle over $P$; then $A^*(P)$ is a free $A^*(X)$-module with basis $1, \xi, \cdots, \xi^{n-1}$, 
where we put 
$$\xi =c_1(\Ost(1)):=s^*s_*(1)$$ 
with the zero-section $s:P \to \Ost(1)$; 
\item[-]  
for an $\bA^n_k$-bundle $p: V\to X$, the pullback 
$$p^*: A^*(X) \to A^*(V)$$
is an isomorphism. 
\end{itemize}

Note that $\CH^*$ and $K^0$ are typical examples. 
Given an oriented cohomology theory $A^*$, one may use Grothendieck's construction of Chern classes of a vector bundle $E \to X$ of rank $n$ over $X \in \Sm$:  for each $0 \le i \le n$ 
there exists a unique element 
$$c_i(E) \in A^i(X)$$ 
such that 
$$c_0(E)=1 \;\; \mbox{and} \;\; \sum_{i=0}^n (-1)^i c_i(E) \xi^{n-i}=0.$$ 
One can easily check all the axioms of Chern classes. 

\subsection{Algebraic cobordism}
A successful theory of {\it algebraic cobordism} $\Omega_*$ on $\Sch$ has been established by Levine-Morel \cite{LevineMorel} following Quillen's axiomatic approach. 
For irreducible $Y \in \Sm$, we set 
$$\Omega^*(Y)=\Omega_{\dim Y - *}(Y).$$
Here we always assume the ground field $k$ to be of characteristic zero. 

\begin{thm}\label{thmLM} 
{\rm (\cite[Thm.7.1.1, Thm.7.1.3]{LevineMorel})} \\
(1) The theory $\Omega_*(-)$ is the universal oriented Borel-Moore homology theory on $\Sch$.\\
(2) The theory $\Omega^*(-)$ is the universal oriented cohomology theory on $\Sm$.
\end{thm}

The universality (2) means that for any oriented cohomology theory $A^*$ on $\Sm$, 
there is a unique natural transformation $\Omega^* \to A^*$. 
In particular, it is known that 
$$\Omega^*\otimes_{\bL^*} \Z \simeq \CH^*, \quad \Omega^*\otimes_{\bL^*} \Z[\beta, \beta^{-1}] \simeq K^0[\beta, \beta^{-1}],$$ 
where $\bL$ denotes the Lazard ring, which is isomorphic to the absolute ring, $\Omega^*(k)$ for $\Spec(k)$ \cite{LevineMorel}. 

A remarkable approach to re-interpreting $\Omega_*$ geometrically was later found by Levine-Pandharipande \cite{LevinePandhari}. 
Briefly, let us review the geometric construction. 

Let $Y \in \Sch$. Two morphisms $f: X \to Y$ and $f': X' \to Y$ are said to be {\it isomorphic} if there exists an isomorphism $\varphi: X \to X'$ so that $f=f'\circ \varphi$. 
Let $\Mcal(Y)$ denote the set of isomorphism classes over $Y$ 
of proper morphisms $f: X \to Y$ with $X \in \Sm$ which is graded by the dimension of $X$. 
The set $\Mcal(Y)$ is a monoid where the addition of isomorphism classes of $f: X_1 \to Y$ and $g: X_2 \to Y$ is defined by the disjoint union 
$$f+g: X_1 \sqcup X_2 \to Y.$$ 
Let $\mathcal{M}_*^+(Y)$ be the graded group completion of $\Mcal_*(Y)$. 

Now, consider a family of morphisms to $Y$; let 
$$F: M \to Y \times \bP^1$$
be a proper morphism with $M \in \Sm$ of pure dimension. Put $pr_1$ and $pr_2$ the projections to the first and second factors of $Y \times \bP^1$, respectively. 
Assume that 
$$\pi:=pr_2\circ F: M \to \bP^1$$ 
is a {\it double-point degeneration over $0 \in \bP^1$}, 
that means that 
$$\pi^{-1}(0)=A \cup B$$ 
where $A$ and $B$ are smooth Cartier divisors intersecting transversely in $M$. 
By restricting the map $pr_1\circ F$, we obtain maps $A \to Y$ and $B \to Y$. 
Let 
$$D=A\cap B$$ 
and $N_{A/D}$ the normal bundles of $D$ in $A$. Set 
$$\pi_P: P:=\bP(N_{A/D}\oplus \Ost_D) \to D.$$
Here it does not matter to distinguish $A$ and $B$, 
because it holds that 
$$N_{A/D}\otimes N_{B/D} \simeq \Ost_D$$ 
and thus the resulting $\Proj^1$-bundles are isomorphic.  
Composing $\pi_P$ with $pr_1\circ F|_D$, we have a map $P \to Y$. 
Let $p \in \bP^1$ be an arbitrary regular value of $\pi$ and $X_1$ the regular fiber at $p$, 
and it yields a map $pr_1\circ F: X_1 \to Y$. 
Then an associated {\it double-point relation over $Y$} is defined by 
$$[A \to Y] + [B \to Y] - [P\to Y] - [X_1 \to Y] \in \mathcal{M}_*^+(Y).$$
Let $\mathcal{R}_*(Y)$ be the subgroup generated by {\em all} double-point relations over $Y$. 
As a special case, if $\pi$ is smooth over $0, \infty \in \bP^1$, we may view $\pi$ as a double-point degeneration with $\pi^{-1}(0)=A \cup \emptyset$ and $D=\emptyset$, then na\"ive cobordism is also a double-point relation: 
$$[\pi^{-1}(0) \to Y]-[\pi^{-1}(\infty) \to Y] \in \mathcal{R}_*(Y).$$

The main theorem of \cite{LevinePandhari} gives a geometric presentation of the algebraic cobordism theory of \cite{LevineMorel}:

\begin{thm}\label{thmLP} 
{\rm (\cite[Thm.0.1]{LevinePandhari})} 
There is a canonical isomorphism: 
$$\Omega_*(Y) \simeq \mathcal{M}_*^+(Y)/\mathcal{R}_*(Y).$$
\end{thm}

\subsection{Formal group law} \label{fgl}
Let $A^*$ be an oriented cohomology theory on $\Sm$. 
For a line bundle $L$ on $X$ in $\Sm$, we have defined $c_1(L):=s^*s_*(1) \in A^*(X)$. 
In case of $A^*=\CH^*$, the first Chern class is additive, i.e., 
$c_1(L\otimes L')=c_1(L)+c_1(L')$, but it is not so for other cases. 
Namely, there is a unique power series 
$$F_A(u, v) =u+v+\sum_{i,j \ge 1} a_{ij} u^iv^j \in A^*(k)[[u, v]]$$ 
with $a_{ij} \in A^{1-i-j}(k)$ such that 
$$c_1(L\otimes M)=F_A(c_1(L), c_1(M))$$
for any line bundles $L$ and $M$ on $X$. 
The series should satisfy $F_A(u, v)=F_A(v,u)$ and $F_A(u, F(v, w))=F_A(F(u, v), w)$ (thus $a_{ij}$ satisfy that $a_{ij}=a_{ji}$ and some additional relations caused by the latter property), and then 
$F$ is called a {\em commutative formal group law} (of rank one). 

In particular, for $A^*=\CH^*$ (additive) and $K^0[\beta, \beta^{-1}]$ (multiplicative), respectively, 
$$F_{ad}(u, v)=u+v, \qquad 
F_{ml}(u, v)=u+v-\beta uv.$$

The universal theory $A^*=\Omega^*$ corresponds to the universal formal group law $F_{\bL}$. 
A computation shows that $a_{11}=-[\Proj^1]$, $a_{12}=a_{21}=[\Proj^1]^2-[\Proj^2]$ and so on \cite[\S 2.5]{LevineMorel}. The double-point relation measures the difference of $F_{\bL}$ from the additive one $F_{ad}=u+v$;  for divisors $A$, $B$ and $X_1$ in the above setting, $X_1$ is linearly equivalent to $A + B$, and thus the double-point relation means 
$$[P \to D]=A+B-F_{\bL}(A, B)=-\sum a_{ij}A^iB^j.$$


\subsection{Cohomology operations} 
We refer to a recent work of A. Vishik \cite{Vishik} on cohomology operations in algebraic geometry; 
results which we will later make use of are algebro-geometric counterpart to well-known facts/theorems in complex cobordism theory. 

Let $A^*$ and $B^*$ be oriented cohomology theories on $\Sm$. Fix $\kappa, \ell \in \Z$. 
A {\it cohomology operation} $G: A^\kappa \to B^\ell$ means a natural transformation when $A^\kappa$ and $B^\ell$ are regarded as functors from $\Sm^{op}$ to pointed sets. That is, for $Y \in \Sm$, 
$$G_Y: A^\kappa(Y) \to B^\ell(Y)$$ 
is a map between two sets so that $G_Y(0)=0$ and for any morphism $f: Y' \to Y$ in $\Sm$, 
the following diagram commutes: 
$$
\xymatrix{A^\kappa(Y') \ar[r]^{f^*} \ar[d]_{G_{Y'}} & A^\kappa(Y) \ar[d]^{G_Y}\\
B^\ell(Y') \ar[r]_{f^*} & B^\ell(Y)
}
$$
An operation $G$ is {\it additive} if the maps $G_Y$ are homomorphisms of abelian groups. 
A main result in \cite{Vishik} is the following theorem: 
%
\begin{thm}\label{Vishik1} 
{\rm (\cite[Thm.1.1]{Vishik})} 
Fix $\kappa, \ell \in \Z$. 
Let $B^*$ be an oriented cohomology theory on $\Sm$. 
Any family of homomorphisms 
$$G: \Omega^\kappa((\bP^\infty)^{\times l}) \to B^\ell((\bP^\infty)^{\times l})$$
for $l \in \Z_{\ge 0}$ 
commuting with  the action of the permutation group $\frak{S}_l$ 
and the pullbacks for specific morphisms listed below 
extends to a unique additive operation $G: \Omega^\kappa \to B^\ell$ on $\Sm$: 
\begin{itemize}
\item[-]  partial diagonals 
$(\bP^\infty \stackrel{\Delta}{\to} \bP^\infty\times \bP^\infty) \times \cdots$; 
\item[-]  partial Segre embeddings 
 $(\bP^\infty \times \bP^\infty \stackrel{{\rm s}}{\to} \bP^\infty) \times \cdots$; 
\item[-]  $(\Spec(k) \hookrightarrow \bP^\infty) \times (\bP^\infty)^{\times r}$ for any $r$; 
\item[-]   partial projection 
$(\bP^\infty)^{\times l} \to  (\bP^\infty)^{\times k}$ for any $k <l$. 
\end{itemize}
\end{thm}
Based on this theorem, a complete description of universal operations $\Omega^\kappa \to \Omega^\ell$ is obtained as follows \cite[Thm.3.18 (Thm.6.1)]{Vishik}: 

\begin{thm}\label{Vishik2} 
{\rm (\cite{Vishik})} 
Let $\varphi: \Omega^\kappa \to \Omega^\ell$ an additive operation on $\Sm$. 
Then  $\varphi: \Omega^\kappa \to \Omega^\ell\otimes_\Z \Q$ is uniquely written as a $\bL\otimes_\Z \Q$-linear combination of Landweber-Novikov operations with multi-indices $I=(i_1, \cdots, i_k)$ 
$$s_I: \Omega^\kappa(Y) \to \Omega^\ell(Y)\otimes_\Z \Q, 
\quad [f: X \to Y] \mapsto f_*(c_I(f)).$$ 
\end{thm}

\begin{rem}\upshape \label{pullback_omega}
{\bf (Pullback in $\Omega^*$)} 
A difficult step in the construction of the algebraic cobordism $\Omega^*(-)$ is to extend the pullback operations for arbitrary morphisms in $\Sm$ other than smooth morphisms. That was done in Levine-Morel \cite{LevineMorel}. 
On the other hand, in geometric construction of  Levine-Pandharipande \cite{LevinePandhari}, the pullback was explicitly described not for all morphisms but for certainly nice ones by {\it moving by translation}. The class of such nice morphisms include  $X \to \prod \bP^{N_i}$  in $\Sm$ \cite[\S8.2-8.3]{LevinePandhari}. 
\end{rem}

\subsection{Multi-singularity loci cobordism classes}\label{multi_cob}
Let $\underline{\eta}=(\eta_1, \cdots, \eta_r)$ be a multi-singularity type. 
To define $m_{\underline{\eta}}(f) \in \CH^*(X)$ and $n_{\underline{\eta}}(f) \in \CH^*(Y)$ (Definition \ref{multi-singularity_class}), we have used the fundamental class of the geometric subset in 
the Chow group of $X^{[[n]]} \in \Sch$: 
$$[\Xi(X; \underline{\eta})] \in \CH_*(X^{[[n]]}).$$ 
Instead, fixing a scheme structure, we may take the class of the structure sheaf 
$$[\Ost_{\Xi(X; \underline{\eta})}]\cdot \beta^d \in G_0(X^{[[n]]})[\beta, \beta^{-1}]$$ 
with $d=\dim \Xi(X; \underline{\eta})$) in the Grothendieck group of coherent sheaves on $X^{[[n]]}$ with parameter $\beta$. 
It is then straightforward to define $K$-theoretic multi-singularity loci classes 
$$m^K_{\underline{\eta}}(f) \in K^0(X)[\beta, \beta^{-1}], \;\;\; 
n^K_{\underline{\eta}}(f) \in K^0(Y)[\beta, \beta^{-1}]$$ 
in entirely the same way as in $\CH^*$ 
(e.g., see \cite[\S 2.2]{KS} for a summary of Intersection Theory in $K$-theory). 

More generally, we are concerned with algebraic cobordism multi-singularity loci classes. 

\begin{definition}\upshape \label{good}
{\bf ($\Omega$-assignment)}
Suppose that to every smooth scheme $X$, one assigns a distinguished bordism class 
$$\mu_{\Xi}(X; \underline{\eta}) \in \Omega_*(X^{[[n]]}).$$ 
We call it an {\em $\Omega$-assignment $\mu$ for the type $\underline{\eta}$}, if the following properties are satisfied: 
\begin{enumerate}
\item 
the image of $\mu_{\Xi}(X; \underline{\eta})$ in $\CH_*(X^{[[n]]})$ via the natural transformation is the fundamental class $[\Xi(X; \underline{\eta})]$; 
\item 
For any closed embedding $\iota: X' \hookrightarrow X$ in $\Sm$, with $X'$ and $X$ connected,  it holds that 
$$\mu_{\Xi}(X'; \underline{\eta})=(\iota^{[[n]]})^!(\mu_{\Xi}(X; \underline{\eta})) \in \Omega_*(X'^{[[n]]})$$
where $\iota^{[[n]]}: X'^{[[n]]} \to X^{[[n]]}$ is the regular embedding obtained in Proposition \ref{emb} 
and $(\iota^{[[n]]})^!$ is the refined Gysin map; 
\item 
The assignment to the disjoint union of $X_1$ and $X_2$ is of the form
$$\mu_{\Xi}(X_1 \sqcup X_2; \underline{\eta}) = \sum_{I \sqcup J} \; \mu_{\Xi}(X_1; \underline{\eta}_I) \times \mu_{\Xi}(X_2; \underline{\eta}_J) $$ 
in certain components of the algebraic cobordism group of the disjoint union
$$\bigoplus_{I \sqcup J}\; \Omega_*(X_1^{[[I]]}\times X_2^{[[J]]}) \subset \Omega_*((X_1\sqcup X_2)^{[[n]]})$$ 
where the summand runs over all partition $I \sqcup J= \{1, \cdots, r\}$ ($I$ or $J$ can be empty), 
$[[I]]:=\sum_{i \in I} n(\eta_i)$, and $\times$ is the cross product in $\Omega_*$. 
\end{enumerate}
\end{definition} 

To see the {\em existence} of such an $\Omega$-assignment, we first need to give $\mu_{\Xi}(X; \underline{\eta}_J)$ (for any $J$) for every {\em connected} smooth scheme $X$ so that (1) and (2) are satisfied. That is achieved by Hironaka's resolution of singularities (see Remark \ref{const_assignment} below). Then, we take the equality in (3) as the definition of $\mu_{\Xi}(X_1 \sqcup X_2; \underline{\eta})$ for the disjoint union of connected $X_1$ and $X_2$, and iteratively use this to define $\mu_{\Xi}(X; \underline{\eta})$ for $X$ having multiple connected components. 

\begin{definition}\upshape \label{multi-singularity_class_Omega}
{\bf (Multi-singularity loci cobordism classes)}
Suppose that we are given an $\Omega$-assignment $\mu$ for type $\underline{\eta}$. 
For an arbitrary proper morphism $f: X \to Y$ of codimension $\kappa$ in $\Sm$, we define the {\em multi-singularity loci cobordism classes associated to $\mu$} by 
\begin{eqnarray*}
m^\Omega_{\underline{\eta}; \mu}(f)&:=&{\rm pr}_{1*}((f^{[[n]]})^!(\mu_{\Xi}(X; \underline{\eta})\times [id_Y]))\, \in \Omega^{\ell-\kappa}(X), \\
n^\Omega_{\underline{\eta}; \mu}(f)&:=&\bar{f}_*((f^{[[n]]})^!(\mu_{\Xi}(X; \underline{\eta})\times [id_Y]))\, \in \Omega^\ell(Y)
\end{eqnarray*}
where $(f^{[[n]]})^!$ is the refined Gysin map of the Hilbert extension map in $\Omega^*$ \cite[\S 6.6]{LevineMorel}. 
\end{definition}
 
For short, we will often write $m^\Omega_{\underline{\eta}}(f)$ and $n^\Omega_{\underline{\eta}}(f)$ by omitting the letter $\mu$ in the subscript if it is easily understood. 
Obviously, changing the coefficients as to be in $\CH^*$ via the natural transformation, we recover $m_{\underline{\eta}}(f)$ and $n_{\underline{\eta}}(f)$ defined in Definition \ref{multi-singularity_class} (also for $K$-theory version as well). 

Emphasized is that multi-singularity loci cobordism classes is not uniquely determined by the type $\underline{\eta}$, for there are many different choices of $\Omega$-assignments. 

\begin{rem}\label{const_assignment}
\upshape 
{\bf (Resolution of singularities)}
Such an $\Omega$-assignment can be constructed by equivariant resolution of singularities as briefly explained below. Let $X$ be a connected smooth scheme. 
Note that $\Xi^\circ(X;\underline{\eta})$ is a locally closed smooth subscheme of $X^{[[n]]}$ given as an orbit of the action of local isomorphisms of $X$ at ordered $r$ points. 
Thus singularities of the Zariski closure $\Xi(X;\underline{\eta})$ are located along 
some adjacent strata in $\Xi(X;\underline{\eta}) - \Xi^\circ(X;\underline{\eta})$, and they are determined by the dimension of $X$ and local data of more complicated multi-singularity types than $\underline{\eta}$. Namely, any desingularization of $\Xi(X; \underline{\eta})$ are essentially indebted on those local data only. For an open set $U \subset \bA^m$, 
choose an embedded resolution of the locus $\Xi(U; \underline{\eta}) \subset U^{[[n]]}$ so that it is equivariant with respect to the groupoid action of local isomorphisms at $r$ points. By glueing such finitely many patches for $X=\bigcup U_\alpha$, we obtain an assignment $\mu_{\Xi}(X; \underline{\eta}) \in \Omega_*(X^{[[n]]})$. When the dimension $m$ goes up, additional desingularizations are created along new adjacent strata. In particular, we may first construct a desingularizaton of $\Xi(\Proj^N;\underline{\eta})$ and go to higher $N$ step by step, that produces a projective limit system $\{\mu_{\Xi}(\Proj^N; \underline{\eta})\}_{N \in \N}$. From this system, we obtain $\mu_{\Xi}(X; \underline{\eta})$ via the Gysin homomorphism of the inclusion of any smooth $X \subset \Proj^N$, that yields an $\Omega$-assignment which we seek for (at least in case that $X$ is quasi-projective). 
\end{rem}

The $\Omega$-assignment itself is just a bit hard technical matter in theory.  Truly interesting is the case that we are given a nice resolution of the geometric subset of a specific singularity type. For instance, there are at least two interesting results on Thom polynomials valued in $\Omega^*$ as far as the author knows. 

\begin{exam}\upshape 
{\bf (Degeneracy loci cobordism class)} Hudson \cite{Hudson12} has established the {\em algebraic cobordism Thom-Porteous formula} by using the Bott-Samelson resolution of the closures of Schubert cells (also see \cite{HPM20}). Those satisfy the above properties, i.e., they define an $\Omega$-assignment. 
Note that since the theory $\Omega^*$ established in \cite{LevineMorel, LevinePandhari}, there have appeared many successful works by several authors on {\em Schubert Calculus} in $\Omega^*$ and also in other oriented cohomology theory, e.g., \cite{HV11, HPM20}.  
\end{exam}

\begin{exam}\upshape 
{\bf (Double-point loci cobordism class)} Based on \cite{McCrory, Ronga73}, Audin \cite{Audin} presented the {\em double-point formula in complex cobordism $MU$}, which should also be valid in $\Omega^*$. In this case, we consider $\underline{\eta}=A_0^2$ ($n=2$) and the geometric subset $\Xi(A_0^2) = X^{[[2]]} = {\rm Bl}(X \times X)$ which is smooth, thus we have an $\Omega$-assignment for the type $A_0^2$. 
The relative Hilbert scheme provides a desingularization of the double-point locus for a generic map (Remark \ref{tri}). 
\end{exam}

\

\section{Proof of Theorems}\label{proof_thm}
In this section we prove Theorems \ref{main_thm1} and \ref{main_thm2}. 
For the simplicity, we first discuss $m_{\underline{\eta}}(f)$ and $n_{\underline{\eta}}(f)$ defined in $\CH^*$. At the last step, we will switch to algebraic cobordism $\Omega^*$ in order to properly employ Vishik's theorem (Theorem \ref{Vishik2}). 
The $\Omega^*$-version will be stated in Theorem \ref{main_thm}. 
A summary of Intersection Theory in \S \ref{Intersection} is available in $\Omega_*$ as well (see \cite[Chap.6]{LevineMorel}), and therefore, most of our arguments work in the $\Omega^*$-setting without any change, unless specifically mentioned. 

\subsection{Double-point degeneration}\label{dpr}
To begin with, we need to precisely know the behavior of multi-singularity loci classes under double-point degeneration. 

Let $\underline{\eta}=(\eta_1, \cdots, \eta_r)$ be a multi-singularity type of maps with codimension $\kappa=l-m$ 
and $\ell:=\ell(\underline{\eta}) \le l$.  
Let 
$$F: M \to \tilde{Y}:=Y \times \bP^1$$ 
be a proper morphism in $\Sm$ with $\dim M=m+1$ and $\dim Y=l$ such that 
$\pi=pr_2\circ F: M \to  \bP^1$ is a double-point degeneration over $0 \in \bP^1$ 
and  $1 \in \Proj^1$ is a regular value. 
Now set 
$$X_0:=\pi^{-1}(0)=A\cup B, \quad X_1=\pi^{-1}(1)$$  
and $D:=A\cap B$, $P:=\bP(N_{A/D}\oplus \Ost_D)$ as same as before. 
We denote the restriction maps by 
$$g_A: A \to Y, \quad g_B: B \to Y, \quad 
h_{P}: P \to Y, \quad h_{X_1}: X_1 \to Y,$$ 
and also, by taking disjoint unions, 
$$g: A \sqcup B \to Y, \quad h: P\sqcup X_1 \to Y.$$  
Those maps are of codimension $\kappa$, thus target multi-singularity loci classes of type $\underline{\eta}$ are defined. 
The following proposition is essential in our argument. 
We will also treat the same claim in the $\Omega^*$-setting.

\begin{prop}\label{degeneration_free} 
Under the above setting of an arbitrary double-point degeneration, we have 
$$n_{\underline{\eta}}(g)=n_{\underline{\eta}}(h)   \in \CH^\ell(Y).$$  
\end{prop}

Before starting the proof, we check the claim in simplest examples, that gives a clue as to what the proof should clarify. 

\begin{exam}\label{ex_dpf}
\upshape
{\bf (Double-point loci class) }
Put $Y=\Proj^2$ and let $M$ be the closure of the family of affine plane curves $xy-t=0$ in $Y \times \Proj^1$.  
Denote the embedding by 
$$F: M \to  \tilde{Y}=Y\times \Proj^1.$$ 
Set $X_0=M \cap (Y\times \{0\})=A\cup B$, where $A$ and $B$ are transverse lines in $Y$, and $X_1=M \cap (Y\times \{1\})$ a smooth conic. 
Put $D=A\cap B$, the double-point, and $P=\Proj^1=\Proj(T_DY)$.  
Now, we apply the double-point formula in \S \ref{multiple} to maps of codimension $\kappa=1$, 
$g: A\sqcup B \to Y$ and $h: P \sqcup X_1 \to Y$. 
Obviously, 
$$n_{A_0^2}(g)=g_*m_{A_0^2}(g)=2 \in \CH^2(Y).$$
On the other hand, there is no double-point of $h_{X_1}$ and $h_P: P \to Y$ is highly degenerate. 
Even so, the double-point formula makes sense for arbitrary maps \cite[\S9.3]{Fulton}. 
Since ${h_P}_*(1)=0$ and $h_P^*TY$ is trivial, we see 
\begin{eqnarray*}
m_{A_0^2}(h_P)&=&h_P^*{h_P}_*(1)-c_1(h_P^*TY-TP)\\
&=&-c_1(-TP)=c_1(T\Proj^1).
\end{eqnarray*}
Thus $n_{A_0^2}(h)=h_*m_{A_0^2}(h)=2$. 
\end{exam}

\begin{exam}\label{ex_degeneracy}
\upshape
{\bf (Degeneracy loci class) }
Let $M$, $Y$, $X_0$, $X_1$ and the map $F$ be the same as in Example \ref{ex_dpf}. 
Composing (the first factor of) $F$ with a generic linear projection to a fixed line $Y'$ in $Y=\Proj^2$, we obtain a proper morphism 
$$F': M \to  \tilde{Y}'=Y' \times \Proj^1.$$
We consider $g: A\sqcup B \to Y'$ and $h: P \sqcup X_1 \to Y'$, which are maps of codimension $\kappa=0$. 
Apply the simplest Thom-Porteous formula $m_{A_1}=c_1$ to those maps. Note that it is available for arbitrary maps \cite[\S 14.4]{Fulton}. 
Since $g_A$ and $g_B$ are non-singular, we see that 
$$n_{A_1}(g)=g_*m_{A_1}(g)=0 \in \CH^1(Z).$$ 
On the other hand, for $h_P^*TY'$ is trivial, we see 
\begin{eqnarray*}
m_{A_1}(h_P)&=&c_1(h_P^*TY'-TP)\\
&=&c_1(-TP)=- c_1(T\Proj^1).
\end{eqnarray*}
Then $n_{A_1}(h_P)={h_P}_*m_{A_1}(h_P)=-2$, which is caused by an excess intersection of $h_P$. 
There are two simple branch points for the projection of the conic $X_1$ to the line $Y'$; 
in other words, the singularity locus for $F': M \to  \tilde{Y}'$ meets $X_1$ at two points. 
Therefore, we see $n_{A_1}(h_{X_1})=2$. 
Hence
$$n_{A_1}(h)=n_{A_1}(h_{X_1})+n_{A_1}(h_P)=2-2=0.$$ 
\end{exam}

\subsection{Refined intersection} \label{refined} 
We prove Proposition \ref{degeneration_free} through several steps. 
We use the same setting and notations as in \S \ref{dpr}. 

To begin with, we fix some terminologies: consider a commutative diagram with  
two vector bundles $E\to B$ and $F\to B'$: 
$$
\xymatrix{
E \ar[r]^{\bar{\varphi}} \ar[d] & F \ar[d] \\
B \ar[r]_{\varphi} & B'}
$$
We call $\bar{\varphi}: E \to F$ a {\em vector bundle morphism} over $\varphi$ 
if it is linearly isomorphic fiberwise 
(i.e., $E$ is identified with the induced bundle $\varphi^*F$ on $B$), while we call $\bar{\varphi}: E \to F$ a {\em vector bundle map} over $\varphi$
if it is linear fiberwise and the rank of the linear maps on fibers can vary,  i.e., it is a section of ${\rm Hom}_k(E, \varphi^*F)\to B$.

\subsubsection{\bf  Quadratic map} \label{quadratic_map}
First, we recall an elementary fact which is equivalent to that the degree $\int c_1(T\Proj^1)=2$: 

\begin{lem}\label{2to1}
There is a vector bundle morphism $\bar{\psi}: T\Proj^1 \to \Ost_{\Proj^1}(1)$ over the quadratic map $\psi: \Proj^1 \to \Proj^1$, $\psi([z_0; z_1])=[z_0^2; z_1^2]$. 
\end{lem}

Set $L:=N_{A/D}$;  its dual $L^*$ is isomorphic to $N_{B/D}$ and 
$L\oplus L^*$ is the normal bundle $N_{M/D}$ of $D$ in $M$. 
As seen before, 
a double-point degeneration yields the $\Proj^1$-bundle 
$$\pi_P: P:=\Proj(L\oplus \Ost_D) \to D.$$ 
Now we consider another $\Proj^1$-bundle 
$$\pi_Q: Q:=\Proj(L\oplus L^*) \to D,$$ 
which is the exceptional fiber of the blow-up of $M$ along $D$. 
We introduce two line bundles, 
$$T\pi_P:=\ker d\pi_P \to P, \quad \nu^*:=\Ost_Q(1) \to Q,$$
the relative (fiber) tangent bundle of $\pi_P$ and the dual to the tautological line bundle $\nu \, (:=O_{Q}(-1))$ of $Q$, respectively. 

\begin{lem}\label{vbmorphism}
There is a vector bundle morphism $\bar{\psi}: T\pi_P \to \nu^*$ over a map $\psi: P \to Q$  
which is a family of fiberwise quadratic maps commuting with projections to $D$. 
$$
\xymatrix{
T\pi_P \ar[rr]^{\bar{\psi}} \ar[d] && \nu^* \ar[d] \\
P \ar[rr]_{\psi} \ar[dr]_{\pi_P}&& Q\ar[dl]^{\pi_Q}\\
&D&
}
$$
\end{lem}

\proof 
We build $\psi$ and $\bar{\psi}$ as follows. 
Take a local trivialization of $L$ and let 
$(z_0, z_1)$ and $(u_0, u_1)$ be linear coordinates of fibers of $L \oplus \Ost_D$ and  $L \oplus L^*$
subject to the splitting, respectively. 
Transitions are 
$(z_0, z_1) \mapsto (\alpha z_0, z_1)$ and $(u_0, u_1) \mapsto (\alpha u_0, \alpha^{-1}u_1)$ 
with $\alpha \in GL(1)$. 
We set locally 
$[u_0;u_1]=\psi([z_0; z_1]):=[z_0^2; z_1^2]$ in homogeneous coordinates of fibers. 
This is compatible with the transition
$$\psi([\alpha z_0; z_1])=[\alpha^2z_0^2; z_1^2]=[\alpha z_0^2; \alpha^{-1} z_1^2],$$ 
thus $\psi: P \to Q$ is well-defined. 
Applying Lemma \ref{2to1} fiberwise, we get the vector bundle morphism $\bar{\psi}$. 
More precisely, putting $z=z_0/z_1$, $u=u_0/u_1$ ($z_1, u_1\not=0$),  take local sections $\rd_z$ of $T\pi_P$ and $s_0:=u_1$ of $\nu^*$, 
then set $\bar{\psi}(\rd_z):=s_0$.  
Put $z'=\alpha z$ and $u_1'=\alpha^{-1}u_1$, then we have 
$\rd_{z'}=\alpha^{-1} \rd_{z}$ and $\bar{\psi}(\rd_{z'})=\alpha^{-1}s_0$. 
That is also for 
$w=z_1/z_0$, $v=u_1/u_0$, $ s_1:=-u_0$ and $\bar{\psi}(\rd_w):=s_1$.  
Hence, $\bar{\psi}$ is compatible with the transition. 
Note that $\rd_w=-z^2\rd_z$ and 
$s_1=-us_0$ ($u=\psi(z)=z^2$). Thus $\bar{\psi}$ is well-defined. 
\qed 

\

There are two global sections corresponding to $z=0, \infty$ in the above construction; 
we denote them by 
$$D_0, \; D_\infty \subset P \to D.$$ 
The fiberwise quadratic map ${\psi}: P \to Q$ sends them to sections 
$$D'_0, \; D'_\infty \subset Q \to D,$$ 
indeed, ${\psi}$ is a double covering branched along $D'_0$ and $D'_\infty$.

\subsubsection{\bf Desingularization} \label{desing} 
Consider the projectivized cotangent bundle 
$$\tilde{\pi}: \Proj(T^*M) \to M$$ 
which is also the Grassmannian bundle of $m$-planes in $TM$ (note that $\dim M=m+1$). 
We denote a point of $\Proj(T^*M)$ 
by the pair $q=(x, \lambda_q)$ of a point $x=\tilde{\pi}(q) \in M$ and a hyperplane $\lambda_q \subset T_x M$. 
At every point $x \in M -D$, we take the tangent space at $x$ of the hypersurface $X_{\pi(x)}$ 
$$T_xX_{\pi(x)}=\ker d\pi_x \subset T_x M.$$
Collect those hyperplanes and also limiting hyperplanes: we denote the Zariski closure by  
$$\tilde{M}:={\rm Closure} \{ \; (x, T_xX_{\pi(x)}) \in \Proj(T^*M) \; | \; x \in M - D \; \}$$
with the natural projection 
$$\tilde{\pi}: \tilde{M} \to M.$$

\begin{lem}\label{blowup} 
$\tilde{\pi}: \tilde{M}-\tilde{\pi}^{-1}(D) \to M-D$ is isomorphic 
and $\tilde{M}$ is non-singular. 
In particular,  $\tilde{\pi}^{-1}(D)$ is identified with $Q$ defined previously, and its normal bundle in $\tilde{M}$ is isomorphic to $\nu=O_{Q}(-1)$.  
\end{lem}

\proof 
Take an affine chart $(U, (x, y, u_2, \cdots, u_m))$ of $M$ centered at a point of $D$ such that $X_0 \cap U$ is defined by $xy=0$, i.e.,  locally, $(u_2, \cdots, u_m)$ are coordinates of $D \cap U$ and $x,y$ are coordinates normal to $D$ in $A$ and $B$, respectively. The tangent space $TX_t$ at $q=(x_0, y_0, u) \in X_t \cap U$ with $t\not=0$ is given by $x_0dy+y_0dx=0$ in $T_q M$, i.e., it is the direct sum of $T_qD$ and the line spanned by $x_0 \rd_x-y_0\rd_y$. Then we can identify $\tilde{\pi}^{-1}(U)$ with the blow-up ${\rm Bl}_{D \cap U} U$ of $U$ along $D \cap U$ by $(q,[x_0; -y_0]) \mapsto  (q, [x_0; y_0])$, and the claim follows. 
\qed 

\

Let $\tilde{X}_0$ denote the schematic preimage $\tilde{\pi}^{-1}(X_0) \subset \tilde{M}$; it consists of irreducible components $\tilde{A}$ and $\tilde{B}$ 
(the strict transforms of $A$ and $B$, respectively) and $Q$: 
$$[\tilde{X}_0]=[\tilde{A}] + [\tilde{B}] + 2[Q].$$
Those components have transverse intersection: 
$$\tilde{A}\cap Q=D_0', \quad \tilde{B}\cap Q=D_\infty', \quad \tilde{A}\cap \tilde{B}=\emptyset.$$  

Put $\tilde{Y}:=Y\times T$ with $T:= \Proj^1$. We set 
$$\tilde{F}:=F \circ \tilde\pi: \tilde{M} \to \tilde{Y}$$
and also 
$F_1:= pr_1\circ \tilde{F}: \tilde{M} \to Y$ and $F_2:=pr_2\circ \tilde{F}: \tilde{M} \to T$. 
 
For any $t\not=0$, $X_t$ is isomorphic to $\tilde{X}_t$. 
Denote the restriction by $$f_t: \tilde{X}_t \to Y\, (=Y\times \{t\}).$$
Let $i_g: \tilde{A} \sqcup \tilde{B} \hookrightarrow \tilde{M}$,  $\iota: \tilde{X}_t \hookrightarrow \tilde{M}$ and $j: Q \to \tilde{M}$ be the inclusion and $i_h= (j\circ \psi) \sqcup \iota: P \sqcup \tilde{X}_t \to  \tilde{M}$ ($t\not=0$) 
where $\psi: P \to Q$ is the quadratic map in Lemma \ref{vbmorphism}.  Put 
$$g=F_1\circ i_g: \tilde{A} \sqcup \tilde{B} \to Y, \quad 
h:=F_1\circ i_h:  P \sqcup \tilde{X}_t \to Y.$$ 
Here $h$ depends on $t\, (\not=0)$, but this simplification may not cause any confusion.

The multi-singularity loci classes of type $\underline{\eta}$ for maps $g$ and $h$ ($t\not=0$)
are defined  (Definition \ref{multi-singularity_class}): 
$$m_{\underline{\eta}}(g) \in \CH^{\ell-\kappa}(\tilde{A})\oplus \CH^{\ell-\kappa}(\tilde{B}), $$
$$m_{\underline{\eta}}(h) \in \CH^{\ell-\kappa}(P)\oplus \CH^{\ell-\kappa}(\tilde{X}_t).$$
We compare these classes in $\CH_{l-\ell}(\tilde{M})=\CH^{\ell-\kappa+1}(\tilde{M})$ through 
the pushforward of $i_g$ and $i_h$. 
Below, as a notational convention, we may drop the letter $\iota_*$ for the pushforward of inclusions $\iota$ to the ambient space, if it is easily understood from the context. 

The latter class $m_{\underline{\eta}}(h)$ is the sum of parts supported on $P$ and $\tilde{X}_t$. 
Since $h_P: P \to Y$ is degenerate along the fibers of $P\to D$, 
the multi-singularity loci class supported on $P$ contains an {\em excess intersection}.  We denote the class by 
$$m_{\underline{\eta}}(h)^P \in \CH_{l-\ell}(P).$$
On the other hand, for any regular value $t$, 
we denote by 
$$m_{\underline{\eta}}(h)^{X_t} \in \CH_{l-\ell}(\tilde{X}_t)$$
the multi-singularity loci class supported on $\tilde{X}_t$. 
The image in $\CH_{l-\ell}(\tilde{M})$ does not depend on $t \not=0$. 
Obviously, we have  
\begin{lem} \label{lem1}
$\displaystyle 
(i_h)_*m_{\underline{\eta}}(h)=
\iota_*m_{\underline{\eta}}(h)^{X_t}+j_*\psi_*m_{\underline{\eta}}(h)^P 
\in \CH_{l-\ell}(\tilde{M})$. 
\end{lem}

\

\subsubsection{\bf Dynamic intersection} \label{dynamic}
Put $T^*=T-\{0\} (=\Proj^1 -\{0\})$ and 
$$\tilde{M}^*:=\tilde{M}-\tilde{X}_0=\tilde{M}\times_T T^*, \quad
\tilde{Y}^*:= Y \times T^*.$$
Also we set 
$$\mathcal{M}_P^\circ:=(P \times_k T^*) \sqcup \tilde{M}^* \; \mbox{ (disjoint)}$$
as a smooth $T^*$-scheme, and define a proper morphism over $T^*$
$$\tilde{F}^\circ:=(h_P\times id_{T^*})\sqcup \tilde{F}:  \mathcal{M}_P^\circ \to \tilde{Y}^*$$
to be the family 
$$
h_t: P \sqcup \tilde{X}_t \to Y \qquad (t \in T^*).
$$ 
Then we have the multi-singularty loci class for the map $\tilde{F}^\circ$ of type $\underline{\eta}$:
$$m_{\underline{\eta}}(\tilde{F}^\circ)  = \alpha + \beta \in \CH_{l-\ell+1}(\mathcal{M}_P^\circ)$$ 
where we denote 
$$\alpha:=m_{\underline{\eta}}(\tilde{F}^\circ)^{\tilde{M}^*} \in \CH_{l-\ell+1}(\tilde{M}^*)$$
and  
$$\beta:=m_{\underline{\eta}}(\tilde{F}^\circ)^{P \times_k T^*} \in \CH_{l-\ell+1}(P \times_k T^*).$$

We take the intersection product 
$$\alpha_t:=\tilde{X}_t\cdot \alpha  = \iota^*(\alpha) \in \CH_{l-\ell}(\tilde{X}_t)$$ 
($t\not=0$) where $\iota: \tilde{X}_t \hookrightarrow \tilde{M}^*$ is the inclusion. 
Also we set $\beta_t:=(P\times_k \{t\})\cdot \beta \in \CH_{l-\ell}(P)$.  

\begin{lem}\label{lim1}
It holds that 
$$\alpha_t = m_{\underline{\eta}}(h)^{X_t}  \in \CH_{l-\ell}(\tilde{X}_t), 
\quad 
\beta_t = m_{\underline{\eta}}(h)^{P}  \in \CH_{l-\ell}(P).$$ 
\end{lem}

\proof 
The inclusion $Y = Y \times \{t\} \hookrightarrow \tilde{Y}^*$ ($t\not=0$)  is transverse to the map $\tilde{F}^\circ: \mathcal{M}_P^\circ \to \tilde{Y}^*$,  
and the fiber square gives $h: P \sqcup \tilde{X}_t \to Y$. 
The lemma follows from Proposition \ref{pullback}. \qed 

\

Recall the basic construction of dynamic intersection in \S \ref{refined_intersection}. 
Write $\alpha=\sum n_i [\mathcal{Y}_i]$ where $\mathcal{Y}_i$ are irreducible subvarieties of $\tilde{M}^*$ of dimension $l-\ell+1$, and take the Zariski closure $\overline{\mathcal{Y}_i} \subset \tilde{M}$. 
Let 
$$\tilde{\alpha}:=\sum n_i [\overline{\mathcal{Y}_i}] \in \CH_{l-\ell+1}(\tilde{M}).$$
Then the limit intersection class is defined,  and it actually satisfies 
$$\lim_{t \to 0} \alpha_t =\tilde{X}_0 \cdot \tilde{\alpha} 
=(\tilde{A}+\tilde{B}+2Q)\cdot \tilde{\alpha} \in \CH_{l-\ell}(\tilde{X}_0).$$
Note that $\tilde{X}_0 \cdot \tilde{\alpha}$ and $\alpha_t=\tilde{X}_t \cdot \tilde{\alpha}$ ($t\not=0$) coincide in $\tilde{M}$, for $[\tilde{X}_0]=[\tilde{X}_t] \in \CH_{l-\ell}(\tilde{M})$. 

Also for $\beta=\sum m_k [\mathcal{Y}'_k]$ where $\mathcal{Y}'_k$ are irreducible subvarieties of $P$, 
we put 
$$\tilde{\beta}:=\sum m_k [\overline{\mathcal{Y}'_k}] \in \CH_{l-\ell+1}(P \times_k T).$$
The limit class is defined and $\lim_{t \to 0} \beta_t=\beta_t \in \CH_{l-\ell}(P)$ ($t\not=0$).

Comparing the limit of $\alpha_t$ and that of $\beta_t$, we obtain the following main lemma, which will be shown later in \S \ref{proof_lem3}. The lemma introduces a certain {\em refined intersection class} $\delta_D$, which is the contribution to multi-singularity loci classes of $h$ and $g$ from bi-singularities lying on the same fiber of $D_0' \sqcup D_\infty' \to D$ (cf. Example \ref{ex_dpf}). 

\begin{lem} \label{lem3}
There is a class $\delta_D \in \CH_{l-\ell}(D'_0 \sqcup D'_\infty)$ satisfying that 
\begin{align*}
&j_*\psi_*m_{\underline{\eta}}(h)^P = j_*\delta_D - 2\, Q \cdot \tilde{\alpha}  \tag{1}\\
&(i_g)_*m_{\underline{\eta}}(g) = j_*\delta_D + (\tilde{A}+\tilde{B})\cdot \tilde{\alpha} \tag{2}
\end{align*}
in $\CH_{l-\ell}(\tilde{M})$ where $j: Q \to \tilde{M}$ is the inclusion. 
\end{lem}

\t
{\it Proof of Proposition \ref{degeneration_free}}: 
It follows from Lemmas \ref{lim1} and \ref{lem3} that 
$$(i_g)_*m_{\underline{\eta}}(g)-j_*\psi_*m_{\underline{\eta}}(h)^P
=(\tilde{A}+\tilde{B}+2Q)\cdot \tilde{\alpha}
=\tilde{X}_0\cdot \tilde{\alpha}=\tilde{X}_t\cdot \tilde{\alpha}=\iota_*m_{\underline{\eta}}(h)^{X_t}$$
in $\CH_{l-\ell}(\tilde{M})$. 
Hence, by Lemmas \ref{lem1}, we see 
$$(i_g)_*m_{\underline{\eta}}(g)=(i_h)_*m_{\underline{\eta}}(h).$$
Since $g=F_1 \circ i_g$ and $h=F_1\circ i_h$, we have 
\begin{eqnarray*}
n_{\underline{\eta}}(g)&=&g_*m_{\underline{\eta}}(g)\\
&=&(F_1)_*  (i_g)_*m_{\underline{\eta}}(g)\\
&=&(F_1)_*  (i_h)_*m_{\underline{\eta}}(h)\\
&=&h_*m_{\underline{\eta}}(h)\\
&=&n_{\underline{\eta}}(h).
\end{eqnarray*}
This completes the proof of Proposition \ref{degeneration_free}. \qed

\

\begin{rem}\label{alpha}
\upshape 
We have defined $\tilde{\alpha} \in \CH_{l-\ell+1}(\tilde{M})$ by dynamic intersection.  
Note that $\tilde{\alpha}$ is entirely different from the multi-singularity loci class $m_{\underline{\eta}}(\tilde{F}) \in \CH_{l-\ell+1}(\tilde{M})$ for the map $\tilde{F}: \tilde{M} \to \tilde{Y}=Y \times T$. 
Indeed, the divisor $Q$ always becomes an irreducible component of the singular locus of $\tilde{F}$ and thus it contributes to some part of $m_{\underline{\eta}}(\tilde{F})$, while $\tilde{\alpha}$ does not have such a component. 
\end{rem}

\subsubsection{\bf Jet section} \label{jets} 
Let $q \in \tilde{M}$ and $k$ sufficiently large. 
Regarding $\lambda_q$ as the germ of a hyperplane in a local chart of $M$ centered at $\tilde{\pi}(q)$, we denote by $\Psi_q$ the $k$-jet of the composed map-germ 
$$\xymatrix{
(\lambda_q, 0) \ar[r]^{incl\;\;\;\;} & (M, \tilde{\pi}(q))  \ar[r]^{F} & (Y, f_1(q)).}
$$
In other words, it is an element of the $k$-jet space 
$$\Psi_q \in J^k(\lambda_q, T_{F_1(q)}Y)=\bigoplus_{i=1}^k \Hom({\rm Sym}^i(\lambda_q), T_{F_1(q)}Y). $$
Let $H \to \tilde{M}$ be the tautological vector bundle of rank $m$, i.e., the bundle whose fiber at $q \in \tilde{M}$ is $\lambda_q$, 
then we have a jet section 
$$\Psi: \tilde{M} \to J^k(H, F_1^*TY).$$
Note that the restriction of $\Psi$ to $\tilde{X}_t$ is nothing but the jet section of $f_t: X_t \to Y$ at $\tilde{\pi}(q) \in X_t$ ($t\not=0$), while the restriction to $Q$ is entirely different from the jet section of $F_1\circ \pi_Q: Q \to Y$. 
Indeed, it follows from Lemma \ref{blowup} that 
$$H|_Q \simeq \nu\oplus \pi_Q^*TD.$$

We compare the section $\Psi|_Q$ over $Q$ with the $k$-jet extension of $h_P=F_1 \circ \psi: P \to Y$, which is a section of $J^k(TP, h_P^*TY)$ over $P$. 
Since the tangent bundle $TP$ splits as 
$$TP=T\pi_P\oplus \pi_P^*TD,$$ 
we associate a bundle morphism 
$$\xymatrix{
\bar{\psi}\oplus id_{TD}: TP \ar[r] & \nu^* \oplus \pi_Q^*TD}
$$ 
over $\psi: P \to Q$ using $\bar{\psi}$ in Lemma \ref{vbmorphism}. 
This induces an isomorphism between the tangent space $T_pP$ and the vector space $\lambda_q$ at every $p \in P$ and $q=\psi(p)$, if we identify $\nu_q$ and $\nu_q^*$ (fibers at $q$).  Combining it with the jet section $\Psi_p$ we have the $k$-jet of a map-germ $(P, p) \to (Y, h_P(p))$
which we may regard as a deformation of the jet of $h_P$ at $p$. 
However, we can not do this globally, for $\nu \not\simeq \nu^*$ in general. 
So we now suppose to take a section $s_\nu$ of $\Hom(\nu^*, \nu)$ over $Q$ -- 
indeed, $s_\nu$ may be the zero section. 
We then think of the bundle map
$$(s_\nu \circ \bar{\psi})\oplus id_{TD}:TP \to  \nu \oplus \pi_Q^*TD=H$$
over $\psi: P \to Q$. Combining this map with the section $\Psi: \tilde{M} \to J^k(H, F_1^*TY)$, we obtian a section 
$$\xymatrix{P \ar[r]& J^k(TP, h_P^*TY).}$$ 
In other words, 
the obtained jet section assigns to each $p \in P$ the $k$-jet of the germ of composed map
$$\xymatrix{
(P, p) \ar[r] & (\lambda_{q}, 0)  \ar[r]^{\Psi\; \;\;} & (Y, h_P(p)).}$$ 
It actually coincides with the $k$-jet of $h_P$ at $p$ if $s_\nu=0$.

For a mono-singularity type $\eta_1$, the singularity loci class $m_{\eta_1}(\Psi) \in \CH^*(\tilde{M})$ is defined. 
Then the intersection product $Q\cdot m_{\eta_1}(\Psi)$ is equal to the loci class $m_{\eta_1}(\Psi|_Q) \in \CH^*(Q)$. 
Also we obtain the loci class $m_{\eta_1}(h_P) \in \CH^*(P)$ of $h_P: P \to Y$, which actually contains excess intersection.  
The following lemma formulates a general property suggested by Example \ref{ex_degeneracy}.

\begin{lem}\label{change_sign}
$(j\circ \psi)_*\beta_t= -2\, Q \cdot \tilde{\alpha} \in \CH^*(\tilde{M})$. 
\end{lem}

\proof 
Let $L_Q$ be the line bundle over $\tilde{M}$ associated to the divisor $Q$, then $c_1(L_Q)=[Q]$ and $\nu=j^*L_Q$. 
The map $\bar{\psi}$ between the total spaces of line bundles $T\pi_P$ and $\nu^* (\subset L_Q^*)$ sends a section to a section (with degree two), and its Gysin-image in $\CH^*(\nu^*)=\CH^*(\tilde{M})$ equals $2c_1(L_Q^*)=-2[Q]$. 
The class $m_{\eta_1}(h_P)$ is obtained by the pullback of $m_{\eta_1}(\Psi|_Q)$ via ${\psi}$ together interchanging $\nu$ and $\nu^*$. Hence, 
by the projection formula, $(j\circ \psi)_*m_{\eta_1}(h_P)=j_*\psi_*(-\psi^* m_{\eta_1}(\Psi))=- j_* (\psi_*(1)\cdot m_{\eta_1}(\Psi))=-2\,Q\cdot m_{\eta_1}(\Psi)$. 
\qed

\begin{rem}\label{change_sign_2}\upshape 
Lemma \ref{change_sign} can also be viewed in another way. 
The map $h_P: P \to Y$ is the composition of $\psi: P \to \tilde{M}$, $\tilde{F}: \tilde{M} \to Y\times T$ and $Y\times T \to Y$. As mentioned in Remark \ref{alpha}, the mono-singularity loci class $m_{\eta_1}(\tilde{F})$ contains some part supported on $Q$.

\end{rem}

\subsubsection{\bf  Refined intersection} \label{refined} 
Recall that $m_{\underline{\eta}}(h)$ is the image via ${\rm pr}_{1*}$ of the intersection product of the Hilbert extension map 
$$h^{[[n]]}: (P \sqcup \tilde{X}_t)^{[[n]]} \to ((P \sqcup \tilde{X}_t)\times Y)^{[[n]]}$$
with the geometric subset $V:=\Xi(P \sqcup \tilde{X}_t; \underline{\eta}) \times Y$ (Definition \ref{multi-singularity_class}). 
The class $m_{\underline{\eta}}(h)^P$ is the part supported on $P$, that is, 
it comes from the intersection product on connected components of $(P \sqcup \tilde{X}_t)^{[[n]]}$ whose first factor (via ${\rm pr}_1$) is mapped to $P$, i.e., the union of 
$$P^{[[J]]} \times \tilde{X}_t^{[[K]]} \quad \mbox{($J\sqcup K=\{1, \cdots, r\}$ with $1 \in J$)}.
\eqno{(c)}$$

A linear vector field of fiber tangents $\bv \in \Gamma(T\pi_P)$ over $P$ is defined by 
$$\bv(z)=z\rd_z\;\;\;  (\mbox{or}\;\;\; \bv(w)=-w\rd_w)$$
using the fiber coordinate $z$ or $w$ appearing in local trivialization of $\pi_P: P \to D$ (indeed, $\bv$ is well-defined on the entire space $P$). 
We also denote by the same letter $\bv$ as a section of $TP$ obtained via the natural inclusion $T\pi_P \subset TP$. 
Notice that the zero locus of $\bv$ coincides with $D_0 \sqcup D_\infty \subset P$. 
Furthermore, $\bv$ is extended to a section on $P \sqcup \tilde{X}_t$ so that it assigns the zero on $\tilde{X}_t$. 
Then we introduce a section on each component 
$$s_{\bbv}: P^{[[J]]} \times \tilde{X}_t^{[[K]]} \to N_{h_P}$$
by assigning to every $(x^{|J|}, z_I) \in  P^{[[J]]}$ the derivation (see \S \ref{DH}) 
$$s_{\bbv}(x^{|J|}, z_I):=Dh_P^{[[|J|]]}\left(\bigoplus_{p \in {\rm Supp}(I)} \bv(p)\right) \in (N_{h_P})_{(x,z_I)}$$
and the zero to other $(x^{|K|}, z_I) \in \tilde{X}_t^{[[K]]}$ 
and by taking the direct sum. 

The normal cone $C_V$ associated to $V$ lives in the normal bundle of the Hilbert extension map $h^{[[n]]}$ which is isomorphic to $(h^*TY)^{[[n]]}$. We denote the normal bundle over the union of components $(c)$ by $N_{h_P}$. Then by the definition,
$$m_{\underline{\eta}}(h)^P={\rm pr}_{1*}s_{\bbv}^*[C_V] \quad  \in \CH_{l-\ell}(P).$$ 

Note that $D_0 \sqcup D_\infty$ is always contained in the preimage of $s_{\bbv}$ of the normal cone $C_V$. 
We distinguish irreducible components 
of $s_{\bbv}^{-1}(C_V)$ as follows. 

Let $Z \subset (P \sqcup \tilde{X}_t)^{[[n]]}$ be the closed subset consisting of $(x, z_I)$ such that the support of $I$ contains two distinct points $p, p' \in D_0 \sqcup D_\infty$ with $\pi_P(p)=\pi_P(p') \in D$ and the first factor ${\rm pr}_1(x) \in P$ is $p$.  
The components of $s_{\bbv}^{-1}(C_V)$ which are contained in $Z$ define the refined intersection class 
$$(s_{\bbv}^*[C_V])^Z \in \CH_{l-\ell}(Z).$$ 
Since ${\rm pr}_1(Z)$ is contained in $D_0 \sqcup D_\infty$, we put 
$$\delta_D:={\rm pr}_{1*}((s_{\bbv}^*[C_V])^Z) \in \CH_{l-\ell}(D_0 \sqcup D_\infty).$$
By the natural identification via $\psi$, we may also think of $\delta_D$ as a class in $\CH_{l-\ell}(D'_0 \sqcup D'_\infty)$. 

\

\subsubsection{\bf Proof of Lemma \ref{lem3}} \label{proof_lem3}  
First we show the equality
\begin{align*}
&j_*\psi_*m_{\underline{\eta}}(h)^P = j_*\delta_D - 2\, Q \cdot \tilde{\alpha}.  \tag{1}
\end{align*}
By the definition, the class $\delta_D$ contributes to $m_{\underline{\eta}}(h)^{P} \in \CH_{l-\ell}(P)$, and the sum of components of $s_{\bbv}^{-1}(C_V)$ which are not contained in $Z$ corresponds to, via ${\rm pr}_{1*}$,  the residual part of $m_{\underline{\eta}}(h)^{P}$ off $\delta_D$. Push the residual class to $Q$ via $\psi$. 
Through the interpretation of $h_P$ by the jet section $\Psi$ and a similar argument as in the proof of Lemma \ref{change_sign}, we can see that the obtained class in $Q$ is just the intersection class of $\tilde{\alpha}$ and $\psi_*[P]=2[Q]$ with the opposite sign:  
$-2\, Q \cdot \tilde{\alpha} \in \CH_{l-\ell}(Q)$. 
This implies (1). 

Next, we show
\begin{align*}
&(i_g)_*m_{\underline{\eta}}(g) = j_*\delta_D + (\tilde{A}+\tilde{B})\cdot \tilde{\alpha}. \tag{2}
\end{align*}
We remark that 
\begin{itemize}
\item[(a)] at every point $p \in D_0$ (resp. $p' \in D_\infty$), 
the tangent space $T_pP$ is naturally identified with $\lambda_q=T_q \tilde{A}$  (resp. $\lambda_{q'}=T_{q'} \tilde{B}$) where $q=\psi(p)$ and $q'=\psi(p')$; 
\item[(b)] thus, punctual Hilbert schemes are isomorphic, $\Hilb^n(\Ost_{P, p}) \simeq \Hilb^n(\Ost_{\tilde{A}, q})$ and $\Hilb^n(\Ost_{P, p'}) \simeq \Hilb^n(\Ost_{\tilde{B}, q'})$; 
\item[(c)] also the jets of $h_P: P \to Y$ at $p$ and $p'$ are identified with 
the jets of $g: \tilde{A} \sqcup \tilde{B} \to Y$ at $q$ and $q'$, respectively. 
\end{itemize}
Let $W$ be a component of $s_{\bbv}^{-1}(C_V)$ contained in $Z$. 
For $(x, z_I) \in W$,  assume that $I$ is supported at $p={\rm pr}_1(x) \in D_0$,  $p'\in D_\infty$ 
(or $p \in D_\infty$,  $p'\in D_0$) 
with $\pi_P(p)=\pi_P(p')$, and other $p_3, \cdots, p_r$ in $P \sqcup \tilde{X}_t$. 
From the third remark (c), the singularity type of $h_P$ at $p$ (resp. $p'$) and that of $g$ at $q$ (resp. $q'$) are the same, say $\eta_1$ (resp. $\eta_2$).  
Other points $p_j$ correspond to other mono-singularity types $\eta_j$'s in $\underline{\eta}$ ($3\le j \le r$). 
If $p_j \in \tilde{X}_t$, then it has the limit in $Q$ or $\tilde{A}\sqcup \tilde{B}$ as $t \to 0$.

\begin{itemize}
\item[(d)] 
If some $p_j$ lies on $P$ and $\psi(p_j) \in Q$ lies on some component $\overline{\mathcal{Y}}_k$ of the cycle representing the class $\tilde\alpha$, then there is some $\hat{p}_j \in \tilde{X}_t$ in the locus satisfying that there are $\hat{p}, \hat{p}' \in D_0 \sqcup D_\infty$ with $\pi_P(\hat{p})=\pi_P(\hat{p}')$ and other points such that the multi-germ of $h$ at those points is of type $\underline{\eta}$. This condition characterizes another component of $s_{\bbv}^{-1}(C_V)$ contained in $Z$, say $W'$, which is different from $W$. From a similar argument as the above proof of the equality (1) using Lemma \ref{change_sign}, we can see that ${\rm pr}_{1*}[W]=-{\rm pr}_{1*}[W']$, and hence, the contributions from those components $W$ and $W'$ to the refined class $\delta_d$ cancel each other out. 
\item[(e)] 
If every $p_j$ lies on $\tilde{X}_t$ and has the limit in $\tilde{A} \sqcup \tilde{B}$, the limit of $p$ (the first entry), regarded as in $D_0' \sqcup D_\infty'$, belongs to the $\underline{\eta}$-multi-singularity locus of $g: \tilde{A} \sqcup \tilde{B} \to Y$. Namely, it is a point in the $\eta_1$-locus of $g$ accompanied by other points in $\tilde{A} \sqcup \tilde{B}$ of prescribed types with the same image in $Y$. This condition characterizes the component $W$, that is, if $W$ consists of points $(x, z_I)$ with this condition, ${\rm pr}_{1*}[W]$ contributes to the part of the loci class $m_{\underline{\eta}}(g)$ supported on $D_0 \sqcup D_\infty$. 
\item[(f)] 
Besides from (d) and (e), if all $p_j$ are paired so that each pair is just the fiber of $D_0 \sqcup D_\infty \to D$ at a point of $D$, then ${\rm pr}_{1*}[W]$ also contributes to the part of $m_{\underline{\eta}}(g)$ supported on $D_0 \sqcup D_\infty$. 
\end{itemize}

By the argument in (d), we may consider only components $W$ in cases (e) and (f). 
Consequently, the refined class $\delta_D$ becomes part of $m_{\underline{\eta}}(g)$. 
The remaining part $m_{\underline{\eta}}(g) - \delta_D$ should come from the $\underline{\eta}$-singularities of $g$ off $D_0 \sqcup D_\infty'$, that is, the intersection $(\tilde{A} + \tilde{B}) \cdot \tilde{\alpha}$. 
Thus we obtain (2). This completes the proof. \qed

\

\subsubsection{\bf  Algebraic cobordism setting} \label{Omega_Q}
We prove a stronger version of Proposition \ref{degeneration_free}. 
Let us choose and fix an $\Omega$-assignment $\mu$ associated to $\underline{\eta}$. 
We will show that  
$$n_{\underline{\eta}; \mu}^\Omega(g)=n_{\underline{\eta}; \mu}^\Omega(h) \;\; \in \Omega^*(Y).$$ 

The proof essentially goes along the same line as seen above but with some technical remarks. 
A main point is to lift up the equality $[\tilde{X}_t]=[\tilde{A}]+[\tilde{B}]+2[Q]$ in $\CH^*(\tilde{M})$ to the level of algebraic cobordism. 
Let  
$$i_0: D'_0 \hookrightarrow Q, \quad i_\infty: D'_\infty \hookrightarrow Q, \quad j: Q \hookrightarrow \tilde{M}$$ 
be the inclusions, and $D'_0, D'_\infty$ are naturally identified with $D$ (e.g., $P \to  D'_0$ is identified with $P \to  D$ through $D_0' \simeq D$). We then set  
$$\Theta:=[P \to  Q]-(i_{0})_*[P \to  D'_0]-(i_{\infty})_*[P\to D'_\infty] \in \Omega^*(Q).$$
Note that the cobordism class $\Theta$ becomes to be $2Q$ after applying $\Omega^* \to \CH^*$. 

\begin{prop}\label{dpr_blowup}
In $\Omega^*(\tilde{M})$, it holds that 
$$[\tilde{X}_t \to \tilde{M}]=[\tilde{A} \to  \tilde{M}] +  [\tilde{B} \to  \tilde{M}]+ j_*\Theta.$$
\end{prop}

\proof 
We refer to \S 3.1 and \S 3.2 in \cite{LevineMorel}. 
We set  
$$E:=\tilde{A}+\tilde{B}+2 \tilde{Q},$$
which is a strict normal crossing divisor in $\tilde{M}$  \cite[Def.3.1.4]{LevineMorel}. 
Note that $\tilde{A} \cap \tilde{B}=\emptyset$. 
Put $L_A:=O_{\tilde{M}}(\tilde{A})$, $L_B:=O_{\tilde{M}}(\tilde{B})$, $L_Q:=O_{\tilde{M}}(Q)$, line bundles of those divisors. 

First, by \cite[Prop.3.1.9]{LevineMorel}, we have 
$$[\tilde{X}_t \to \tilde{M}]=|O_{\tilde{M}}(E)|
=c_1(L_A\otimes L_B \otimes L_Q^{\otimes 2}) \;\; \in \Omega^*(\tilde{M}).$$ 
Since the singular locus $|E|_{\rm sing}$ is just $Q$, it follows from \cite[Lem.3.2.2]{LevineMorel} that there exists some $\theta \in \Omega^*(Q)$ such that 
$$
[\tilde{X}_t \to \tilde{M}] = [\tilde{A} \to \tilde{M}]+ [\tilde{B} \to \tilde{M}] + j_*\theta
\;\; \in \Omega^*(\tilde{M}). 
\eqno{(1)}
$$
We show $\theta=\Theta$ below.  
Expanding $\theta$ by the formal group law, 
we may write\footnote{
Actually, $\theta'$ is the sum of terms $a_{ijk} c_1(j^*L_A)^ic_1(j^*L_B)^jc_1(j^*L_Q)^k$ 
with $a_{ijk} \in \bL$, $i+j\ge 1$, $k\ge1$ and $a_{i00}=a_{0j0}=0$ in the expansion by the formal group law $F_{\bL}$. } 
$$
\theta = 2 [id_Q] + \theta' \;\;\; \in \Omega^*(Q)
$$ 
for some $\theta' \in \Omega^*(Q)\cdot \gamma$, where $\gamma:=c_1(O_Q(1))=c_1(j^*L_Q)$, 
and $id_Q$ is the identity map $Q \to Q$. 
Note that 
$$
\tilde{\pi}_*[\tilde{X}_t \to \tilde{M}]
=[{X}_t \to M]= [A \to M]+[B \to M]-[P \to M] \;\; \in \Omega^*(M)
$$
by the double-point relation over $M$. 
Comparing this equality with $(1)$, we get 
$$\tilde{\pi}_*j_*\theta=-[P \to M].$$ 
Here, we may take the normal bundle of $D$ in $M$, instead of $M$, 
and then
$$(\pi_{Q})_*\theta=-[P \to D]\eqno{(2)}$$ 
where $\pi_{Q}=\tilde{\pi}|_Q: Q \to D$ is the projection. 

Since $\pi_{Q}: Q \to D$ is a projective bundle, $\Omega^*(Q)$ is a free $\Omega^*(D)$-module of rank $2$ \cite{HV11, Vishik07}:  
$$\Omega^*(Q)\simeq \Omega^*(D)\cdot [id_Q]\oplus \Omega^*(D)\cdot \gamma\eqno{(3)}$$ 
Now $\psi: P \to Q$ is a $2$-to-$1$ surjective map branched along $D_0' \sqcup D_\infty'$. 
It then follows from (3) that there is an element $\xi \in \Omega^*(D)\cdot \gamma$ such that  
$$
[P \to Q] = 2 [id_Q] + \xi \;\;\; \in \Omega^*(Q)
$$
(this also follows from the degree formula in \cite{LevineMorel}). 
Hence, $\theta=[P \to Q]+\theta'-\xi$. For $(\pi_Q)_*[P \to Q]=[P \to D]$ and $(2)$, we see 
$$[P \to Q]+\theta=2[P \to Q] +\theta' -\xi\in \ker(\pi_Q)_*.$$ 
Here $\theta'-\xi \in \Omega^*(D)\cdot \gamma$ and $(\pi_Q)_*(\theta'-\xi)=-2[P \to D]$. 
Thus we may write it by 
$$
\theta'-\xi=-2(i_0)_*[P \to D_0']
=-(i_0)_*[P \to D_0']-(i_\infty)_*[P \to D_\infty']. 
$$
That implies that $\theta=\Theta$. 
\qed 

\

Note that every ingredient defined in $\CH^*$ in the previous subsections can be formulated in algebraic cobordism setting -- not only multi-singularity loci classes $m^\Omega_{\underline{\eta}; \mu}$ but also the classes $\tilde{\alpha}$,  $\tilde{\beta}$ and $\delta_D$ are defined  in $\Omega^{\ell-\kappa}(\tilde{M})$, $\Omega^{\ell-\kappa}(P)$ and $\Omega^*(D_0' \sqcup D'_\infty)$, respectively. 
In entirely the same way as in \S \ref{refined}, we get a cobordism version of Lemma \ref{lem3}: 

\begin{lem}
In $\Omega^{\ell-\kappa}(\tilde{M})$, it holds that 
\begin{align*}
&j_*\psi_*m^\Omega_{\underline{\eta}; \mu}(h)^P = j_*\delta_D - \Theta  \cdot \tilde{\alpha} \tag{1}\\
&(i_g)_*m^\Omega_{\underline{\eta}; \mu}(g) = j_*\delta_D + ([\tilde{A}\to \tilde{M}]+[\tilde{B}\to \tilde{M}])\cdot \tilde{\alpha}  \tag{2} 
\end{align*}
\end{lem}

Consequently, we see that 
$$(i_g)_*(m^\Omega_{\underline{\eta}; \mu}(g))= (i_h)_*(m^\Omega_{\underline{\eta}; \mu}(h)) \in \Omega^{\ell - \kappa+1}(\tilde{M}).$$
Pushing the both sides to $Y$ via $(F_1)_*$ yields 
$$n^{\Omega}_{\underline{\eta}; \mu}(g) = n^{\Omega}_{\underline{\eta};\mu}(h) \in \Omega^*(Y).$$
This completes the proof of Proposition \ref{degeneration_free} valued in $\Omega^*$.  
\qed

\

\subsection{Target multi-singularity loci classes} \label{target_tp}
We prove Theorem \ref{main_thm1} along the line mentioned at the step (v) in Introduction. 

\subsubsection{\bf Residual class operation} 

\begin{definition}\label{delta_n}\upshape
Let $f: X \to Y$ be a morphism of codimension $\kappa$ in $\Sm$. 
\begin{enumerate}
\item 
For a type $\eta$ of mono-germs, 
we put 
$$\Delta n_{\eta}(f):=f_*(m_{\eta}(f))=n_{\eta}(f) \;\; \in \CH^{\ell(\eta)}(Y).$$
\item 
For a multi-singularity type $\underline{\eta}=(\eta_1, \cdots, \eta_r)$ with $\ell:=\ell(\underline{\eta})$ $(r\ge 2)$, 
we recursively define 
$$\Delta n_{\underline{\eta}}(f):=
n_{\underline{\eta}}(f)-\sum  \, \Delta n_{J_1}(f)\cdots  \Delta n_{J_s}(f) \;\; \in \CH^\ell(Y)$$
where the summand runs over all partitions into a disjoint union of at least two non-empty unordered subsets, 
$\{1, \cdots, r\}= J_1\sqcup \cdots \sqcup J_s$ $(s\ge 2)$. 
\end{enumerate}
 \end{definition}
Let $f: X_1 \to Y$ and $g: X_2 \to Y$ be morphisms in $\Sm$ with 
$\dim X_1=\dim X_2$, and 
$f + g: X_1 \sqcup X_2 \to Y$ be the disjoint sum of two morphisms. 

For a mono-singularity type $\alpha$ ($r=1$), 
the geometric subset $\Xi(X_1 \sqcup X_2; \alpha)$ is the disjoint union of 
$\Xi(X_1; \alpha) \subset X_1^{[[n]]}$ and 
$\Xi(X_2; \alpha) \subset X_2^{[[n]]}$ (Lemma \ref{H_disjoint}), 
and thus it follows from the definition of degeneracy class of maps (Defnition \ref{multi-singularity_class}) that 
$$n_\alpha(f + g) = n_\alpha(f)+n_\alpha(g).$$ 
As a double-check, 
this equality follows from the fact that the Thom polynomial for mono-singularity type is a polynomial of $c_i(f)$'s. 

In case of $r=2$, for a bi-germ $\underline{\eta}=(\alpha, \beta)$, we can show that 
$$n_{\alpha\beta}(f+g)=n_{\alpha\beta}(f)+n_{\alpha\beta}(g)+n_\alpha(f)n_\beta(g)+n_\alpha(g)n_\beta(f)$$
and therefore  
\begin{eqnarray*}
\Delta n_{\alpha\beta}(f+g)&=&
n_{\alpha\beta}(f+g)-n_\alpha(f+g) \cdot n_\beta(f+g)\\
&=& n_{\alpha\beta}(f+g)-(n_\alpha(f)+n_\alpha(g))(n_\beta(f)+n_\beta(g))\\
&=& n_{\alpha\beta}(f)-n_\alpha(f)n_\beta(f)+n_{\alpha\beta}(g)-n_\alpha(g)n_\beta(g)\\
&=& \Delta n_{\alpha\beta}(f)+\Delta n_{\alpha\beta}(g).
\end{eqnarray*}
More generally, the following proposition holds: 

\begin{prop}\label{prop1}
The class $\Delta n_{\underline{\eta}}(f)$ depends only on the cobordism class of $f$, and it yields a well-defined group homomorphism  
$$\Delta n_{\underline{\eta}}:\Omega^{\kappa}(Y) \to \CH^\ell(Y), \quad [f] \mapsto \Delta n_{\underline{\eta}}(f).$$ 
In particular, $n_{\underline{\eta}}(f)$ also depends only on the cobordism class of $f$. 
\end{prop}

\proof 
Note that $\Delta n_{\underline{\eta}}: \Mcal^\kappa(Y) \to \CH^\ell(Y)$ is well-defined as a map between sets, because 
$n_{\underline{\eta}}(f)$ is invariant under the action of isomorphisms of the source spaces (Lemma \ref{iso} (2)). 
First, we show  
\begin{center}
Claim: 
$\Delta n_{\underline{\eta}}$ is additive on the monoid $\Mcal^\kappa(Y)$. 
\end{center}
We use the induction on $r$. 
In mono-singularity case ($r=1$), the claim is clear as seen above. 
Let $r \ge 2$ and 
suppose that the claim holds for any $r'$-tuple multi-singularities with $r'<r$. 

Now let $\underline{\eta}=(\eta_1, \cdots, \eta_r)$, a multi-singularity type. 
Let $f: X_1 \to Y$ and $g: X_2 \to Y$ be maps in $\Sm$ with the same codimension $\kappa$ and the same target space. 

We recall Lemma \ref{H_disjoint}: in $(X_1\sqcup X_2)^{[[n]]}$ 
$$\Xi(X_1\sqcup X_2;\underline{\eta}) = 
\bigsqcup_{I\sqcup J} \Xi(X_1; \underline{\eta}_{I}) \times \Xi(X_2; \underline{\eta}_{J}).$$
Consider the intersection product of the Hilbert extension map 
$$(f + g)^{[[n]]}: (X_1 \sqcup X_2)^{[[n]]} \to ((X_1 \sqcup X_2) \times Y)^{[[n]]}$$
with the subscheme $\Xi(X_1 \sqcup X_2; \underline{\eta})\times Y$.   
Pushing it out to $Y$, we have 
$$n_{\underline{\eta}}(f+g)=n_{\underline{\eta}}(f)+n_{\underline{\eta}}(g)+\sum_{I\sqcup J} n_I(f)\cdot n_J(g),$$
where $I \sqcup J = \{1, \cdots, r\}$ and both $I$ and $J$ are {\it non-empty}. It is for $n_I(f)\cdot n_J(g) \in \CH^\ell(Y)$ 
is the pull back of $n_I(f) \times n_J(g)\in \CH^*(Y \times Y)$ 
via the diagonal embedding $Y \hookrightarrow Y\times Y$. 
By Definition \ref{delta_n}, it holds that 
$$n_I(f)=\sum \Delta n_{I_1}(f) \cdots \Delta n_{I_{t}}(f)$$ 
where the summand runs over all partitions of $I$ into non-empty unordered subsets ($t\ge 1$), 
and $n_J(g)$ is also. Therefore, we see that 
$$\sum_{I\sqcup J} n_I(f)\cdot n_J(g)=\sum_{(J_i; h_i)} \Delta n_{J_1}(h_1)\cdots  \Delta n_{J_s}(h_s)$$
where $J_1\sqcup \cdots \sqcup J_s \; (s\ge 2)$ are partitions of $\{1, \cdots, r\}$ 
into non-empty unordered subsets and for each $i$, 
$h_i$ denotes either $f$ or $g$ but all $h_1, \cdots, h_s$ are not the same simultaneously. 
Now, since each $J_i$ satisfies $|J_i| < r$, 
the induction hypothesis implies that for each $1\le i \le s$, 
$$\Delta n_{J_i}(f+g)=\Delta n_{J_i}(f)+\Delta n_{J_i}(g).$$ 
Thus we see 
\begin{eqnarray*}
&&\Delta n_{\underline{\eta}}(f+g)\\
&&=
n_{\underline{\eta}}(f+g)-\sum  \, \Delta n_{J_1}(f+g)\cdots  \Delta n_{J_s}(f+g)\\
&&=
n_{\underline{\eta}}(f)+n_{\underline{\eta}}(g)+\sum_{I\sqcup J} n_I(f)\cdot n_J(g) 
- \sum \prod (\Delta n_{J_i}(f)+\Delta n_{J_i}(g)) \\
&&= \left(n_{\underline{\eta}}(f)-\sum  \, \prod \Delta n_{J_i}(f)  \right)
+ \left(n_{\underline{\eta}}(g)-\sum  \, \prod \Delta n_{J_i}(g)  \right)\\
&&=\Delta n_{\underline{\eta}}(f)+\Delta n_{\underline{\eta}}(g).
\end{eqnarray*}
Hence, the claim is proven. 

It immediately yields a group homomorphism 
$$\Delta n_{\underline{\eta}}: \Mcal^{\kappa+}(Y) \to \CH^\ell(Y).$$ 
Recall that $\Rcal(Y)$ is generated by the double-point relation: 
$$[g: A \sqcup B \to Y]  - [h: X_1 \sqcup P \to Y]  \in \Mcal^{\kappa+}(Y).$$
Proposition \ref{degeneration_free} says that $n_{\underline{\eta}}(g)=n_{\underline{\eta}}(h)$. 
From Definition \ref{delta_n}, by the induction on $r$, we see that  
$$\Delta n_{\underline{\eta}}(g-h)=\Delta n_{\underline{\eta}}(g)-\Delta n_{\underline{\eta}}(h)=0,$$
i.e., the homomorphism $\Delta n_{\underline{\eta}}$ vanishes on $\Rcal(Y)$. 
Thus we have a well-defined group homomorphism 
$$\Delta n_{\underline{\eta}}: \Omega^\kappa(Y) \to \CH^\ell(Y), \quad 
[f] \mapsto \Delta n_{\underline{\eta}}(f).$$ 
This completes the proof. \qed

\

\begin{prop}\label{prop2}
$\Delta n_{\underline{\eta}}:\Omega^{\kappa}(Y) \to \CH^\ell(Y)$ is an additive cohomology operation. 
\end{prop}

\proof 
Thanks to Vishik's result (Theorem \ref{Vishik1}), 
it suffices to show that for $Y=\prod \Proj^{n_j}$, product of projective spaces,  
$\Delta n_{\underline{\eta}}$ commutes with the pullback 
$$\varphi^*: \Omega^*(Y) \to \Omega^*(Y')$$ 
for any morphisms $\varphi: Y' \to Y$ listed in Theorem \ref{Vishik1}. 
Let $Y$ and $\varphi$ be so. 
Take an arbitrary cobordism class $[f: X \to Y] \in \Omega^\kappa(Y)$. 
By moving by translation, 
$f$ can be linearly deformed to be transverse to $\varphi$ so that the cobordism class does not change (Remark \ref{pullback_omega}). 
So we may assume $f$ itself admits the transversality, 
then the fibre square 
$$f': X \times_Y Y'\to Y'$$ 
exists in $\Sm$ and satisfies $\varphi^*([f])=[f']$ (Remark \ref{pullback_omega}). 
Lemma \ref{pullback} says that $n_{\underline{\eta}}$ satisfies the naturality for base-change.  
From Definition \ref{delta_n} combined with the induction on $r$, 
we see that $\Delta n_{\underline{\eta}}$ also satisfies the same base-change property: 
$\Delta n_{\underline{\eta}}([f'])=\varphi^*(\Delta n_{\underline{\eta}}([f]))$. 
Hence we have 
$$\Delta n_{\underline{\eta}}\circ \varphi^*=\varphi^*\circ \Delta n_{\underline{\eta}}$$
 for any $[f]\in \Omega^\kappa(Y)$. This completes the proof. 
\qed 

\


We immediately have

\begin{cor}\upshape \label{prop3}
For each multi-singularity type $\underline{\eta}$, 
$$n_{\underline{\eta}}:\Omega^{\kappa}(Y) \to \CH^\ell(Y)$$ 
is a well-defined cohomology operation (not additive if $r\ge 2$). 
\end{cor}

\subsubsection{\bf Theorem in cobordism setting} \label{main_thm}
Now we switch from $\CH^*$ to $\Omega^*$. 
Let $\underline{\eta}=(\eta_1, \cdots, \eta_r)$ be a multi-singularity type of map-germ with codimension $\kappa$. 
Suppose that we are given an $\Omega$-assignment $\mu$ for the type $\underline{\eta}$  (Definition \ref{good}). 
For an arbitrary proper morphism $f: X \to Y$, $n_{\underline{\eta}; \mu}^\Omega(f) \in \Omega^\ell(Y)$ is defined (Definition \ref{multi-singularity_class_Omega}). 
We show the following result: 

\begin{thm}\label{main_thm}
{\bf (Target cobordism Thom polynomials)} 
For every $J \subset \{1, \cdots, r\}$, there is a residual polynomial $R_{J}$ of $c_1, c_2, \cdots$ with coefficients in $\bL\otimes_Z \Q$ such that for an arbitrary proper morphism $f: X \to Y$ with codimension $\kappa$, the target multi-singularity loci cobordism class is expressed by 
$$n_{\underline{\eta}; \mu}^\Omega (f) =\sum f_*(R_{J_1})\cdots f_*(R_{J_s})$$
in $\Omega^\ell(Y)_\Q$, where the summand runs over all partitions. 
\end{thm}

\proof
It immediately follows that 
\begin{itemize}
\item[-] 
Proposition \ref{degeneration_free} holds also in $\Omega^*$-setting (proven in \S \ref{Omega_Q}); 
\item[-] 
the class $\Delta n_{\underline{\eta}; \mu}^\Omega(f)$ is also defined (Definition \ref{delta_n}); 
\item[-] 
Propositions \ref{prop1} and \ref{prop2} in $\Omega$-setting are shown in entirely the same way, and thus 
$$\Delta n_{\underline{\eta}; \mu}^\Omega:\Omega^{\kappa}(Y) \to \Omega^\ell(Y)$$
is a well-defined additive operation.  
\end{itemize}
Now, we are ready to employ Theorem \ref{Vishik2} (Vishik  \cite[Thm.1.2]{Vishik}). 
We conclude that 
the operation $\Delta n_{\underline{\eta}; \mu}^\Omega$ is uniquely written by a linear combination of Landweber-Novikov operations: 
$$\Delta n_{\underline{\eta}; \mu}^\Omega=\sum a_I s_I \quad (a_I \in \bL\otimes_\Z \Q).$$ 
Put $R_{\underline{\eta}}:=\sum a_Ic^I$, then 
$$f_*(R_{\underline{\eta}}(f))=\Delta n_{\underline{\eta}; \mu}^\Omega([f]).$$ 
Hence, by the induction on $r$ and the definition of $\Delta n_{\underline{\eta}; \mu}^\Omega([f])$, we get the desired expression of the class $n_{\underline{\eta}; \mu}^\Omega (f)$. Namely, we have shown the existence of the (target) cobordism Thom polynomial for an $\Omega$-assignment $\mu$ associated to $\underline{\eta}$. \qed

\subsubsection{\bf Proof of Theorem \ref{main_thm1}} 
For any cohomology theory $A^*$ obtained from $\Omega^*$, we define the $A^*$-multi-singularity loci classes associated to an $\Omega$-assignment $\mu$. 
In particular, in case of $A^*=\CH^*$ (and also $K^0[\beta, \beta^{-1}]$), the classes depend only on $\underline{\eta}$, not the choice of $\mu$: 
$$
\xymatrix{
\Omega^\kappa \ar[r]^{n_{\underline{\eta}; \mu}^\Omega} \ar[dr]_{n_{\underline{\eta}}} & \Omega^\ell \ar[d]^{\otimes_{\bL} \Z}\\
& \CH^\ell
}
$$
In Theorem \ref{main_thm}, it is shown that $n_{\underline{\eta}; \mu}^\Omega$ is expressed by a universal multi-singularity Thom polynomial. After reducing to $\CH^*$ (or  $K^0[\beta, \beta^{-1}]$), 
such universal polynomials for any choices of $\mu$ coincide with each other.  
This completes the proof of Theorem \ref{main_thm1}. 
\qed

\

\subsection{Source multi-singularity loci classes} \label{source_tp}
We prove Theorem \ref{main_thm2} by using basic materials of the restriction method due to R. Rim\'anyi \cite{Rimanyi01, Rimanyi02, Kaz03, Kaz06, RimanyiSzucs, Ohmoto06, Ohmoto12}. 
The proof uses Thoerem \ref{main_thm1} but also relies on a different nature which comes from {\em local symmetry of singularities and equivariant cohomology}.  
That suggests there should be a comprehensive new theory, see \S \ref{classifying_stack} and \S \ref{motivic}. 
Note that in the statement of Theorem \ref{main_thm2}, the case of $\kappa < 0$ is excluded (the main reason is mentioned in \S \ref{critical}). 

Now, fix $\kappa=l-m \ge 0$. We work with $\CH^*_\Q$ throughout this section. 

\begin{rem}\upshape \label{qh}
{\bf (Quasi-homogenious normal forms)}
We say that a multi-singularity type $\underline{\eta}=(\eta_1, \cdots, \eta_r)$ is {\em quasi-homogeneous} if every entry $\eta_i$ admits quasi-homogeneous polynomial normal forms in suitable local coordinates (under $\K$-equivalence). For instance, considering the classification list  \cite[\S 5]{Mather} and a standard construction of stable unfoldings of finitely determined map-germs \cite{MatherIV}, it turns out that every locally stable map-germ in {\em Mather's nice dimensions} is of quasi-homogeneous type, see  \cite{MondNuno}. 
Below we only consider such singularity types (cf. admissible maps in Theorem \S \ref{main_thm2}). 

\end{rem}

\subsubsection{\bf Symmetry of singularity type} 
We denote the group of formal right-left equivalence by 
$$\hat{\A}_S:=Aut(\hat{\Ost}_{\bA^m, S}) \times Aut(\hat{\Ost}_{\bA^l, 0})$$ 
($S \subset \bA^m$ with $|S|=s$), which acts on the space of map-germs $f: (\bA^m, S) \to (\bA^l, 0)$ by $(\sigma, \varphi).f:= \varphi \circ f \circ \sigma^{-1}$. If $s=1$ (the case of mono-germs), we simply write it by $\hat{\A}$. Note that any quasi-homogeneous germ has the isotropy subgroup in $\hat{\A}_S$ which contains the {\em non-trivial} maximal torus, because there is at least a non-trivial $GL(1)$-action on the source and the target space preserving the map-germ \cite{Rimanyi01, Kaz03}. 

\

\t
{\it Mono-singularities}: 
Let $\xi$ be a quasi-homogeneous mono-singularity type of codimension $\kappa$,  and set $m':=\ell(\xi)-\kappa$ and $l':=\ell(\xi)$. Consider a mono-germ of locally stable map of that type; abusing the notation, we denote it by 
$$\xi: (\bA^{m'}, 0) \to (\bA^{l'}, 0).$$
Abusing the notation, we use the same letter $\xi$ to mean the singularity type as well. 
The {\em symmetry group $G(\xi)$ of the singularity type $\xi$} is defined to be the isotopy subgroup of $\hat{\A}$ for $\xi$, i.e., it consists of elements $(\sigma, \varphi) \in \hat{\A}$ with $(\sigma, \varphi).\xi=\xi$; in other words, it is defined by the group of automorphims of the quotient formal algebra of $\xi$. 
Since any locally stable map-germ is $(l'+1)$-determined \cite{MatherIV}, we only have to work at the level of $(l'+1)$-jets. Thus we may assume that the symmetry group is realized as a subgroup of a certain linear algebraic group. 
There are linearized representations $\lambda_0: G(\xi) \to GL(m')$ and $\lambda_1: G(\xi) \to GL(l')$ into general linear groups acting on the source and the target spaces $\bA^{m'}$ and $\bA^{l'}$, respectively. 
By the algebraic Borel construction \cite{Totaro}, there exist the classifying space $BG(\xi)$ of the group $G(\xi)$ as an ind-variety, and it admits the universal principal $G(\xi)$-bundle $EG(\xi) \to BG(\xi)$. 
Put 
$$E_0(\xi):= EG(\xi) \times_{\rho_0} \bA^{m'}, \quad E_1(\xi):= EG(\xi) \times_{\rho_1} \bA^{l'},$$ 
which are vector bundles of rank $m'$ and $l'$ over $BG(\xi)$ associated to the representations $\lambda_0$ and $\lambda_1$, respectively. 
Since the germ $\xi$ is invariant under the action of $G(\xi)$, we can construct a morphism $f_\xi$ between total spaces of those bundles so that the restriction on each fiber is $\A$-equivalent to $\xi$: 
$$
\xymatrix{
E_0(\xi) \ar[rr]^{f_{\xi}} \ar[dr]_{\pi_0} && E_1(\xi) \ar[dl]^{\pi_1}\\
& BG(\xi)&\\
}
$$

\begin{rem}\upshape \label{max_torus}
By Remark \ref{qh}, $G(\xi)$ has the non-trivial maximal torus, thus especially, {\em the top Chern classes of $E_0(\xi)$ and $E_1(\xi)$ are non-zero divisors}. Since $BGL(1)=\Proj^\infty$, we may think of $BG(\xi)$ as a projective smooth variety of dimension high enough. 
\end{rem}

\

\t
{\it Multi-singularities}:  
Let $\underline{\xi}$ be a quasi-homogenous multi-singularity type of codimension $\kappa$,  and set $m':=\ell(\underline{\xi})-\kappa$ and $l':=\ell(\underline{\xi})$, i.e., 
$$\underline{\xi}=(\xi_1, \xi_2, \cdots, \xi_s), \quad \xi_i: (\bA^{m'}, 0) \to (\bA^{l'}, 0).$$
Here, especially the first entry $\xi_1$ is distinguished from others; 
denote the rest $(s-1)$-tuple by $\underline{\xi}/\xi_1:=(\xi_2, \cdots, \xi_s)$. 
Let $Aut(\underline{\xi}) \subset \mathfrak{S}_s$ denote the subgroup of permutations of $\{1, \cdots, s\}$ preserving the collection $\underline{\xi}$, and we regard $Aut(\underline{\xi}/\xi_1) \subset Aut(\underline{\xi})$ as a subgroup fixing $1$  (see Proposition \ref{multi-sing})\footnote{For instance, if $\underline{\xi}=(A_0, A_1, A_0, A_0, A_2)$, $Aut(\underline{\xi})=\{\, \sigma \in  \mathfrak{S}_5\, |\, \sigma(2)=2, \sigma(5)=5\, \} \simeq \mathfrak{S}_3$ and $Aut(\underline{\xi}/\xi_1)=\{\, \sigma  \, | \,  \sigma(1)=1, \sigma(2)=2, \sigma(5)=5\, \} \simeq \mathfrak{S}_2$. }. 

The semi-direct product $ \hat{\A}_S \ltimes  \mathfrak{S}_s$ acts on the space of all ordered $s$-tuples of mono-germs. Let $G(\underline{\xi}) \subset \hat{\A}_S \ltimes Aut(\underline{\xi})$ denote the isotropy subgroup at $\underline{\xi}$. 
Also we define the smaller isotropy subgroup $G(\underline{\xi}/\xi_1)$ by the intersection of $G(\underline{\xi})$ with $\hat{\A}_S \ltimes Aut(\underline{\xi}/\xi_1)$. 
Let $\lambda_0$ be the linearized representation of $G(\underline{\xi}/\xi_1)$ on the source domain $\bA^{m'}$ of the first (fixed) mono-germ $\xi_1$, and $\lambda_1$ be the representation of $G(\underline{\xi})$ on the common target domain $\bA^{l'}$ of the multi-germ $\underline{\xi}$. 
By the Borel construction again, we associate the universal vector bundles 
$$X_{\underline{\xi}}:=EG(\underline{\xi}/\xi_1)\times_{\lambda_0} \bA^{m'}, \quad 
Y_{\underline{\xi}}:=EG(\underline{\xi})\times_{\lambda_1} \bA^{l'}$$
with the base spaces being the classifying spaces $BG(\underline{\xi}/\xi_1)$ and $BG(\underline{\xi})$, respectively. 
The inclusion map between the two groups induces a morphism 
$$\rho_{\underline{\xi}}: BG(\underline{\xi}/\xi_1) \to BG(\underline{\xi}),$$ 
which is a $d$-to-$1$ covering with 
$$d:=\deg_1(\underline{\xi})=\#\Aut(\underline{\xi})/\#\Aut(\underline{\xi}/\xi_1).$$ 
In fact, the algebraic Borel construction uses finite dimensional approximation of the classifying spaces \cite{Totaro}, so we may regard $\rho_{\underline{\xi}}$ as a covering of schemes of sufficiently high dimension.  
Furthermore, over the map $\rho_{\underline{\xi}}$, we obtain a universal singular map $f_{\underline{\xi}}$ between the total spaces of the above vector bundles 
$$
\xymatrix{
X_{\underline{\xi}} \ar[r]^{f_{\underline{\xi}}} \ar[d]_{\pi_0} & Y_{\underline{\xi}}\ar[d]^{\pi_1}\\
BG(\underline{\xi}/\xi_1) \ar[r]_{\rho_{\underline{\xi}}} &  BG(\underline{\xi})\\
}
$$
by the gluing construction such that the restriction of $f_{\underline{\xi}}$ to every fiber is $\A$-equivalent to the mono-germ $\xi_1$. Note that $f_{\underline{\xi}}$ itself is proper, for $\kappa \ge 0$. 
Thus the pushforward $(f_{\underline{\xi}})_*$ is defined. 

\begin{lem}\upshape \label{injectivity}
$(f_{\underline{\xi}})_*: \CH^{*}(X_{\underline{\xi}})_\Q \to \CH^{*+\kappa}(Y_{\underline{\xi}})_\Q$ 
is injective.\end{lem}

\proof 
The groups $G(\underline{\xi}/\xi_1)$ and $G(\underline{\xi})$ share the common non-trivial maximal torus, and thus 
$\CH^*(X_{\underline{\xi}})_\Q \simeq \CH^*(Y_{\underline{\xi}})_\Q$ via the pullback $f_{\underline{\xi}}^*$, moreover, those are isomorphic to $\CH^*(BG(\underline{\xi}))_\Q$. Note that $f_{\underline{\xi}* }(f_{\underline{\xi}}^*\omega)=f_{\underline{\xi}*}(1)\cdot \omega$ and the top Chern classes of $X_{\underline{\xi}}$ and $Y_{\underline{\xi}}$ are non-zero divisors (Remark \ref{max_torus}). 
Further, the zero section of $X_{\underline{\xi}}$ is sent via $f_{\underline{\xi}}$ to the zero section of $Y_{\underline{\xi}}$, i.e., the restriction map between the zero sections is identical to $\rho_{\underline{\xi}}$. 
Hence, we have an obvious relation between the top Chern classes 
$$f_{\underline{\xi}*}(1) \cdot c_{m'}(X_{\underline{\xi}})=d\, c_{l'}(Y_{\underline{\xi}})$$
in $\CH^{l'}(BG(\underline{\xi}))_\Q$ (cf. \cite[Lem.4.1]{Kaz06}), and especially, $f_{\underline{\xi}*}(1)$ is a non-zero divisor. 
The injectivity of $f_{\underline{\xi}*}$ follows from this relation. 
 \qed

\

The constructed universal maps $f_\xi$ and $f_{\underline{\xi}}$ are building blocks of the restriction method (Remark \ref{restriction}).

\subsubsection{\bf Stratifications} \label{stratification}
From now on, we consider an {\em admissible} map in $\Sm$ 
$$f:X \to Y \quad (m=\dim X, \;\; l=m+\kappa=\dim Y,\;\; \kappa \ge 0).$$
That is, $f$ is a proper locally stable map which has only quasi-homogeneous singularity types and satisfies `stratawise normal Euler condition'. 
The latter condition will be used in a Mayer-Vietoris argument arising in the final step of \S \ref{proof_thm2}. 

Since $f$ is locally stable, we obtain a stratification of each of $X$ and $Y$ into several multi-singularity loci in a coherent way  subject to $f$.  Every stratum is reduced, locally closed, but {\em not assumed to be connected}.  Also every stratum satisfies the frontier condition -- its closure is the union of some strata of lower dimension (called boundary strata). 
These statifications are indexed by multi-singularity types in the following sense\footnote{In case of $\kappa>0$, $M_{A_0}$ for mono-type $A_0$ of immersion-germ is the unique open stratum in $X$, and every positive codimensional stratum (of type $A_0^2$, $A_1$ and so on) is adjacent to it; especially, $M_{A_0}^\circ \to N_{A_0}^\circ$ is isomorphic. In case of $\kappa=0$, we also include the type $A_0$, although our original definition of $\underline{\xi}$ excludes this type, and we always assume that our equidimensional proper map $f: X \to Y$ is dominant; especially,  $M_{A_0}^\circ \to N_{A_0}^\circ$ is a $d$-to-$1$ covering, where $d$ is the degree of $f$. }. 
The source manifold is stratified as 
$$X=\bigsqcup M_{\underline{\xi}}^\circ$$
where the union is taken over all mono or {\em semi-ordered} multi-singularity types of codimension $\kappa$ (that means that the first entry is distinguished and the others are non-ordered). 
The target manifold is stratified as 
$$Y= (Y-f(X)) \; \sqcup \; \bigsqcup N_{\underline{\xi}}^\circ$$
where the union is taken over all {\em non-ordered} mono/multi-singularity types. 
Here, we abuse the notation so that for each semi-ordered/non-ordered type, an ordered multi-singularity type $\underline{\xi}$ is selected as a representative. 
The number of multi-singularity types $\underline{\xi}$ with $\ell(\underline{\xi}) \le l$ are finite, so the strata are finite. 

Note that 
$\dim M_{\underline{\xi}}^\circ = \dim N_{\underline{\xi}}^\circ = l-\ell(\underline{\xi})$ and 
the restriction of $f$ to every stratum
$$f: M_{\underline{\xi}}^\circ \to N_{\underline{\xi}}^\circ$$
is either an isomorphism or a finite-to-one covering (the degree is $\deg_1(\underline{\xi})$). 
If an entry $\xi_j$ with $j\not=1$ is of different type from $\xi_1$, the preimage $f^{-1}(N_{\underline{\xi}}^\circ) \subset X$ contains another stratum $M_{\underline{ \sigma\xi}}^\circ$ with some $\sigma \in Aut(\underline{\xi})$ with $\sigma(1)=j$. 

We denote the embeddings of those loci by 
$$i_{\underline{\xi}}: M_{\underline{\xi}}^\circ \hookrightarrow X, \quad j_{\underline{\xi}}: N_{\underline{\xi}}^\circ \hookrightarrow Y \quad (j_{\underline{\xi}}=f\circ i_{\underline{\xi}})$$ 
and their normal bundles, respectively, by 
$$p_{0}: {\rm N}_0(\underline{\xi})\to M_{\underline{\xi}}^\circ, 
\quad p_{1}: {\rm N}_1(\underline{\xi}) \to N_{\underline{\xi}}^\circ$$
with the zero-sections $s_{0}: M_{\underline{\xi}}^\circ \to {\rm N}_0(\underline{\xi})$ and $s_{1}: N_{\underline{\xi}}^\circ \to {\rm N}_1(\underline{\xi})$. 
Below we often denote the normal bundles by ${\rm N}_0$ and  ${\rm N}_1$ for short if it does not make any confusion. 
Since the structure groups of those vector bundles are reduced to the groups $G(\underline{\xi}/\xi_1)$ and $G(\underline{\xi})$, respectively,  there are classifying maps\footnote{We simplify our notation of classifying maps as noted in the footnote of \S \ref{TP}.} such that those bundles are induced from $X_{\underline{\xi}}$ and $Y_{\underline{\xi}}$, respectively; let $\rho_0$ and $\rho_1$ denote the bundle morphisms. 
In a canonical way, we obtain (the germ at the zero-section of) a locally stable map\footnote{Intuitively, $\varphi_{f,\underline{\xi}}$ is thought of as an algebro-geometric alternative to the restriction map of $f: X\to Y$ to a tubular neighborhood of  $M_{\underline{\xi}}^\circ$ in $X-\rd M_{\underline{\xi}}$. } 
$$\varphi_{f,\underline{\xi}}: {\rm N}_0 \to {\rm N}_1$$
such that the germ of $\varphi_{f,\underline{\xi}}$ at $x_1 \in M_{\underline{\xi}}^\circ$ (zero-section) is isomorphic to the germ $f: (X, x_1) \to (Y, f(x_1))$. 
In fact, $\varphi_{f,\underline{\xi}}$ is induced from the universal singular map $f_{\underline{\xi}}$, i.e., we have the following commutative diagram 
(we set $B_0:=BG(\underline{\xi/\xi_1})$, $B_1:=BG(\underline{\xi})$):  
\[
   \xymatrix{
& & {\rm N}_0  \ar[rd]_{\varphi_{f, \underline{\xi}}} \ar[dd]_{p_0} \ar[rr]^{\rho_0} & & X_{\underline{\xi}} \ar[dd]|{\hole} \ar[rd]^{f_{\underline{\xi}}} &  \\
& & & {\rm N}_1 \ar[dd] \ar[rr]^{\rho_1\qquad}  & &  Y_{\underline{\xi}} \ar[dd] \\
X\ar[rd]_f &  &  M_{\underline{\xi}}^\circ  \ar[ll]_{i_{\underline{\xi}}}\ar[rr]|{\hole} \ar[rd]_f & & B_0\ar[rd]_{\rho_{\underline{\xi}}} & \\
& Y  & &  N_{\underline{\xi}}^\circ \;\; \ar[ll]^{j_{\underline{\xi}}}\ar[rr] & &B_1 } \]

At every $p \in M_{\underline{\xi}}^\circ$, the tangent space at $s_0(p)$ of the total space ${\rm N}_0$ splits into the tangent $T_pM_{\underline{\xi}}^\circ$ and the normal ${\rm N}_0|_p$, thus, especially, we have
$$s_0^*T{\rm N}_0 = i_{\underline{\xi}}^*TX, \quad 
s_0^*\varphi_{f,\underline{\xi}}^*T{\rm N}_1 = j_{\underline{\xi}}^*f^*TY.$$ 

Set $\rd M_{\underline{\xi}}:=M_{\underline{\xi}}-M_{\underline{\xi}}^\circ$ (the union of boundary strata). 
We let $i_{\underline{\xi}}^*$ denote the composition $\CH^{q}(X)_\Q \to \CH^{q}(X-\rd M_{\underline{\xi}})_\Q \to \CH^{q}(M_{\underline{\xi}}^\circ)_\Q$ of the flat pullback (of the open embedding) and the Gysin map (of the closed embedding). 

\begin{lem}\upshape \label{qc1}
It holds that 
$$s_0^*\rho_0^*c(f_{\underline{\xi}})=i_{\underline{\xi}}^*c(f) \in \CH^*(M_{\underline{\xi}}^\circ)_\Q.$$ 
Also $s_I(f_{\underline{\xi}})$ and $s_I(f)$ coincide through the pullback of $\rho_1\circ s_1$ and $j_{\underline{\xi}}$. 
\end{lem}

\proof 
Since ${\rm N}_0$ and ${\rm N}_1$ are induced from the universal bundles $X_{\underline{\xi}}$ and $Y_{\underline{\xi}}$, respectively, we see 
$$\rho_0^*c(\rho_{\underline{\xi}}^*Y_{\underline{\xi}}-X_{\underline{\xi}})
=c(f^*{\rm N}_1-{\rm N}_0) \in \CH^*(M_{\underline{\xi}}^\circ)_\Q$$
(here we abuse the notation $\rho_0$ to mean the classifying map $M_{\underline{\xi}}^\circ \to B_0$ as well). 
Over the total space ${\rm N}_0$, the tangent bundle $T{\rm N}_0$  (resp. $\varphi_{f,\underline{\xi}}^*T{\rm N}_1$) splits into the pullback via $p_0^*$ of two factors, the tangent bundle $TM_{\underline{\xi}}^\circ$ of the base space and the normal bundle ${\rm N}_0$ itself (resp. $f^*{\rm N}_1$). 
Then we may write 
$$c(\varphi_{f,\underline{\xi}}):=c(\varphi_{f,\underline{\xi}}^*T{\rm N}_1-T{\rm N}_0) 
= p_0^*c(f^*{\rm N}_1-{\rm N}_0)$$
in $\CH^*({\rm N}_0)_\Q$ by canceling the factor $TM_{\underline{\xi}}^\circ$. 
That is the same for the tangent bundle $TX_{\underline{\xi}}$ and $f_{\underline{\xi}}^*TY_{\underline{\xi}}$ over $X_{\underline{\xi}}$, and we get 
$$c(f_{\underline{\xi}}):=c(f_{\underline{\xi}}^*TY_{\underline{\xi}} - TX_{\underline{\xi}})
=\pi_0^*c(\rho_{\underline{\xi}}^*Y_{\underline{\xi}}-X_{\underline{\xi}}) $$
in $\CH^*(X_{\underline{\xi}})_\Q$  
(here we may regard the base space $B_0$ as a finite dimensional approximation in the algebraic Borel construction \cite{Totaro}). 
In particular, we have 
$$c(\varphi_{f,\underline{\xi}})=\rho_0^*c(f_{\underline{\xi}}) \in \CH^*({\rm N}_0)_\Q.$$ 
The rest is a straightforward calculation: 
\begin{eqnarray*}
s_0^*\rho_0^*c(f_{\underline{\xi}})
&=&s_0^*c(\varphi_{f,\underline{\xi}})\\
&=&s_0^*c(\varphi_{f,\underline{\xi}}^*T{\rm N}_1 - T{\rm N}_0)\\
&=&c(s_0^*\varphi_{f,\underline{\xi}}^*T{\rm N}_1 - s_0^*T{\rm N}_0)\\
&=&c(j_{\underline{\xi}}^*TY-i_{\underline{\xi}}^*TX)\\
&=&i_{\underline{\xi}}^*c(f^*TY-TX)\\
&=&i_{\underline{\xi}}^*c(f). 
\end{eqnarray*}
Furthermore, LN classes are also, for $\rho_1^*\circ (f_{\underline{\xi}})_*=(\varphi_{f,\underline{\xi}})_*\circ \rho_0^*$. 
 \qed


\

\subsubsection{\bf Restriction equations} 
Let $f: X \to Y$ be as above in \S \ref{stratification} and take the stratifications of $X$ and $Y$ according to multi-singularity types. 

Pick a type $\underline{\eta}=(\eta_1, \cdots, \eta_r)$ such that 
$M_{\underline{\eta}}^\circ\not=\emptyset$ for $f$, and 
put $\ell:=\ell(\underline{\eta})$. 
The closure $M_{\underline{\eta}}$ has the right dimension $m-\ell+\kappa$ and represents 
the source loci class  $m_{\underline{\eta}}(f) \in \CH^{\ell-\kappa}(X)_\Q$ up to the multiplicity  (Proposition \ref{multi-sing}). 
For the simplicity, we may assume that $M_{\underline{\eta}}$ is one of the highest codimensional strata, so it is a smooth closed subscheme $M_{\underline{\eta}}=M_{\underline{\eta}}^\circ$. 
Then we have the so-called {\em restriction equations} \cite{Rimanyi01, FR02} 
$$i_{\underline{\xi}}^* m_{\underline{\eta}}(f)=
\left\{
\begin{array}{cr}
\#\Aut(\underline{\eta}/\eta_1) \cdot c_{\ell-\kappa}({\rm N}_0(\underline{\xi})) & (\underline{\xi}=\underline{\eta}) \\ 
0 & (\underline{\xi}\not=\underline{\eta})
\end{array}
\right. 
\eqno{(r)}
$$
in $\CH^{\ell-\kappa}(M_{\underline{\xi}}^\circ)_\Q\simeq \CH^{\ell-\kappa}({\rm N}_0(\underline{\xi}))_\Q$ 
(this is almost obvious, for $m_{\underline{\eta}}(f)$ is represented by $[M_{\underline{\eta}}] \in \CH_{m-\ell+\kappa}(X)_\Q$). 
Recall that $\varphi_{f,\underline{\xi}}: {\rm N}_0(\underline{\xi}) \to {\rm N}_1(\underline{\xi})$ is deduced from the universal singular map $f_{\underline{\eta}}: X_{\underline{\eta}} \to Y_{\underline{\eta}}$. 

\begin{lem}\upshape \label{qc2}
For any $\underline{\xi}$, it holds that 
$$i_{\underline{\xi}}^*m_{\underline{\eta}}(f)= s_0^*\rho_0^* m_{\underline{\eta}}(f_{\underline{\xi}}) 
\in \CH^{\ell-\kappa}(M_{\underline{\xi}}^\circ)_\Q.$$
\end{lem}

\proof 
If $\underline{\xi}=\underline{\eta}$, the zero section is the stratum of type $\underline{\eta}$ for the map $f_{\underline{\eta}}$, thus 
$m_{\underline{\eta}}(\varphi_{f,\underline{\eta}})={s_{0}}_*(1) \in \CH^{\ell-\kappa}({\rm N}_0(\underline{\eta}))_\Q$, and thus 
$$i_{\underline{\eta}}^*m_{\underline{\eta}}(f)=c_{top}({\rm N}_0(\underline{\eta}))=s_0^*{s_{0}}_*(1)=s_0^*m_{\underline{\eta}}(\varphi_{f,\underline{\eta}})=s_0^*\rho_0^* m_{\underline{\eta}}(f_{\underline{\eta}}).$$ 
If $\underline{\xi} \not=\underline{\eta}$, then the both sides are zero, because $f_{\underline{\xi}}$ does not have $\underline{\eta}$-singular points. 
\qed


\

\subsubsection{\bf Proof of Theorem \ref{main_thm2}} \label{proof_thm2}
Let $f: X \to Y$ be an admissible map of codimension $\kappa (\ge 0)$ in $\Sm$. 
Take the stratification as above, and 
put 
$$\hat{m}_{\underline{\eta}}(f) := \sum R_{J_1}(f)f^*n_{J_2}(f)\cdots f^*n_{J_s}(f)$$
in $\CH^{\ell-\kappa}(X)_\Q$, where the summand runs over all possible partitions 
$\{1, \cdots, r\}= J_1 \sqcup \cdots \sqcup J_s$ with $1 \in J_1$. 
By Lemma \ref{iso} (1) and a remark after Theorem \ref{main_thm1}, we see 
$$f_*m_{\underline{\eta}}(f)=n_{\underline{\eta}}(f)=f_*\hat{m}_{\underline{\eta}}(f) \in \CH^\ell(Y)_\Q.$$ 
In general $f_*$ may not be injective, so this formula is insufficient to imply the theorem. 

First, we claim that for any types $\underline{\xi}$ it holds that 
$$i_{\underline{\xi}}^*{m}_{\underline{\eta}}(f) = i_{\underline{\xi}}^*\hat{m}_{\underline{\eta}}(f)\; \in \CH^{\ell-\kappa}(M_{\underline{\xi}}^\circ)_\Q. \eqno{(*)}$$ 
In fact, applying the above relation to the universal singular map $f_{\underline{\xi}}$, we have 
$$f_{\underline{\xi}*}m_{\underline{\eta}}(f_{\underline{\xi}})=n_{\underline{\eta}}(f_{\underline{\xi}})=
f_{\underline{\xi}*}\hat{m}_{\underline{\eta}}(f_{\underline{\xi}}) \in \CH^\ell(Y_{\underline{\xi}})_\Q.$$ 
Since $f_{\underline{\xi}*}$ is injective by Lemma \ref{injectivity}, we see  
$$m_{\underline{\eta}}(f_{\underline{\xi}})=\hat{m}_{\underline{\eta}}(f_{\underline{\xi}}) \in \CH^{\ell-\kappa}(X_{\underline{\xi}})_\Q.$$
Pulling back the both sides via $\rho_0\circ s_0$, 
$$i_{\underline{\xi}}^*{m}_{\underline{\eta}}(f) 
= s_0^*\rho_0^* m_{\underline{\eta}}(f_{\underline{\xi}}) 
=s_0^*\rho_0^* \hat{m}_{\underline{\eta}}(f_{\underline{\xi}}) 
= i_{\underline{\xi}}^*\hat{m}_{\underline{\eta}}(f)$$
where the first equality follows from Lemma \ref{qc2} and the third one from Lemma \ref{qc1}.  
So the claim $(*)$ is proven. 

Recall the excision property \cite[\S1.8]{Fulton}: for a closed subscheme $W$ of $X$, we have the exact sequence$$\xymatrix{\CH_{q}(W)_\Q \ar[r]^{i_*} & \CH_{q}(X)_\Q \ar[r]^{j^*} & \CH_{q}(X-W)_\Q \ar[r] & 0}$$where $i_*$ is the pushforward of the inclusion and $j^*$ is the flat pullback. 
Initially we take $W$ to be the complement to the open stratum $M_{A_0}^\circ$, and then replace $X$ and $W$, respectively, by $M_{\underline{\xi}}$ and $\rd M_{\underline{\xi}}:= M_{\underline{\xi}} - M_{\underline{\xi}}^\circ$ (the union of boundary strata) according to the adjacency of strata. 

Inductively using the excision, it turns out that the class $\hat{m}_{\underline{\eta}}(f) \in  \CH^{\ell-\kappa}(X)_\Q=\CH_{m-\ell+\kappa}(X)_\Q$ is supported on $\CH_{m-\ell+\kappa}(M_{\underline{\eta}})_\Q$, i.e., for any other type $\underline{\xi}$ than $\underline{\eta}$,  the class vanishes in $\CH_{m-\ell+\kappa}(M_{\underline{\xi}}^\circ)_\Q$. Indeed, if $\hat{m}_{\underline{\eta}}(f)$ is localized to a class $c \in \CH_{m-\ell+\kappa}(M_{\underline{\xi}})_\Q$, then it is pulled back to $j^*c \in \CH_{m-\ell+\kappa}(M_{\underline{\xi}}^\circ)_\Q$ and 
$$s_0^*s_{0*}j^*c=i_{\underline{\xi}}^*\hat{m}_{\underline{\eta}}(f)=0$$ 
by $(r)$ and $(*)$, where $s_0: M_{\underline{\xi}}^\circ \hookrightarrow {\rm N}_0(\underline{\xi})$ is the zero section. Then $j^*c=0$ by our assumption on $f$ that $c_{top}({\rm N}_0(\underline{\xi}))$ is a non-zero divisor 
(indeed, $s_0^*s_{0*}$ is just the multiplication by the top Chern class).  
Hence, by the excision,  $c \in \CH_{m-\ell+\kappa}(\rd M_{\underline{\xi}})_\Q$. 
Since $\rd M_{\underline{\xi}}$ consists of several strata of codimension higher than that of $\underline{\xi}$, the induction proceeds. 

Finally,  $(r)$ and $(*)$ in the case of $\underline{\xi}=\underline{\eta}$ shows that $i_{\underline{\eta}}^*\hat{m}_{\underline{\eta}}(f)$ equals $\#\Aut(\underline{\eta}/\eta_1)\cdot [M_{\underline{\eta}}^\circ]$ in  $\CH_{m-\ell+\kappa}(M_{\underline{\eta}})_\Q$. 
Thus we conclude that $m_{\underline{\eta}}(f)$ and $\hat{m}_{\underline{\eta}}(f)$ coincide. 
 \qed


\begin{rem}\upshape \label{restriction}
{\bf (Restriction method)} 
In case of $\kappa \ge 0$, we have seen that $m_{\underline{\eta}}(f)$ is expressed by a unique polynomial of variables $c_i(f)$ and $s_I(f)$, that is the source multi-singularity Thom polynomial of type $\underline{\eta}$. 
This fact makes it possible to determine the precise form of the Thom polynomial practically.
Namely, the above formula $(r)$ taken for all $\underline{\xi}$ gives an overdetermined system of linear equations of unknown coefficients of the polynomial for $\underline{\eta}$,  and all that remains is to solve it.
That is the strategy of the restriction method, and there have already been done many explicit computations of Thom polynomials, see e.g., \cite{Rimanyi01, FR02, Kaz03, Kaz06, Ohmoto16}. 
This method brought a striking breakthrough in this field, because people had previously thought that to get the Thom polynomial for a mono/multi-singularity type is more or less equivalent to explicitly find a nice desingularization of the singularity locus of that type (in general, the latter problem is hard). 
Indeed, the problem is now understood separately in two aspects, the hidden combinatorics in computation and the existence of universal expression. 
The computation strategy well performs with equivariant geometry, while our proof of the existence of multi-singularity Thom polynomials heavily relies on the existence of desingularizations of geometric subsets ($\Omega$-assignments) that is guaranteed by Hironaka's resolution of singularities. 
\end{rem}


\section{Perspectives}\label{perspective}
\subsection{Further problems}
Finally, we collect some of further research directions and explain each of them briefly. 

\subsubsection{\bf  Thom series and geometric subsets} \label{series} 
For mono-singularity types, Thom polynomials have been considered so that $\kappa=\dim Y - \dim X$ is involved as a {\em parameter} -- those are often called {\em Thom series} \cite{FR12, BercziSzenes, Kaz17}. 
Also for multi-singularity types, it is very natural to expect that there is a universal series containing $\kappa$ as a parameter, say {\em multi-singularity Thom series} (cf. \S \ref{multiple}). Indeed, by our definition, all multi-singularity loci classes of type $\underline{\eta}$ for any $\kappa$ are determined at once by the invariant stratum $\Xi(X; \underline{\eta}) \subset X^{[[n]]}$ which depends on the type of local algebra. Here $Y$ and $\kappa$ are sort of auxiliary data. 
In this sense,  the expected theory of multi-singularity Thom series is tied with the geometry of the invariant stratification of Hilbert schemes.

\subsubsection{\bf Universal Segre clsses}\label{higher}
Let $\underline{\eta}$ be a multi-singularity type of map-germs. 
Suppose that $f: X \to Y$ (in $\Sm$) is a proper locally stable map. 
Set $V:=N_{\underline{\eta}}(f)$ for short and $\iota: V \hookrightarrow Y$ the inclusion. 
If $V$ is smooth, we may consider the Chern classes $c(TV)$ of the tangent bundle and its pushforward into $Y$.  In general, however, the locus $V$ is singular along adjacent strata, and 
then, an alternative to $c(TV)$ is the {\em Chern-Schwartz-MacPherson class} (CSM class) for singular varieties \cite{Mac}, denoted by $c^\SM(V) \in \CH_*(V)$.  
More suitably,  we make use of the {\em Segre-Schwartz-MacPherson class} (SSM class) for $\iota: V \hookrightarrow Y$; 
$$s^{\SM}(V, Y):=\iota^*c(TY)^{-1}\frown c^\SM(V) \; \in \CH_*(V)$$ 
(in case that $V$ is smooth, then $s^{\SM}(V, Y)=c(TV-\iota^*TY)$). 
Consider the SSM class in $\CH^*(Y)$ via the pushforward $\iota_*$ (for $Y$ is smooth), 
the leading term is nothing but the class $[V]$, thus it is expressed by the Thom polynomial of type $\underline{\eta}$ in $Y$. Namely, taking account of the multiplicity (Proposition \ref{multi-sing}), we have 
$$\# \Aut(\underline{\eta})\cdot \iota_*s^{\SM}(V, Y) =n_{\underline{\eta}}(f)+h.o.t. \in \CH^*(Y).$$
In particular, the Euler characteristic of the locus is given by 
$$\chi(V)=\frac{1}{\# \Aut(\underline{\eta})} \int_Y c(TY)\cdot s^{SM}(V, Y).$$

A general theory of universal SSM classes for singularity loci has been proposed in the author's previous work \cite{Ohmoto16, Ohmoto06, Ohmoto08}. 
It was well established for mono-singularity case in \cite{Ohmoto06}. Afterwards, the theory has been featured in variants of {\em Schubert Calculus}, e.g., equivalriant conormal Schubert Calculus, by several authors in relation with the theory of stable envelopes due to Okounkov and others, see e.g., \cite{AMSS17, FRW21}. 

Also for multi-singularity case, it is natural to expect that each higher order term of the SSM class (in $Y$) is universally expressed in terms of Landweber-Novikov classes, which we call {\em target higher Thom polynomial}, provisionally. 
That would possibly be verified in a similar way as in the present paper, that will be discussed somewhere else. In fact, assuming such universal expressions, some higher Thom polynomials have been computed by the restriction method \cite{Ohmoto16, Pallares, NekardaOhmoto2} and applied successfully to the study on vanishing topology of quasi-homogeneous map-germs \cite{Ohmoto16, Pallares}. Also this is related to Chern class of the sheaf of logarithmic vector fields for a free divisor (discriminant) \cite{Liao, Ohmoto16}.

\subsubsection{\bf Generating series}\label{generating}
In \cite{Kaz06} Kazarian presented a generating series for target multi-singularity Thom polynomials, that comes from the combinatorial structure of residual polynomials: 
$$1+\sum_{\underline{\eta}} n_{\underline{\eta}} \cdot 
\frac{t^{\underline{\eta}}}{|\Aut(\underline{\eta})|}
=\exp \left( 
\sum_{\underline{\eta}} f_*(R_{\underline{\eta}}) \cdot 
\frac{t^{\underline{\eta}}}{|\Aut(\underline{\eta})|}
\right).
$$
Also, in a similar way, the author made a conjectural generating series for target higher Thom polynomials (universal SSM classes) in \cite[\S 6.4]{Ohmoto16}. 
It would be very natural to expect that those generating series are directly related to an invariant stratification of $X^{[n]}$ mentioned in \S \ref{series}. 
For example, in \cite{Ohmoto08, CMOSY}, a generating series (zeta function) of the CSM classes of $X^{[[n]]}$ is presented in relation with the symmetric product $S^nX$. 
Its relative version might give some clue for future work.

\subsubsection{\bf Bivariant theory}\label{bi}
As mentioned before, 
we expect that Theorem \ref{main_thm2} on the existence of source multi-singularity Thom polynomials should be true for any codimension $\kappa$ and for arbitrary maps. 
An interesting approach is related to Grothendieck's {\em bivariant theory}. 
From our approach using cohomology operations, an optimistic hope is that there is a suitable {\em bivariant algebraic cobordism theory} which admits a bivariant analogue to Vishik's theorem so that there exists a certain cohomology operation 
$$m_{\underline{\eta}, \mu}^\Omega: \Omega^{0}(X \stackrel{f}{\to} Y) \to \Omega^{\ell-\kappa}(X \stackrel{f}{\to} Y)$$
expressed in terms of Chern classes and Landweber-Novikov classes which satisfies 
that the class $[id: X \to X]$ is sent to the source loci class $m_{\underline{\eta}, \mu}^\Omega(f) \in \Omega^{\ell-\kappa}(X)$. 
That leads to the desired solution for the existence of the source multi-singularity Thom polynomial in the same principle as in the proof of Theorem \ref{main_thm1}. 
This picture may be discussed better in a wider perspecitve of classifying stacks or derived algebraic geometry  (\S \ref{classifying_stack}, \S \ref{motivic}).

\subsubsection{\bf Critical cobordism}\label{critical} 
We have seen in Remark \ref{restriction} that in case of non-negative codimension $\kappa\ge  0$, the restriction method is available for computing certain muti-singularity Thom polynomials. 
However, in case of $\kappa < 0$, Theorem \ref{main_thm1} is not enough for fully supporting the use of this method -- indeed, in this case, the universal singular map $f_{\underline{\xi}}$ is always improper, so Theorem \ref{main_thm1} can not be applied to it. 

This trouble was noticed by Kazarian in \cite{Kaz06} and he introduced {\em localized Thom polynomials} (still its existence is a conjecture). 
Indeed, our true object is the {\em restriction map on the critical scheme of $f$} 
$$f_C: C(f) \to Y$$ 
and it should be enough to assume the properness of $f_C$, not that of $f$. 
Note that $C(f)$ is not always smooth, thus it does not directly yield a cobordism class in $Y$. 
Then it is natural to pick up a specified cobordism class $[\tilde{f}_C: \tilde{C}(f) \to Y]$ by using a canonical desingularization of the first Thom-Boardman singularity, i.e., we take the simplest degeneracy loci class in $\Omega^*(X)$ for the bundle map $df: TX \to f^*TY$ in the sense of \cite{Hudson12} and push it out to $\Omega^*(Y)$. 
We need to characterize this cobordism class and connect it with $m_{\underline{\eta}}(f)$ and $n_{\underline{\eta}}(f)$, and seek for a new variant of algebraic cobordism theory, say it provisionally,  {\em critical cobordism}, for this restricted class of proper maps $\tilde{f}_C: \tilde{C}(f) \to Y$ without starting with $f: X \to Y$.

\subsubsection{\bf Legendre singularities}\label{leg}
A particular case of the expected critical cobordism is an algebraic version of {\em Legendre cobordism theory} (cobordism theory for Legendre maps due to Arnol'd and Vassiliev \cite{Vassiliev}). 
Let $Y$ and $L$ be in $\Sm$ such that $\dim Y=\dim L+1$. 
The projectived cotangent bundle $\pi: \Proj(T^*Y) \to Y$ possesses a canonically defined {\em contact structure} (i.e., maximally non-integrable hyperplane field).  
An immersion $\bar\varphi: L \to \Proj(T^*Y)$ is said to be {\em Legendrian} if it is integrable (i.e., 
the image $d\varphi(T_xN)$ at every point $x \in N$ is included in the contact hyperplane), and 
the composed map $\varphi:=\pi\circ \bar\varphi: L \to Y$ is called a {\em Legendre map} \cite{Arnold}. 
Legendre maps form a special class of maps. 
For instance, if $f: X \to Y$ ($\kappa<0$) is appropriately generic and $\dim 
{\rm coker}\, df_x \le 1$, then the critical locus $L:=C(f)$ is smooth and admits a Legendre immersion $L \to \Proj(T^*Y)$ in a canonical way so that its Legendre map, i.e., the composition with $\pi$, coincides with the restriction $f_C: L \to Y$. 
Legendre singularities has been well explored in Arnol'd's school in the context of the classification of hypersurface singularities (note that it differs from the classification of ICIS). Legendrian Thom polynomials for mono-singularities, that live in the Chow ring of Lagrange Grasmannian,  have been well studied, see \cite{Vassiliev, Kaz03, PragczWeber}. 
So we have enough material to build new theories of {\em algebraic Legendre cobordism} and muti-singularity Legendrian Thom polynomials.

\subsubsection{\bf Multi-singularities in $\A$-classification} \label{A-classfication}
There are many different versions of classification theories of map-germs. Legendre singularities mentioned above is one of them.  Also $\A$-classification is quite important in singularity theory of maps. For locally stable germs, $\A$-classification and $\K$-classification are equivalent \cite{MatherIV},  but for unstable germs, $\A$-classification gets to be much finer. Such unstable singularities appear in families of mappings, e.g., morphisms between total spaces of manifold bundles over the common base, and one seeks for universal expression of the cohomology class represented by the {\em bifurcation locus} of prescribed type defined as a subscheme of the base. 
Thom polynomials in $\A$-classification of mono-germs has been studied in \cite{Sasajima18}, while there is a large room for developing multi-singularity Thom polynomials in $\A$-classification. 
For instance, bifurcation loci of a map from a universal family of curves \cite{Mumford} to another manifold are quite common objects in algebraic geometry (epsecially the theory of Gromov-Witten invariants). 

\subsubsection{\bf Classifying stacks for multi-singularities} \label{classifying_stack}
The set of all (semi-ordered) multi-singularity types of map-germs of codimension $\kappa$ admits a poset-structure via the adjacency relation: $\underline{\xi} < \underline{\tau}$ if  $\underline{\xi}$ is adjacent to $\underline{\tau}$, i.e., $M_{\underline{\xi}}^\circ \subset {\rm Closure}(M_{\underline{\tau}}^\circ)-M_{\underline{\tau}}^\circ$ for a stable unfolding of type $\underline{\xi}$. Let $\kappa>0$ be fixed. 
With respect to this partial order, suppose that we are given a poset $\mathcal{T}$ consiting of several multi-singularity types. A {\em $\mathcal{T}$-map} $f: X \to Y$ is defined to be a proper map in $\Sm$ having only singularities of type in $\mathcal{T}$. 
 We may expect that there is a `universal $\mathcal{T}$-map' 
$$F_{\mathcal{T}}: \mathcal{X}_{\mathcal{T}} \to \mathcal{Y}_{\mathcal{T}}$$
between certain Artin stacks with global stratifications of quotient stacks 
$$\mathcal{X}_{\mathcal{T}}:=\bigsqcup_{\underline{\xi} \in \mathcal{T}}\,  X_{\underline{\xi}}, \;\; \mathcal{Y}_{\mathcal{T}}:=\bigsqcup_{\underline{\xi} \in \mathcal{T}}\,  Y_{\underline{\xi}}$$
such that for any $\mathcal{T}$-map $f: X \to Y$ there exists a morphism (classifying map) $\rho: Y \to \mathcal{Y}_{\mathcal{T}}$ 
whose fiber product with $F_{\mathcal{T}}$ recovers $f$ in a certain sense (perhaps, up to cobordism): 
$$
\xymatrix{
X \ar[r]^{\bar{\rho}} \ar[d]_f & \mathcal{X}_{\mathcal{T}} \ar[d]^{F_{\mathcal{T}}}\\
Y \ar[r]_{\rho} &  \mathcal{Y}_{\mathcal{T}}
}
$$
In differential topology, it has been made up as the {\em Thom-Pontrjagin construction for singular maps} due to Rim\'anyi-Sz\"ucs \cite{RimanyiSzucs}. 
In algebraic geometry,  the counterpart has not yet been found so far -- building blocks, i.e., quotient stacks $X_{\underline{\xi}}$, $Y_{\underline{\xi}}$ and the universal singular map $f_{\underline{\xi}}$, have been obtained, but it is totally unclear how to glue them (perhaps, for this glueing, the Hilbert schemes should take a key role). 
As for Intersection Theory on Artin stacks, see e.g., \cite{Kresch}. 

In this picture, one may assign to the cobordism class of $f: X \to Y$ the pullback $\rho^*[BG_{\underline{\eta}}] \in \CH^\ell(Y)$, that should define our cohomology operation $n_{\underline{\eta}}$. 
It is also interesting to relate the universal map $F_{\mathcal{T}}$ with an expected bivariant theory (cf. \S \ref{bi}). 

\subsubsection{\bf Motivic homotopy theory} \label{motivic}
There is a `bigger' theory of algebraic cobordism due to Morel and Voevodsky, denoted by ${\rm MGL}^{p, q}(-)$ with bi-degree, as an algebraic counterpart to stable homotopy theory in topology. 
If the base field is of characteristic $0$, there is a natural transformation 
$$\Omega^\kappa(Y) \to {\rm MGL}^{2\kappa, \kappa}(Y)$$ 
(for $Y \in \Sm$) and it is known to be isomorphic (M. Levine). 
Beyond this background, it would be very meaningful to 
pursue a direct analogy of Thom's original idea \cite{Thom71} on the Thom-Pontrjagin construction in algebraic geometry side, as Kazarian suggested in \cite{Kaz03}. 
In fact, as the classifying space of algebraic cobordism, the {\em infinite loop space of motivic Thom spectrum} $\Omega_{\rm T}^\infty \Sigma^\kappa_{\rm T} {\rm MGL}_k$ is recently argued in the context of derived algebraic geometry, e.g., see \S 1.4 of \cite{EHKSY}. 
The cohomology of the loop space must contain all the information of universal polynomials associated to multi-singularity types. 
The space should also be related to  $\mathcal{Y}_{\mathcal{T}}$ mentioned in \S \ref{classifying_stack}; it approximates the loop space as an Artin stack presentation.

\subsubsection{\bf  Direct approach via tautological integrals} \label{tautological}
There is still room for thinking about a constructive approach which directly generalizes the proof of double-point formula (cf. \S \ref{multiple}). Consider inductively the class 
$$R_{\underline{\eta}}(f) := m_{\underline{\eta}}(f) - \sum R_{J_1}(f)f^*n_{J_2}(f)\cdots f^*n_{J_s}(f)$$ 
in $\CH^{\ell-\kappa}(X)$, where the sum runs over the partition with $J_1, J_2\not=\emptyset$. 
This class should be understood as the `unnecessary contribution from the smallest diagonal', which is given by the excess intersection of $f^{[[n]]}$ along fibers over the diagonal of the Hilbert-Chow morphism, and the problem is to show that $R_{\underline{\eta}}$ is indeed uniquely written as a polynomial of Chern classes $c_i=c_i(f)$. 
That should be directly related to tautological integrals on Hilbert schemes  \cite{Rennemo, Berczi, BercziSzenes, Berczi23} 
and also non-associative Hilbert schemes \cite{Kaz17}.

\subsubsection{\bf  Real singularities} 
Enumerative theory of multi-singularities for real algebraic, real analytic or $C^\infty$ maps is somewhat mysterious. First we apply the classical Borel-Heafliger theorem (cf. Matszangosz \cite{Akos}). That is, taking fixed point sets of the involution for a complex algebraic variety $X(\C)$, 
we may obtain mod.2 enumerative formulas for multi-singularities of real maps $f: X(\R) \to Y(\R)$ after replacing Chern classes $c_i=c_i(f)$ by the Stiefel-Whitney classes $w_i=w_i(f)$ in $H_*(X(\R); \Z_2)$. 
Also it is quite meaningful to seek for formulas with integer coefficients using Pontryagin classes, but  known results are less, comparing with the case of complex singularities.  
Moreover, as a very new direction, it sounds interesting to seek for `multiple-point formulas' in the Chow-Witt cohomology \cite{Levine19}. 

\subsubsection{\bf  Complex links and quantum invariants} 
For the link of a complex plane curve singularity, the so-called {\em Hilbert scheme invariants} comes up in relation with quantum polynomial invariants of the link in the $3$-sphere $S^3$. There has not been found any higher dimensional analogue so far. That may be related to the Euler product-type expression of the generating series of weighted Euler characteristics of strata of Hilbert schemes, or of the CSM classes of those strata \cite{Ohmoto08, CMOSY} (see \S \ref{generating}). 
A hope is that for a quasi-homogeneous isolated complete intersection singularity (or hypersurface singularity), we may find Hilbert scheme invariants for the complex link via torus-localization of higher Thom polynomials (universal SSM classes for multi-singularities) mentioned in \S \ref{higher}.

\subsection{Concluding remark}
We have developed a general enumerative theory of singularities of maps by solving affirmatively Kazarian's conjecture \cite{Kaz03, Kaz06} in the context of algebraic geometry over an algebraically closed filed of characteristic $0$. We have shown the existence of unique universal polynomials which express the source and the target multi-singularity loci classes for prescribed type (Theorems \ref{main_thm1} and \ref{main_thm2}). Our result also holds suitably for any oriented cohomology theory $A^*$ (Theorem \ref{main_thm}). Main technical ingredients are Hilbert extension maps,  desingularization of geometric subsets ($\Omega$-assignments), cohomology operations on algebraic cobordism $\Omega^*$, and equivariant geometry with symmetries of singularities. 
We remark that in Theorem \ref{main_thm2} about the source multi-singularity loci classes, we only consider a certainly nice case with $\kappa \ge 0$, thus the full scope of the conjecture, including the setup over a ground field of positive characteristic, still remains open. 

Thom polynomial theory connects local and global natures of singularities of maps. 
Local classification theory of singularities of maps is essentially the matter of commutative algebra and 
has been thoroughly studied for several decades, while the global homotopical nature in algebraic geometry has been developed much later -- it was quite recent that key materials have become available, such as algebraic Borel construction, algebraic version of Thom spectra, algebraic coboridisms and algebraic cohomology operations, and so on. 
Now, on the basis of these recent achievements in algebraic geometry, we have viewed a new horizon in enumerative geometry of singularities.


\subsection*{Acknowledgement} 
The author would like to thank organizers and participants of the conference ``Characteristic classes for singular spaces" (Kiel, Sept. 2023), where he had an opportunity to give a series of lectures on the topic of this paper.  
This work was partly supported by JSPS KAKENHI Grant Numbers JP17H02838 and JP23H01075.


\end{document}